\documentclass[10pt,reqno]{amsart}

\usepackage{amsmath}
\usepackage{amssymb}
\usepackage{amsfonts}
\usepackage{epsfig}

\newtheorem{proposition}{Proposition}[section]
\newtheorem{theorem}{Theorem}[section]

\newtheorem{corollary}{Corollary}[section]
\newtheorem{remark}{Remark}[section]

\numberwithin{equation}{section}

\begin{document}

\qquad\qquad\qquad\qquad\qquad\qquad{Dedicated to the memory of V.N.Diesperov}
\title[Correct Statement, Analysis and Numerical Solution \ldots]{Correct Statement, Analysis and Numerical
Solution of Singular Nonlinear Problems for Self-Similar Solutions to the Boundary Layer Equations
with Zero Pressure Gradient}
\author{N.B.Konyukhova}
\address{Dorodnicyn Computing Center of Federal Research Center "Computer Science and Control",
Russian Academy of Sciences, ul. Vavilova 40, Moscow, 119333 Russia}
\email{nadja@ccas.ru}
\author{S.V.Kurochkin}
\address{Dorodnicyn Computing Center of Federal Research Center "Computer Science and Control" \,,
Russian Academy of Sciences, ul. Vavilova 40, Moscow, 119333 Russia}
\email{kuroch@ccas.ru}
\author{M.B.Soloviev}
\address{Dorodnicyn Computing Center of Federal Research Center "Computer Science and Control"\,,
Russian Academy of Sciences, ul. Vavilova 40, Moscow, 119333 Russia}
\email{solmb@mail.ru}
\begin{abstract} For the problems indicated in the title,
a further development of a new approach (different from those applied before) is given.
A basic problem under consideration arises in viscous incompressible fluid dynamics
and describes self-similar solutions to the boundary layer equation for a stream function
with zero pressure gradient (connected with the plane-parallel laminar flow in a mixing
layer). Some previous results concerning singular nonlinear Cauchy problems, smooth stable initial manifolds,
and parametric exponential Lyapunov series are used to state correctly and analyze
the singular "initial-boundary-value" \,problem for a third-order nonlinear ordinary differential
equation defined on the entire real axis. The detailed analysis of this singular nonlinear problem leads,
in particular, to efficient methods for solving it approximately and gives a possibility
to obtain numerically the particle trajectories in the plane of flow. Some results of the numerical
experiments are displayed and their physical interpretation is discussed. A connection of this basic problem
with some known physical and mathematical problems, arising for self-similar solutions to the boundary layer
equations with zero pressure gradient, is described, namely the "flooded jet"\,, the plane "semi-jet"\, and  the "near-wall jet"\, problems are considered which are of interest by themselves.
\end{abstract}
\maketitle
\tableofcontents
\section{Introduction}
For correct statement and study of the problems indicated in the paper title, we use a new approach of \cite{dks_07}, \cite{kss_09} (different from that applied in \cite{dies1_86}, \cite{dies2_86}) and give its further development. Along with new results, we present the main results of \cite{dks_07}, \cite{kss_09} in a revised, extended and more correct form (certain misprints and inaccuracies admitted in \cite{dks_07}, \cite{kss_09} are also corrected here). A more detailed comparative analysis of the methods and results from \cite{dks_07}, \cite{kss_09} and this paper with the previous ones from  \cite{dies1_86}, \cite{dies2_86} is also given.

The basic problem under consideration arises in incompressible viscous fluid dynamics and describes
self-similar flows in mixing layers. Being formulated in \cite{dies1_86} and \cite{dies2_86}, it has the form
\begin{equation}\label{d1_in}
\Phi^{\prime\prime\prime} + \Phi\Phi^{\prime\prime} - [(m-1)/m]\,(\Phi ^\prime)^{2}=0,
\quad  - \infty < \tau < \infty,
\end{equation}
\begin{equation}\label{con1_in}
\lim_{\tau\to -\infty}{\Phi^\prime(\tau)}=0,
 \end{equation}
\begin{equation}\label{con2_in}
\Phi(0)=0,
\end{equation}
\begin{equation}\label{con3_in}
\lim_{\tau\to\infty}{\left(\Phi(\tau )/ \tau^{m}\right)}= b,
\end{equation}
where \,$m>0$ and \,$b>0$ are given parameters.

In Subsection 2.1, the mathematical description of the original physical model is given. It concerns a flow in a mixing layer, resulting from the interaction of two unbounded layers of viscous incompressible fluid of which the upper one moves (with a power dependence on hight of the horizontal velocity component) and the lower one is at rest. For self-similar flows in the form of \cite{dies1_86}, \cite{dies2_86}, it leads at first sight to the singular nonlinear problem \eqref{d1_in}--\eqref{con3_in}.
But, for mathematically correct statement of this problem defined on the whole real line, a concept of admissible limit conditions at infinity for systems of nonlinear autonomous ODEs (see \cite{k1_94}--\cite{kony_01} and references therein) should be take into account.

As a result of the study in \cite{kony_96} (see also \cite{dks_07} and Subsection 2.2 here), condition \eqref{con1_in} should be replaced by more accurate limit condition with a parameter \,$a>0$:
\begin{equation}\label{ccon1_in}
\lim_{\tau\to -\infty}{\exp{(-\varepsilon \tau)}\{\Phi(\tau) + a, \Phi^\prime (\tau), \Phi^{\prime\prime}(\tau )\}}=
\{0,0,0 \} \quad \forall \varepsilon: 0 < \varepsilon <a.
\end{equation}
This condition corresponds to the property of a solution to tend to the stationary point \,$(-a,0,0)$ in the phase space of ODE \eqref{d1_in}. For any fixed \,$a>0$, this point is a pseudohyperbolic equilibrium point with a one-dimensional stable separatrix. Condition \eqref{ccon1_in}, for the solutions of ODE \eqref{d1_in}, is equivalent in a finite point \,$\tau=-T$, $T\ge 1$, to two nonlinear relations which specify a stable saddle separatrix. Thus, provided that \,$-\infty<\tau\le 0$, two-point BVP \eqref{d1_in}, \eqref{ccon1_in}, \eqref{con2_in} with the parameter \,$a>0$ is defined. The parameter \,$a=a(b)$ then can be found from \eqref{con3_in} if such behavior of the BVP solutions is valid (as follows below, it is correct when \,$m:\,1/2<m<\infty$).

In \cite{dies1_86}, \cite{ dies2_86}, in the study of the singular nonlinear problem \eqref{d1_in}--\eqref{con3_in}, the order reduction methods are used that lead to rather a complicated analysis on "the Poincar\'e sphere"\, of a two-dimensional nonlinear dynamical system involving some nonphysical variables, and the return to the initial variables makes rather a difficult procedure (also for some fixed physically interesting values of \,$m$, analogous methods for problem \eqref{d1_in}--\eqref{con3_in}  were applied earlier in \cite{dies_84} and \cite{dies_85}). Brief discussions of this original approach to the singular problem  \eqref{d1_in}--\eqref{con3_in} are given in \cite{dks_07} and here in Appendix B.
 
The different approach presented in \cite{dks_07}, \cite{kss_09} and in this paper takes advantage of the results on
singular Cauchy problems (CPs), smooth stable initial manifolds (SIMs), and parametric exponential Lyapunov series (see the classical monograph \cite{lyapunov}), and allows to obtain a more accurate statement of the problem in terms of the initial physical variables: the problem is now split into (1) a singular two-point BVP (with an unknown parameter) on non-positive real semi-axis and (2) a CP on the positive semi-axis that gives the continuation of the solution of the BVP.
The constraints on the self-similarity parameter are then formulated that guarantee the existence and uniqueness of the solution to the basic IBVP under consideration. The two-sided estimates for the solution are given and its properties are investigated along with the properties of other (regular and singular) solutions of the nonlinear ODE involved for various self-similarity parameter values. Corresponding numerical methods and results of computations are presented (unlike \cite{dks_07} and this paper, no numerical results are given in \cite{dies1_86}, \cite{dies2_86}).

As an addition to the results of \cite{dks_07}--\cite{dies2_86}, the statement is given here for the first time of a  basic problem for the two-dimensional stationary boundary layer equations (with zero pressure gradient) considered in the whole space whose solution   in a class of self-similar functions leads to problem \eqref{d1_in}--\eqref{con3_in}.
Along with numerical simulations of stream function (as a function of self-similar variable), the  particle trajectories in the plane of flow are  obtained numerically which were not presented before by other authors. Of separate mathematical and physical interest is an auxiliary two-point BVP on the non-negative real semi-axis having, in particular, exact solutions for certain self-similarity parameter values; numerically obtained flow patterns for the corresponding problems were not also presented before.

Moreover, we give a more detailed and correct analysis of mathematically interesting accompanying auxiliary singular nonlinear problems arising due to the approach of \cite{dies1_86}, \cite{dies2_86} (see Appendix B).

Concerning the boundary layer theory and the problems of fluid and gas dynamics discussed in the present paper, monographs \cite{shlich}--\cite{ol_sam} are used. Certain earlier publications \cite{kony_96}--\cite{kony_01}, \cite{ks_2001}--\cite{ks_2008} either contain some preliminary results for the above basic problem or use it as an example to the correct statement and solution of singular BVPs for autonomous systems of nonlinear ODEs.

The papers \cite{dies1_86}, \cite{dies2_86}, \cite{dies_84}, \cite{dies_85}  by a known specialist in fluid and gas dynamics
give a considerable staff of (evidently) new formulas, results and conclusions. Yet we believe that both the original model and
the concomitant basic and auxiliary singular nonlinear ODE problems are worth a more detailed mathematical and numerical treatment. We think that our approach to these tasks is more simple and that it makes possible to give the full answers to good many questions, which is rather a rare opportunity in case of singular nonlinear problems. On the other hand,
the problems under study give us an opportunity to demonstrate that self-similar solutions to hydrodynamics problems
are themselves difficult enough and demand thorough analysis.

\section{Statement of the Basic Singular Nonlinear Problem and Preliminary Propositions and Remarks}
\subsection{Mathematical Description of the Original Physical Model}
We discuss a mathematical model of a flow in a mixing layer, resulting from the interaction
of two unbounded layers of viscous incompressible fluid of which the upper one moves
(with a power dependence on hight of the horizontal velocity component) and the lower one
is at rest. For a model description, the steady-state
boundary layer equations for a plane-parallel laminar flow with
zero pressure gradient are used:
\begin{equation}\label{pran_eq}
u\frac{\partial u}{\partial x}+v\frac{\partial u}{\partial y}= \nu
\frac{\partial ^2 u}{\partial y^2}, \quad x>0, \quad y \in
\mathbb R,
\end{equation}
\begin{equation}\label{eq_inc}
\frac{\partial u}{\partial x}+\frac{\partial v}{\partial y}=0, \qquad x>0,
\qquad y \in \mathbb R
\end{equation}
(see, e.g., \cite{shlich}, ch. IX, and  \cite{ol_sam}, ch. I). Here, \eqref{pran_eq}
is the Prandtl equation, and \eqref{eq_inc} is the equation of continuity
(incompressibility); the \,$x$ axis is aligned with
the flow and coincides with the free streamline, \,$u$ and \,$v$ are
the velocity components along and perpendicular to the flow
direction correspondingly, and \,$\nu$ is  kinematic viscosity
parameter (in dimensionless variables, we can set \,$\nu=1$; see below Remark \ref{r_2}).

Due to the physical interpretation, the form of self-similar solutions
stated below and the limit relation \eqref{con3_in}, the following limit boundary conditions
must be fulfilled \,$\forall x>0$:
\begin{equation}\label{con_u}
\lim_{y\to -\infty}{u(x,y)}=0,
\end{equation}
\begin{equation}\label{con_v}
v(x,0)=0,
\end{equation}
\begin{equation}\label{con_uv}
\lim_{y\to \infty}{\left[u(x,y)/y^{m-1}\right]=U_0}, \qquad \lim_{y\to \infty}{v(x,y)}=0,
\end{equation}
where \,$m$ and \,$U_0$ are given parameters, \,$m>0$, \,$U_0>0$.

A stream function $\psi(x,y)$ is introduced to satisfy the
equation of continuity \eqref{eq_inc}. Then, taking into account that the
\,$x$ axis coincides with the free streamline, the following
relations hold:
\begin{equation}\label{uv_psi}
u(x,y)=\frac{\partial \psi}{\partial y}(x,y),\quad
v(x,y)= -\frac{ \partial \psi}{\partial x}(x,y), \quad
\psi(x,0)=0 \quad \forall x >0.
\end{equation}

For \,$\psi (x,y)$, we obtain the following problem:
\begin{equation}\label{psi_eq}
\frac{\partial \psi}{\partial y} \, \frac{\partial ^2 \psi}{ \partial
x \partial y} - \frac{\partial \psi}{\partial x} \,
\frac{\partial ^2 \psi}{\partial y^2} = \nu \, \frac{\partial ^3
\psi}{\partial y^3}, \qquad x >0, \qquad y \in \mathbb R,
\end{equation}
\begin{equation}\label{con1_psi}
\lim_{y \to -\infty}{\frac{\partial{\psi(x,y)}}{\partial y}}=0
\qquad \forall x>0,
\end{equation}
\begin{equation}\label{con2_psi}
\psi(x,0)=0 \qquad \forall x >0,
\end{equation}
\begin{equation}\label{con3_psi}
\lim_{y \to \infty}{\left (\frac{\partial{\psi(x,y)}}{\partial y}\Big/y^{m-1}\right )}=U_0,
\qquad \lim_{y \to \infty}{\frac{\partial{\psi(x,y)}}{\partial x}}=0   \qquad \forall x>0.
\end{equation}

The solutions to the problem \eqref{psi_eq}--\eqref{con3_psi} are sought in the class of the self-similar
functions of the form suggested in \cite{dies1_86}, \cite{dies2_86}:
\begin{equation}\label{ss_psi}
\psi (x,y)=\omega^{-1/2} x^{\nu \omega} \Phi (\tau),
\end{equation}
\begin{equation}\label{ss_var}
\tau=\omega^{1/2} y/x^{1/(m+1)}, \quad \omega >0, \quad m>0,
\quad \nu\omega =m/(m+1).
\end{equation}

The self-similar variable \,$\tau$ depends on the
self-similarity parameter \,$m$, and the following relations are
valid (a prime denotes a derivation on \,$\tau$):
\begin{equation}\label{u_ss}
u(x,y)=\frac{\partial \psi}{\partial y}(x,y) =
x^{(m-1)/(m+1)} \Phi^ \prime(\tau),
\end{equation}
\begin{equation}\label{v_ss}
v(x,y)=-\frac{\partial \psi}{\partial x}(x,y) =
\sqrt{\frac{\nu}{m(m+1)}}\, x^{-1/(m+1)}[\tau \Phi^\prime
(\tau)-m \Phi (\tau)].
\end{equation}

For the unknown function \,$\Phi(\tau)$, we obtain the singular nonlinear
problem, with the self-similarity parameter \,$m>0$:
\begin{equation}\label{phi_eq1}
\Phi^{\prime\prime\prime}+\Phi \Phi^{\prime\prime}-[(m-1)/m](\Phi ^\prime)^{2}=0,
\quad  - \infty<\tau<\infty,
\end{equation}
\begin{equation}\label{phi_con1}
\lim_{\tau\to -\infty} \Phi ^\prime(\tau)=0,
\end{equation}
\begin{equation}\label{phi_con2}
\Phi(0)=0,
\end{equation}
\begin{equation}\label{phi_con3}
\begin{array}{c}
\lim_{\tau\to \infty}{\left[\Phi^\prime(\tau)/\tau^{m-1}\right]}=U_0/\left[m/(\nu(m+1))\right]^{(m-1)/2},
\\
\\
\lim_{\tau\to \infty} \left[\tau\Phi^\prime(\tau)-m \Phi(\tau)\right]=0.
\end{array}
\end{equation}
Conditions \eqref{phi_con3}  imply the limit relation \eqref{con3_in}, where we obtain
\begin{equation}\label{b_U0}
b=(U_0/m)/\left[m/(\nu(m+1))\right]^{(m-1)/2}.
\end{equation}

Finally we have formally the singular nonlinear problem as in \cite{dies1_86}, \cite{dies2_86}, with the
parameters \,$m>0$, and \,$b>0$:
\begin{equation}\label{eqf_phi}
\Phi^{\prime\prime\prime}+\Phi \Phi^{\prime\prime}-[(m-1)/m]\,(\Phi ^\prime)^{2}=0,
\quad  - \infty < \tau < \infty,
\end{equation}
\begin{equation}\label{con1f_phi}
\lim_{\tau \to -\infty} \Phi ^\prime(\tau)=0,
\end{equation}
\begin{equation}\label{con2f_phi}
\Phi(0)=0,
\end{equation}
 \begin{equation}\label{con3f_phi}
 \lim_{\tau\to \infty} \left(\Phi (\tau )/ \tau ^{m} \right)= b>0.
\end{equation}

Let us remark preliminary: for \,$m=1$, analogous (more simple) problem is discussed
in \cite{shlich}, pp.180--181 (see also references therein) as a problem for laminar layer on an interface between
two flows (known as a plane "semi-jet"\, problem). In Subsection 4.3.2 of Section 4, we give a comparison
of our numerical results with ones from \cite{shlich} and demonstrate in addition a picture
of flows in the plane of \,$\{x,y\}$ (this example, for \,$m=1$, doesn't discuss in \cite{dks_07}, \cite{kss_09} and
the initial editions of the Schlichting monograph don't contain this problem).

\begin{remark}\label{r_1}(about a setting of original problem).
A setting of the original problems in form \eqref{pran_eq}--\eqref{con_uv}, for the velocity components,
and as a consequence in form \eqref{psi_eq}--\eqref{con3_psi} for the stream function,  has not been
considered in our previous publications and might has not been discussed in the papers by other authors as well
(at least in {\rm\cite{dies1_86}, \cite{dies2_86}} such statements are also missed).
\footnote{The need for statements of the original problems has been noted by V.N.Samokhin \,(XX Intern. Conf. "Mathematics. Economics. Education." -- 2012, Rostov-na-Donu, Russia).}
It is possible that there are some questions: {\rm 1)} Do there exist different solutions
to the original singular nonlinear problem \eqref{pran_eq}--\eqref{con_uv} except
the self-similar ones (the existence of the latter under certain constraints on the
self-similarity parameter $m$ is known from the publications cited above and
is discussed in more detail below), or the solutions' self-similarity results from
the absence of conditions at $x=0$ and as $x\to\infty$ ? {\rm 2)} Is it correct to
state the following  problem (analogous to similar ones in {\rm\cite{ol_sam}}): provided
a velocity profile $u(0,y)$ is given at $x=0$, such that \eqref{con_u} and \eqref{con_uv} hold,
and, e.g., $v(0,y)\equiv 0$, does there exist a solution to the initial problem
\eqref{pran_eq}--\eqref{con_uv} with such addition,  and, if it does, will the fluid flow
corresponding to the solution tend to some self-similar mode at large $x$ values,
and under which conditions?

We emphasize again that this paper deals only with problem
\eqref{eqf_phi}--\eqref{con3f_phi} for self-similar solutions and discusses some
consequences from it. As far as the correctness of the above-discussed model
with flow pattern in upper layer of form \eqref{con_uv} is concerned, such a model seems
correct to us at least because for fixed $m : 1/2< m<\infty$ and $b>0$ problem
\eqref{eqf_phi}--\eqref{con3f_phi} is uniquely solvable, as it will be seen from the analysis below.
For a more detailed discussion of a physical meaning of the problem under consideration, see
{\rm\cite{dies1_86}, \cite{dies2_86}, \cite{dies_84}, \cite{dies_85}}
and references therein (see also {\rm\cite{shlich}} and bibliography therein
for the discussion of analogous problem with $m=1$).
\end{remark}
\begin{remark}\label{r_2}(on the dimensionless variables). Any approach how to do the variables,
e.g., in Eq.\eqref{psi_eq} without dimensions is not discussed in {\rm\cite{dies1_86}, \cite{dies2_86}},
and parameter $\nu$ is kept, though in dimensionless variables it may be set to unity
(in general, parameter $\nu$ is kept in other publications and monographs on the subject).
As all the variables and magnitudes are considered dimensionless in the analysis below,
let us consider one particular technique of making equation {\rm\eqref{psi_eq}} dimensionless.
Following a conventional approach used in boundary layer theory for the variables $x$ and $y$,
we set $x=L\widetilde x$, where $L$ is a characteristic length along which the flow is considered,
$y=(L/\sqrt{{\rm Re}})\widetilde y$, where ${\rm Re}$ is Reynolds number, ${\rm Re} = U_\infty L/\nu$
(for the characteristic velocity $U_\infty$, we use magnitude $U_\infty = U_0 Y^{m-1}$,
where $Y$ is a characteristic height along which the upper flow is considered);
$\psi = (U_\infty L/\sqrt{{\rm Re}})\widetilde \psi$. Particularly, if we take $Y=L/\sqrt{{\rm Re}}$,
then we get ${\rm Re}=(U_0 L^m /\nu)^{2/(m+1)}$. Then, for dimensionless magnitude
$\widetilde\psi(\widetilde x, \widetilde y)$, we obtain a problem of the same form as
in \eqref{psi_eq}--\eqref{con3_psi}, where parameters $\nu$ and $U_0$ are set to unity.
\end{remark}
\begin{remark}\label{r_3}(to a flow description in the initial variables).
Recall that self-similar variable $\tau=\tau(x,y, \nu, m)$ ($\nu>0$,  $m>0$) has the form:
\begin{equation}\label{tau}
\tau(x,y)=\sqrt{\frac{m}{\nu(m+1)}}\, y/x^{1/(m+1)}, \quad  x>0, \quad  y\in \mathbb{R}.
\end{equation}
For fixed \,$m$ and \,$b$, let \,$\Phi_m(\tau, b)$ be a solution to the singular nonlinear
problem \eqref{eqf_phi}--\eqref{con3f_phi}. Then to study the steady state motion of
fluid particles in the plane \,$\{x,y\}$, we can use the following
nonlinear CPs with the parameters (due to singularities when
$x\to +0$ and other evident enough reasons, the numerical treatment of
these CPs is rather difficult):
\begin{equation}\label{x_t}
\frac{dx}{dt}=u(x,y)=x^{(m-1)/(m+1)}\Phi_m^\prime(\tau,b),
\end{equation}
\begin{equation}\label{y_t}
\frac{dy}{dt}=v(x,y)=\sqrt{\frac{\nu}{m(m+1)}}\, x^{-1/(m+1)}[\tau \Phi_m^\prime (\tau,b)-m
\Phi_m (\tau,b)],
\end{equation}
\begin{equation}\label{x0_y0}
x(t_0)=x_0, \quad y(t_0)=y_0 \quad  (x_0>0,  \quad y_0\in\mathbb{R});
\end{equation}
these CPs are equivalent to the CPs with the parameters \,$x_0$, \,$y_0$ of the following form:
\begin{equation}\label{y_x}
\frac{dy}{dx}=\frac{v(x,y)}{u(x,y)}=\sqrt{\frac{\nu}{m(m+1)}}\, x^{-m/(m+1)}
\left[\tau-m\Phi_m(\tau,b)/\Phi_m^\prime(\tau,b)\right],
\end{equation}
\begin{equation}\label{y0_x0}
y(x_0)=y_0  \quad (x,x_0>0, \quad y,y_0\in \mathbb{R}),
\end{equation}
where various values of \,$x_0$ and \,$y_0$  must be considered.

A different approach is to define the level lines of the stream
function
\begin{equation}\label{psi_xy}
\psi(x,y)= \sqrt{\frac{\nu(m+1)}{m}}\, x^{m/(m+1)} \Phi_m
(\tau,b), \quad x>0, \quad y\in \mathbb{R},
\end{equation}
which is in some sense equivalent to solving the above nonlinear
CPs for ODEs but more convenient.
\end{remark}
\begin{remark}\label{r_4} (on the flows in the upper layers). For the
horizontal velocity component \,$u(x,y)$ in the upper layers, we
derive from \eqref{u_ss}, \eqref{phi_con3} and \eqref{con3f_phi}, for large \,$y$ and
\,$\forall x>0$:
\begin{equation}\label{u_up}
u(x,y) \sim U _0(m, b, \nu) y^{m-1}, \qquad y \gg 1,
\end{equation}
\begin{equation}\label{U0}
U_0(m, b, \nu) = mb \Bigl\{ m/[(m+1)\nu]\Bigr\}^{(m-1)/2}.
\end{equation}
Thus, the specification of \,$b>0$ in condition \eqref{con3f_phi} is equivalent
to the specification of \,$U_0$ in \eqref{u_up} which describes the
\,$y$-dependence of the horizontal velocity component in the upper
flow, for large \,$y$,  and there exist three different cases corresponding to
the values \,$m<1$, \,$m=1$, or $m>1$.
\end{remark}

According to \cite{dies1_86}, \cite{dies2_86}, \cite{dies_84}, \cite{dies_85} (see also references therein), the third-order nonlinear autonomous ODE \eqref{eqf_phi} and some similar ones have been extensively studied, for different fixed values of the
self-similarity parameter \,$m$. In our view, there exist more incompressible fluid dynamics problems connected with this and other ODEs, which need a more rigorous mathematical statement as singular nonlinear BVPs (or, in a sense, IBVPs) and require a more comprehensive and rigorous mathematical analysis. Some simple examples are indicated below, in addition to the main problem under consideration.

For the new approach suggested in \cite{dks_07}, \cite{kss_09} and developed in the paper presented, the results concerning singular nonlinear CPs, smooth SIMs of solutions, parametric exponential Lyapunov series, and asymptotic expansions are used. Analogous approach may be developed for study of different singular nonlinear BVPs arising in incompressible fluid dynamics including the non-autonomous ODE systems.
\subsection{Correct Mathematical Statement of the Singular Nonlinear IBVP and  Axillary Singular Nonlinear BVP
with a Parameter}
Given an ODE system with a singular boundary point (finite or at infinity), then the limit conditions at this
point must be formulated for all the solution components defining locally a singular CP. When there exists
a family of solutions to the singular CP, the dimension of such set and the relations generated by the values
of the solutions in the phase space are important, in particular for the accurate transfer of the
boundary conditions from the singular point to a close nonsingular one.

For autonomous nonlinear ODE systems with (pseudo)hyperbolic equilibrium points, the above problems are investigated
in detail (see \cite{k1_94}--\cite{kony_01}), and the corresponding results are used here to correct the previous
formulation of the singular nonlinear problem in form \eqref{eqf_phi}--\eqref{con3f_phi}. In particular, in addition
to the condition \eqref{con1f_phi}, we have to impose conditions for \,$\Phi(\tau)$ and \,$\Phi^{\prime\prime}(\tau)$ as
\,$\tau\to -\infty$. To this end, first we must study the stationary points of the nonlinear autonomous ODE \eqref{eqf_phi}.
In the phase space \,$\mathbb R^3$ of the variables \,$(z_1, \,z_2, \, z_3)=(\Phi, \,\Phi^\prime, \,\Phi^{\prime\prime})$,
Eq.\eqref{eqf_phi} has an infinite set of stationary points (equilibria):
\begin{equation}\label{s_p}
(z_1, \,z_2, \,z_3)_s (a)=(\Phi, \,\Phi^\prime, \,\Phi^{\prime\prime})_s (a)=
(-a, \,0, \,0), \quad a\in\mathbb R.
\end{equation}
In terms of \,$z =(z_1,\,z_2,\,z_3)^T$,  where
\begin{equation}\label{z}
z_1(\tau)=\Phi(\tau), \quad z_2(\tau)=\Phi^\prime(\tau), \quad z_3(\tau)=
\Phi^{\prime\prime}(\tau),
\end{equation}
we obtain the ODE system:
\begin{equation}\label{z_eq}
z^\prime = Q(z), \quad  -\infty<\tau<\infty,
\end{equation}
\begin{equation}\label{Q}
Q(z) = \left(
\begin{array}{c}
z_2 \\
z_3 \\
z_2^2 (m-1)/m - z_1 z_3
\end{array}
\right).
\end{equation}
The Jacobian matrix, for \eqref{Q}, taken at the fixed stationary
point \eqref{s_p} has the form
\begin{equation}\label{J_Q}
\frac{\partial Q}{\partial z} \left(z_s(a)\right) = \left (
\begin{array}{ccc}
0 & 1 & 0 \\
0 & 0 & 1 \\
0 & 0 & a
\end{array}
\right).
\end{equation}
Then \,$\forall a>0$ stationary point \eqref{s_p} (for system \eqref{z_eq}, \eqref{Q}, considered on
\,$\mathbb{R}_-$), is a pseudo-hyperbolic saddle with a one-dimensional stable separatrix (or, by
different definition, a saddle-node of the \,$(-\varepsilon,1)$-type, where \,$\varepsilon$ is an
arbitrary number in the interval $(0,a)$).

This means that condition \eqref{con1f_phi} has to be replaced by the more precise limit condition
with the parameter \,$a>0$:
\begin{equation}\label{ccon1f_phi}
\lim_{\tau\to-\infty}{\exp{(-\varepsilon \tau)}\left\{\Phi(\tau)+a, \,\Phi^\prime(\tau),
\,\Phi^{\prime\prime}(\tau)\right\}}=\{0, \,0, \,0 \} \quad \forall\varepsilon: \,0<\varepsilon<a.
\end{equation}
The delicacy of setting boundary condition \eqref{ccon1f_phi}, with the additional parameter
\,$\varepsilon: 0<\varepsilon<a$,  is associated with the fact that Jacobian matrix \eqref{J_Q} has a
second-order Jordan block corresponding to a zero eigenvalue. Anyway, whether there exists an analytical
SIM for the solutions to Eq.\eqref{eqf_phi} in the neighborhood of point \eqref{s_p} and what its
dimension is can be answered only when the limit conditions are set in the form of \eqref{ccon1f_phi}.

The problem \eqref{eqf_phi}, \eqref{ccon1f_phi} is treated as a singular nonlinear CP. The following
result succeeds from \cite{lyapunov}, section 23.
\begin{proposition}\label{p_1}(the Lyapunov series).
For fixed \,$a>0$ and \,$m \ne 0$, the singular
nonlinear CP  \eqref{eqf_phi}, \eqref{ccon1f_phi} has a one-parameter family of solutions
\,$\Phi _m (\tau, a, d)$. These solutions can be represented in the form of the exponential Lyapunov series
\begin{equation}\label{l_ser}
\Phi_m (\tau, a, d)= -a + d\exp{(a\tau)} + \sum_{l=2}^\infty{h_l\,d^l
\exp{(la\tau)}}, \quad \tau\le \widetilde\tau, \quad \widetilde\tau \in \mathbb R,
\end{equation}
where \,$d$ is a parameter, \,$|d\exp{(a\widetilde\tau)}|$ is small, and
the coefficients \,$h_l$ are independent of \,$d$  \,($l\ge 1$, \,$h_1\doteq 1$):
\begin{equation}\label{hl}
h_l = \left [\sum _{k=1}^{l-1}{k\left(\frac{(m-1)(l-k)}{m} - k \right) h_kh_{l-k}}\right]\Bigr/[al^2(l-1)],
\qquad l=2,3,\ldots \,.
\end{equation}
In particular, it follows from \eqref{hl} that
\begin{equation}\label{h2_3}
h_2  = -1/(4am), \qquad  h_3 = (m+4)/(72a^2m^2), \quad \ldots\, .
\end{equation}
\end{proposition}
\begin{remark}\label{r_5}(to the limit case \,$m\to\infty$).
For \,$m\to \infty$ (i.e., setting \,$m = +\infty$ in \eqref{eqf_phi} and  \eqref{l_ser}--\eqref{h2_3}),
we obtain that singular CP \eqref{eqf_phi}, \eqref{ccon1f_phi} has the two-parameter set of the exact
solutions \,$\Phi _\infty (\tau,a,d)$ existing on the entire real line:
\begin{equation}\label{phiInf}
\Phi _\infty (\tau, a, d) = - a + d\exp{(a\tau)}, \quad \tau\in \mathbb{R},
\end{equation}
where \,$a$ and \,$d$ are the parameters, \,$a>0$, \,$d\in \mathbb{R}$.
\end{remark}

Taking into account the results of \cite{k1_94}--\cite{kony_96}, we obtain the following
statement.
\begin{proposition}\label{p_2}(on analytic a one-dimensional SIM).
For the fixed \,$m \ne 0$ and \,$a>0$, in the neighborhood of stationary point \eqref{s_p} in the phase
space \,$\mathbb R^3$ of the variables \,$\left(\Phi, \,\Phi^\prime, \,\Phi^{\prime\prime}\right)$, the values
of the solutions to the singular nonlinear CP \eqref{eqf_phi}, \eqref{ccon1f_phi} form a \,$\tau$-invariant
one-dimensional analytical SIM \,${\bf M}_{-} ^{(1)}(a, m)$ that is specified by two nonlinear relations
\begin{equation}\label{SIM}
{\bf M}_{-}^{(1)}(a, m): \quad \Phi + a = \rho _1
(\Phi^{\prime\prime}, a, m), \quad \Phi^\prime = \rho _2 (\Phi^{\prime\prime}, a, m).
\end{equation}
Here, \,$\{\rho _1(y), \,\rho _2(y) \}$ is a solution to the
Lyapunov-type singular nonlinear problem:
\begin{equation}\label{L_eq1}
\frac{d\rho_1}{dy}\left [ay + \frac {m-1}{m}\rho _2^2 - \rho _1y \right] =
\rho_2,\quad
\end{equation}
\begin{equation}\label{L_eq2}
\frac{d\rho_2}{dy}\left[ay + \frac{m-1}{m}\,\rho_2^2-\rho_1 y \right] = y, \quad  |y|<\Delta, \quad \Delta>0,
\end{equation}
\begin{equation}\label{conL}
\rho_1(0) = \rho_2(0)=0.
\end{equation}
The solution \,$\{\rho_1(y, a, m), \,\rho_2(y, a, m)\}$ to this problem exists and is unique and
holomorphic at the point \,$y=0$:

\begin{equation}\label{serS}
\rho _1(y)=\sum_{k=1}^{\infty}b_k y^k , \quad \rho
_2(y)=\sum_{k=1}^{\infty}c_k y^k , \qquad |y|< \Delta _0, \quad
\Delta _0 >0,
\end{equation}
\begin{equation}\label{b1_c1}
b_1 = 1/a^2, \qquad c_1=1/a,
\end{equation}
\begin{equation}\label{ck}
c_k= \left[\sum_{l=1}^{k-1} \left( l c_l b_{k-l} - \frac{m-1}{m}
\sum_{s=1}^{k-l} l c_l c_s c_{k-l-s+1} \right)\right]\Big/ (ak),
\end{equation}
\begin{equation}\label{bk}
b_k = \left [ c_k+ \sum_{l=1}^{k-1} \left ( l b_l b_{k-l} -
\frac{m-1}{m} \sum_{s=1}^{k-l} l b_l c_s c_{k-l-s+1} \right )
\right ] \Big / (ak), \quad k=2,3, \ldots \,;
\end{equation}
in particular \eqref{b1_c1}--\eqref{bk} imply
\begin{equation}\label{c2_b2}
c_2=1/(2ma^4), \qquad b_2=3/(4ma^5), \,\ldots \,.
\end{equation}
For the limit case \,$m\to\infty$,  the problem
\eqref{L_eq1}--\eqref{conL} has the exact solution,
\begin{equation}\label{ro_inf}
 \rho_1(y, a, \infty) =y/a^2, \qquad \rho_2(y, a, \infty)=y/a,
\end{equation}
so that, due to \eqref{SIM}, the one-dimensional SIM \,${\bf M}_-^{(1)}(a, \infty)$ becomes linear
and exists globally in \,$\mathbb{R}^3$, it is generated by the values of solutions  \eqref{phiInf}.
\end{proposition}
\begin{corollary}\label{c_1}. For any fixed \,$m\not=0$ and \,$a>0$, there
exists \,$T_0\gg 1$ such that, for the solutions of  Eq.\eqref{eqf_phi},
the limit condition \eqref{ccon1f_phi} is equivalent to the following two
nonlinear relations \,$\forall T \ge T_0$:
\begin{equation}\label{simT}
\Phi(-T) + a = \rho_1\Bigl (\Phi^{\prime\prime}(-T), a, m\Bigr),
\qquad \Phi^\prime(-T) = \rho_2 \Bigl(\Phi ^{\prime\prime}
(-T), a, m \Bigr).
\end{equation}
Here \,$\rho_1(y, a, m)$ and \,$\rho_2(y, a, m)$ are the same as in \eqref{SIM}, and
\,$|\Phi ^{\prime\prime} (-T_0)|$ is small enough.
For the limit case \,$m\to\infty$, taking into account \eqref{ro_inf} we
have the exact formulas
\begin{equation*}
\rho_1=\Phi^{\prime\prime}(-T)/a^2,\qquad \rho_2=\Phi^{\prime\prime}(-T)/a
\end{equation*}
so that the relations \eqref{simT} become linear.
\end{corollary}
\begin{remark}\label{r_6}. For the curve \eqref{simT}, the tangent taken at the
stationary point \,$(\Phi, \,\Phi^\prime, \,\Phi^{\prime\prime})_s (a) = (-a, \,0, \,0)$ is defined  by the
relations
\begin{equation}\label{rol_T}
\Phi(-T) + a - \Phi ^{\prime\prime}(-T)/a^2 = 0, \qquad
\Phi ^\prime (-T) - \Phi ^{\prime\prime} (-T)/a = 0.
\end{equation}
These relations give the linear approach to the nonlinear SIM
\eqref{simT}.

For  \,$m\to\infty$, \eqref{simT} coincides with \eqref{rol_T} (linear SIM for
a nonlinear ODE!).
\end{remark}
\begin{corollary}\label{c_2}. For any fixed \,$m\not=0$ and \,$a>0$, singular
nonlinear BVP \eqref{eqf_phi}, \eqref{ccon1f_phi},
\eqref{con2f_phi}, defined on \,$\mathbb{R}_-$,
is equivalent to the regular nonlinear two-point BVP \eqref{eqf_phi}, \eqref{simT}, \eqref{con2f_phi}
defined on a finite interval \,$[-T, 0]$, generally with a variable left endpoint.
\end{corollary}
\begin{corollary}\label{c_3}. For any fixed  \,$m\not=0$ and \,$a>0$, the
singular nonlinear IBVP \eqref{eqf_phi}, \eqref{ccon1f_phi}, \eqref{con2f_phi}, \eqref{con3f_phi},
defined on the entire real axis is equivalent to the nonlinear  IBVP
 \eqref{eqf_phi}, \eqref{simT}, \eqref{con2f_phi}, \eqref{con3f_phi} defined on a semi-infinite
interval \,$[-T, \infty)$, generally with a variable left endpoint, where the unknown parameter \,$a=a(b)>0$ in \eqref{simT} must be chosen to satisfy condition \eqref{con3f_phi} (with fixed \,$b>0$) when
this limit behavior is valid.
\end{corollary}
\subsection{Preliminary Conclusions and Remarks}
\subsubsection{Transfer of the limit boundary conditions from infinity to a finite point and the
shooting methods for solving an approximate nonlinear two-point BVP}
For the approximate transfer of the limit conditions \eqref{ccon1f_phi} to a
finite point, we retain in \eqref{simT} the principal terms on \,$\Phi^{\prime\prime} (-T)$
using \eqref{serS} with \,$y=\Phi^{\prime\prime} (-T)$.
For linear approach, we have relations \eqref{rol_T}. If we retain in \eqref{serS}
the terms on \,$y$ up to the second order inclusive, then due to formulas \eqref{b1_c1},
\eqref{c2_b2}, \eqref{simT} we obtain the following two approximate boundary conditions at
the point \,$\tau=-T$:
\begin{equation}\label{eq1_ST_ap}
\Phi(-T) + a - \Phi ^{\prime\prime}(-T)/a^2 = [3/(4a^5m)][\Phi^{\prime\prime}(-T)]^2,
\end{equation}
\begin{equation}\label{eq2_ST_ap}
\Phi ^\prime (-T) - \Phi ^{\prime\prime} (-T)/a = [1/(2a^4m)][\Phi^{\prime\prime}(-T)]^2.
\end{equation}
Thus, instead of singular nonlinear BVP \eqref{eqf_phi}, \eqref{ccon1f_phi},
\eqref{con2f_phi} defined on \,$\mathbb{R}_-$ and depending  on a parameter \,$a>0$,  we obtain
the approximate two-point BVP on \,$[-T,0]$, e.g., in the form \eqref{eqf_phi}, \eqref{eq1_ST_ap}, \eqref{eq2_ST_ap},
\eqref{con2f_phi} (or, for the linear approach, in the form \eqref{eqf_phi}, \eqref{rol_T}, \eqref{con2f_phi})
with the same parameter \,$a>0$.

For the fixed \,$a>0$ and \,$m\not=0$, when the singular nonlinear BVP \eqref{eqf_phi}, \eqref{ccon1f_phi},
\eqref{con2f_phi} is uniquely solvable, the stable shooting methods may be used to solve it numerically:

1) Starting from the approximate one-dimensional SIM defined by \eqref{eq1_ST_ap}, \eqref{eq2_ST_ap}
(or \eqref{rol_T}) for fixed \,$T\gg 1$, and solving rightwards the auxiliary CPs with
\,$\Phi ^{\prime\prime}(-T)$ as the parameter for shooting \,($|\Phi^{\prime\prime}(-T)|\ll 1$),
we find the value of this parameter  to satisfy condition \eqref{con2f_phi} at the point $\tau=0$.

2) As an alternative equivalent method we use the Lyapunov series \eqref{l_ser}, \eqref{hl} taken at
the point \,$\tau = -T$, and by solving rightwards the corresponding auxiliary CPs with \,$d>0$ as
the parameter for shooting, we choose the value of this parameter to satisfy condition \eqref{con2f_phi}
(such approach is rather convenient for shooting method).
\subsubsection{The scaling transformation to the singular nonlinear IBVP}
Let singular nonlinear IBVP \eqref{eqf_phi}, \eqref{ccon1f_phi}, \eqref{con2f_phi}, \eqref{con3f_phi}
defined on \,$\mathbb{R}$ be uniquely solvable, for some fixed \,$m>0$ and \,$b>0$, and has a solution
\,$\Phi_m(\tau,a)$, where \,$a=a(b)>0$.

In order to find this solution it suffices to solve the above IBVP with $a=1$. Indeed,
first we solve the associated singular nonlinear BVP  \eqref{eqf_phi}, \eqref{ccon1f_phi}, \eqref{con2f_phi}
defined on \,$\mathbb{R}_-$ with \,$a=1$ and obtain \,$d=d_m(1)>0$ and the corresponding
solution \,$\Phi _m(\tau, 1)$ (here and further \,$d=d_m(a)>0$ is a parameter in the Lyapunov series
\eqref{l_ser}, \eqref{hl}). Taking the solution \,$\Phi_m(\tau, 1)$ and extending it to \,$\tau>0$, as a
solution of CP with the obtained initial data at the point \,$\tau=0$, we derive the constant \,$b=b_m(1)>0$.

Then the needed \,$a=a(b)>0$, with given \,$b>0$ in \eqref{con3f_phi}, can
be obtained from the relations (due to the scaling transformations):
\begin{equation}\label{b_m}
b_m(a) = b_m(1)a^{m+1}>0, \qquad a=a(b)=[b/b_m(1)]^{1/(m+1)}>0.
\end{equation}

The desired solution \,$\Phi _m(\tau ,a)$ and the corresponding value of the Lyapunov parameter
\,$d_m(a)>0$, where \,$a=a(b)>0$ is defined in \eqref{b_m}, can be obtained by the scaling
transformations:
\begin{equation}\label{phim_dm}
\Phi _m(\tau, a) = a\Phi _m(a\tau, 1), \quad \tau\in \mathbb{R},
\quad d_m(a)=a d_m(1)>0.
\end{equation}

The values \,$b_m(1)>0$ and \,$d_m(1)>0$ cannot be found by
local analysis methods and have been determined numerically.
\subsubsection{On the convergence of the accompanying improper integrals}
The setting of conditions in the form of \eqref{ccon1f_phi} as \,$\tau\to-\infty$
ensures the convergence of the improper integrals
\begin{equation}\label{i_k}
I_k (\tau) =\int _{-\infty }^\tau \Phi^{(k)}(s)ds, \qquad k =
0,1,2,3,
\end{equation}
(and other ones required later on) for the solutions \,$\Phi (\tau, a)$ to the singular nonlinear
BVP  \eqref{eqf_phi}, \eqref{ccon1f_phi}, \eqref{con2f_phi}. Here, \,$\Phi^{(0)}(\tau)\equiv \Phi (\tau)$ and,
for fixed \,$k\ge 1$, \,$\Phi^{(k)}(\tau)$ is the \,$k$th derivative of \,$\Phi(\tau)$. This
allows us to transform Eq.\eqref{eqf_phi} by the integration over the interval \,$(-\infty, \tau )$ in order
to give  the overall analysis of the above singular nonlinear BVP (see further Section 3).

Using the Lyapunov indices, conditions \eqref{ccon1f_phi} can be rewritten in the form
\begin{equation}\label{l_ind}
\limsup_{\tau \to -\infty}{\frac {\ln{|\Phi (\tau) +a|}}{\tau}}=\limsup_{\tau \to -\infty}
{\frac {\ln{|\Phi ^\prime (\tau)|}}{\tau}}= \limsup_{\tau \to -\infty}{\frac {\ln{|\Phi ^{\prime\prime}
(\tau)|}}{\tau}}>\varepsilon
\end{equation}
$$\forall \varepsilon: \,  0<\varepsilon <a.$$
\subsubsection{The families of the exact regular and singular (blow-up) solutions to the initial
third-order nonlinear ODE}
For some fixed values of \,$m$, there exist (both the well known and partially new)
exact solutions to ODE  \eqref{eqf_phi} which are not the solutions to the
singular nonlinear IBVP  \eqref{eqf_phi}--\eqref{con3f_phi}.

Namely, along with the obvious solutions \,$\Phi(\tau) \equiv {\rm
const}$ \,$\forall m \in \mathbb R$, ODE \eqref{eqf_phi} has:

1) for each \,$m:(m \ne 0) \wedge (m \ne -1)$, the
one-parameter family of the blow-up solutions
\begin{equation}\label{phi1_sing}
\Phi^{(1)}_ {{\rm sing}, m} \left(\tau - \tau_p\right) = \frac {6m}{(m+1)(\tau-\tau_p)}\,,
\qquad \tau_p \in \mathbb{R};
\end{equation}
2) for \,$m=1/2$, the two-parameter family of the blow-up
solutions
\begin{equation}\label{phi2_sing}
\Phi^{(2)}_{{\rm sing}, 1/2}\left(\tau - \tau_p,a\right) = a\coth{\Bigl (a(\tau - \tau_p)/2
\Bigr)}, \qquad a,\tau_p\in \mathbb{R},
\end{equation}
which becomes the same as \eqref{phi1_sing} when $a=0$:
\begin{equation*}
 \Phi^{(2)}_{{\rm sing}, 1/2}\left(\tau - \tau_p,0\right) \equiv \Phi^{(1)}_ {{\rm
sing}, 1/2} \left(\tau - \tau_p\right);
\end{equation*}
3) for \,$m\in\{1/2; 1; 2; \infty\}$, the two-parameter families of
solutions \,$\Phi_{m}(\tau-\tau_s, a)$  (${a,\tau_s\in \mathbb{R}}$) existing on the entire
real line:
\begin{equation}\label{phi_05}
\Phi_{1/2}\left(\tau - \tau_s,a\right) = a\tanh{\Bigl (a(\tau - \tau_s)/2 \Bigr)},
\end{equation}
\begin{equation}\label{phi_1_2}
\Phi _1\left(\tau - \tau _s,a\right) = a\left(\tau - \tau _s\right), \qquad
\Phi _2\left(\tau - \tau _s,a\right) = a\left(\tau - \tau _s\right)^2,
\end{equation}
\begin{equation}\label{phi_inf}
\Phi_{\infty}\left(\tau - \tau_s,a\right) =
a\Bigl[\exp{\Bigl (a (\tau - \tau_s) \Bigr)} -1 \Bigr];
\end{equation}
4) for \,$m=1/3$, the two-parameter family of solutions obtained in the implicit form:
\begin{equation*}
\Bigl(\Phi_{1/3}(\tau)/|\Phi_{1/3}(\tau)|\Bigr)(\tau - \tau_s) =  \frac{1}{2 b^2}\,
\ln{\left(\frac{b^2 + b\,\sqrt{|\Phi_{1/3}(\tau)|} +
|\Phi_{1/3}(\tau)|}{\left(b -
\sqrt{|\Phi_{1/3}(\tau)|}\right)^2}\right)} +
\end{equation*}
\begin{equation}\label{phi_1/3}
 + \frac{\sqrt{3}}{b^2} \arctan{\left(\frac{2\, \sqrt{|\Phi_{1/3}(\tau)|} + b}{b\,
\sqrt{3}}\right)}, \quad b\not=0, \quad \tau_s\in \mathbb{R}
\end{equation}
(in \eqref{phi1_sing}--\eqref{phi_1/3},  the values \,$\tau_p$ and \,$\tau_s$ are the shift
parameters, and \,$a$ and \,$b$ are arbitrary numbers, \,$b\not=0$).

These particular solutions are obtained by simple order reduction in ODE \eqref{eqf_phi}
(see also here Appendix B).

For each fixed \,$a>0$ and \,$b=\sqrt{a}$, the functions
\,$\Phi^{(2)}_{{\rm sing}, 1/2}\left(\tau - \tau_p,a\right)$,
\,$\Phi_{1/2}\left(\tau - \tau_s,a\right)$, \,$\Phi_{\infty}\left(\tau - \tau_s,a\right)$ and
\,$\Phi_{1/3}\left(\tau - \tau_s,a\right)$ are the exact solutions
of singular CP \eqref{eqf_phi}, \eqref{ccon1f_phi}; they are represented by the Lyapunov series
\eqref{l_ser}, \eqref{hl} with the respective parameters \,${d=d_{{\rm sing}, 1/2}(a,
\tau_p)}$,  \,$d=d_{1/2}(a, \tau_s)$,   \,$d=d_{\infty}(a, \tau_s)$,
and \,$d=d_{1/3}(a, \tau_s)$ ($\tau_p, \tau_s\in \mathbb{R}$) defined by the formulas:
\begin{equation}\label{d_1/2_s}
d_{{\rm sing}, 1/2}(a, \tau_p)= - 2a\exp{(-a\tau_p)};
\end{equation}
\begin{equation}\label{d1/2_inf}
d_{1/2}(a, \tau_s)=2a\exp{(-a\tau_s)}; \qquad
d_{\infty}(a, \tau_s)=a\exp{(-a\tau_s)};
\end{equation}
\begin{equation}\label{d_1/3}
d_{1/3}(a, \tau_s)=2\sqrt{3}\,a\exp{\left(\sqrt{3}\pi/3 - a\tau_s\right)}.
\end{equation}

For singular CP \eqref{eqf_phi}, \eqref{ccon1f_phi} with \,$m=1/3$, if we put \,$\tau_s
= \widetilde\tau_s + \sqrt{3}\pi/(6a)$, then it is not difficult
to show that \,$\Phi_{1/3}\left(\widetilde\tau_s,a\right) =
\Phi_{1/3}^\prime\left(\widetilde\tau_s,a\right) = 0$,
\,$\Phi_{1/3}^{\prime\prime}\left(\widetilde\tau_s,a\right) =
-2a^3/9$, and \,$\Phi_{1/3}\left(\tau - \tau_s,a\right)$  does not
exist globally: it has a finite pole-point \,$\tau=\tau_p$, where
\,$\tau_p = \widetilde\tau_s + 2\sqrt{3}\pi/(3a) = \tau_s +
\sqrt{3}\pi/(2a)$ (in detail, see an example in Subsection 3.2.3 of Section 3).
\subsubsection{A necessary condition for the Lyapunov family of solutions
to exist globally (on the entire real line)}
In general, for any fixed \,$\tau_p \in \mathbb{R}$ and \,$m>0$,
there exists two-parameter family of blow-up solutions which tends to the exact
solution (1.61) when these two parameters vanish (in detail, see here Appendix A).

The existence of such singular solutions having the pole points
\,$\tau=\tau_p\in\mathbb{R}$ is closely related to the problem of the  global
existence of the solutions to Eq.\eqref{eqf_phi} (further we omit index $^{(1)}$ in
the exact solutions \eqref{phi1_sing}). Namely, it follows from  \eqref{phi1_sing} that:
\begin{equation}\label{phi_sing}
\Phi _{\rm sing}(\tau, \tau_p) <0, \quad \Phi ^{\prime}
_{\rm sing}(\tau, \tau_p) <0, \quad \Phi^{\prime\prime} _{\rm
sing}(\tau, \tau_p)<0, \quad \tau < \tau_p;
\end{equation}
\begin{equation}\label{phid_sing}
\Phi_{\rm sing}(\tau, \tau_p)>0, \quad \Phi^{\prime}
_{\rm sing}(\tau, \tau_p)<0, \quad \Phi^{\prime\prime} _{\rm
sing}(\tau, \tau_p)>0, \quad \tau>\tau_p.
\end{equation}

Then any solution to Eq.\eqref{eqf_phi} for which similar inequalities
hold at some point \,$\tau \in \mathbb {R}$ moves onto a pole, i.e.,
does not exist globally (blows up at some finite point \,$\tau =
\tau _p$).

In particular, it is valid for the solutions to the singular
nonlinear CP \eqref{eqf_phi}, \eqref{ccon1f_phi} of the Lyapunov series
\eqref{l_ser}, \eqref{hl} type  with \,$m>0$ and \,$d<0$ (see, e.g., \eqref{phi2_sing}, where \,$d$
is defined by \eqref{d_1/2_s}).
\begin{corollary}\label{c_4}.  For fixed \,$m>0$ and \,$a>0$, let \,$\Phi _m (\tau,
a, d)$ be a solution to the singular nonlinear CP \eqref{eqf_phi}, \eqref{ccon1f_phi},
where \,$d$ is the parameter in the Lyapunov expansion \eqref{l_ser}, \eqref{hl}.
Then the requirement \,$d>0$ is a necessary condition for a
solution \,$\Phi _m (\tau, a, d)$ to exist globally, i.e., on the
entire real line. In addition, we obtain \,$\Phi _m (\tau)<0$, \,$\Phi
_m^\prime (\tau)>0$, \,$\Phi _m^{\prime\prime}(\tau)>0$ at least for \,$\tau\ll -1$.
\end{corollary}
\section{The Accompanying Singular Nonlinear BVP on the Non-Positive
Semi-Axis and Its Univalent Solvability}
\subsection{Statement of Singular Nonlinear BVP  with Two Parameters and
Associated Integro-Differential Relations}
First, the following singular nonlinear BVP, with parameters \,$a>0$ and \,$m>0$, must be studied:
\begin{equation}\label{eq2_phi}
(\Phi ^{\prime\prime} + \Phi \Phi ^\prime)^\prime = [(2m-1)/m]\, (\Phi ^\prime )^2,
\qquad -\infty <\tau \le 0,
\end{equation}
\begin{equation}\label{eq2_inf}
\begin{array}{c}
\lim_{\tau \to -\infty}{ \{ \exp{(-\varepsilon \tau)} [\Phi (\tau
)+a] \}} = \lim_{\tau \to -\infty}{[\exp{(-\varepsilon \tau)}
\Phi ^\prime(\tau )]}= \\ \\
= \lim_{\tau \to -\infty}{ [ \exp{(-\varepsilon \tau)} \Phi
^{\prime\prime}(\tau )]}= 0 \qquad \forall \varepsilon : 0 <
\varepsilon < a,
\end{array}
\end{equation}
\begin{equation}\label{eq2_0}
 \Phi(0) = 0.
\end{equation}
Here, Eq.\eqref{eqf_phi} is rewritten in form \eqref{eq2_phi} because it is
convenient to analyze a global behavior of the solutions \,$\Phi _m
(\tau, a, d)$ to the singular nonlinear CP \eqref{eq2_phi}, \eqref{eq2_inf}; these
solutions are represented by the Lyapunov  exponential series
\eqref{l_ser}, \eqref{hl} with a parameter \,$d$, where \,$|d\exp{(a\tau)}|$ is
small enough.

For the above solutions, we integrate twice both sides of \eqref{eq2_phi}
from \,$-\infty$ to \,$\tau$ and take into account that
\begin{equation}\label{int_dphi}
\int\limits_{-\infty}^\tau {\int\limits_{-\infty}^s {[\Phi ^\prime
(t)] ^2 \,dt} \,ds}= \int\limits_{-\infty}^\tau {(\tau-s)[\Phi
^\prime (s)] ^2 \,ds}
\end{equation}
(this equality is obtained by the integration by parts of the integral in the left hand-side
of  \eqref{int_dphi}).

Thus, for any fixed \,$a>0$ and \,$m\not=0$, the following relations
are valid for the solutions \,$\Phi_m (\tau, a, d)$ to the singular
nonlinear CP \eqref{eq2_phi}, \eqref{eq2_inf}, where $d$ is a parameter in the
Lyapunov exponential series \eqref{l_ser}, \eqref{hl}:
\begin{equation}\label{intdif1}
\Phi_m ^{\prime\prime}(\tau, a, d) + \Phi_m (\tau, a, d) \Phi_m^\prime (\tau, a, d) =
[(2m-1)/m] \int\limits_{-\infty}^\tau {[\Phi_m ^\prime (t, a, d)]^2  dt}\,,
\end{equation}
\begin{equation}\label{intdif2}
\Phi_m ^\prime (\tau, a, d) = [a^2 - \Phi_m ^2 (\tau, a, d)]/2 +
[(2m-1)/m]\,\int\limits_{-\infty}^\tau {(\tau-s) [\Phi_m ^\prime (s, a, d)] ^2 ds}\,.
\end{equation}

Then, for \,$m>1/2$ and $d>0$, due to these relations and Corollary \ref{c_4}, we obtain
\,$\Phi_m^\prime (\tau, a, d)>0$ \,$\forall\tau\in\mathbb{R}$. Indeed, otherwise there exists
\,$\tau=\widetilde\tau$ such that \,$\Phi_m ^\prime(\widetilde\tau)=0$ and \,$\Phi_m
^{\prime\prime}(\widetilde\tau) \le 0$ but \eqref{intdif1}  implies \,$\Phi_m ^{\prime\prime}(\widetilde\tau) >0$.
\begin{corollary}\label{c_5}.  For any fixed  \,$m>1/2$, \,$a>0$, and \,$d>0$, let
\,$\Phi _m (\tau, a, d)$ be a solution to the singular nonlinear CP
\eqref{eq2_phi}, \eqref{eq2_inf}, where \,$d>0$ is a parameter in the Lyapunov
exponential series \eqref{l_ser}, \eqref{hl}. Then such solution exists
globally on \,$\mathbb{R}$ and is a strictly increasing  function.
Moreover, \,$\Phi _m (\tau, a, d)$ is a convex function at least
while it remains non-positive.
\end{corollary}

Now let \,$\Phi_m(\tau, a)$  be a solution to the singular nonlinear BVP \eqref{eq2_phi}--\eqref{eq2_0}.
Then it is represented by the Lyapunov exponential series \eqref{l_ser}, \eqref{hl} with some
\,$d=d(a, m)>0$, and relations \eqref{intdif1}, \eqref{intdif2} together with condition \eqref{eq2_0} imply:
\begin{equation}\label{intd01}
\Phi_m^{\prime\prime} (0, a)=[(2m-1)/m] \int\limits_{-\infty}^0
{[\Phi_m ^\prime (t, a)]^2\,dt},
\end{equation}
\begin{equation}\label{intd02}
\Phi _m  ^\prime (0, a) = a^2/2 - [(2m-1)/m] \int\limits_{-\infty}^0 {s\,[\Phi_m ^\prime
(s, a)]^2\,ds}.
\end{equation}

From these relations, for the solutions \,$\Phi_m(\tau, a)$ to the
singular nonlinear BVP \eqref{eq2_phi}--\eqref{eq2_0}, we obtain:
\,(i) if \,$m: 0<m<1/2$, then  \,$\Phi_m^{\prime\prime}(0)<0$ and
\,$\Phi_m^{\prime}(0)<a^2/2$; \,(ii) \,$\Phi_{1/2}^{\prime\prime}(0)=0$ and
\,$\Phi_{1/2}^{\prime}(0)=a^2/2>0;$ \,(iii) if \,$m>1/2$, then  \,$\Phi_m^{\prime\prime}(0)>0$ and
\,${\Phi_m^{\prime}(0)>a^2/2>0}$.
\begin{remark}\label{r_7}.  For fixed $m: 0<m<1/2$, if there exists a
solution $\Phi _m (\tau)$ to BVP \eqref{eq2_phi}--\eqref{eq2_0}, then it must be
$\Phi _m ^{\prime} (0) \ge 0$. Otherwise, if $\Phi _m ^{\prime}
(\tau)$ changes a sign, for some $\tau<0$, then such solution moves
onto a pole-point (blows up at some finite point $\tau=\tau_p$)
because the inequalities \eqref{phi_sing} hold.
\end{remark}
\subsection{Some Exact Solutions of Singular Nonlinear BVP and Their Independent
Physical Meaning}
For some values of $m$, namely $m\in\{1/3,1/2,\infty\}$, there exist the exact solutions
$\Phi_m(\tau,a)$ to the singular nonlinear BVP \eqref{eq2_phi}--\eqref{eq2_0} having in particular
an independent physical meaning.

Mathematical settings of these physical problems in terms of
BVPs of this type are more accurate than those available in physical papers.

Moreover, Remark \ref{r_3} with the solution \,$\Phi_m (\tau,b)$ of basic singular nonlinear IBVP
\eqref{eq2_phi}--\eqref{eq2_0} replaced by a solution \,$\Phi_m (\tau,a)$ of the indicated singular nonlinear BVP
remains true for the solution \,$\Phi_m(\tau,a)$ of singular nonlinear BVP  \eqref{eqf_phi}--\eqref{con3f_phi}.
Hence the formulas are valid for expressing the stream function and the flow velocity components
via the solution \,$\Phi_m (\tau,a)$ and its derivative, and the methods described in Remark \ref{r_3} remain applicable
for calculating flows in \,$\{x,y\}$ plane (we would not give any more details in this section).

Here, for considered below Problems 1--3, a double enumeration of figures is given, where the first number is a problem number. It is convenient to distinguish them from main figures relating to the initial basic IBVP.
\subsubsection{Problem 1: a plane laminar "flooded jet"\, ($m=1/2$)}
For \,$m=1/2$, the exact solution \,$\Phi _{1/2} (\tau, a)$ exists
globally on \,$\mathbb{R}$ and is given by the formula
\begin{equation}\label{phi05}
\Phi_{1/2}(\tau,a) = a \tanh{\left(a\tau/2\right)}, \qquad \tau\in \mathbb{R}.
\end{equation}
For data at the point $\tau=0$ and for the parameter $d=d_{1/2}(a)$ in the Lyapunov series
\eqref{l_ser}, \eqref{hl}, we have
\begin{equation}\label{phi05_0}
\Phi_{1/2}(0)=\Phi_{1/2}^{\prime\prime}(0)=0, \qquad
\Phi_{1/2}^\prime(0)=a^2/2, \qquad d_{1/2}(a)=2 a.
\end{equation}
For this case, singular nonlinear BVP \eqref{eq2_phi}--\eqref{eq2_0} considered on
the entire real axis (for $\tau\to\infty$, the limit conditions of
\eqref{eq2_inf}-type are posed with the replacement $\tau$ by $-\tau$ and
$\Phi$ by $-\Phi$) corresponds to the problem for \textbf{an
unbounded plane laminar jet ("a flooded jet")}; it is the well-known physical problem
(see, e.g., \cite{shlich}--\cite{loit_73} and references therein).

In this problem, a shift parameter is fixed by the condition
\eqref{eq2_0}, and the parameter \,$a>0$ is defined by the fixed value of the
integral
\begin{equation}\label{phi05_int}
I_{1/2}(a)=\int_{-\infty}^{\infty}{[\Phi^\prime_{1/2}(\tau, a)]^2 \,
d\tau}=2a^3/3.
\end{equation}

For the variables \,$x,y$: \,$x>0$, \,$y\in\mathbb{R}$, we obtain:
\,$\tau_{1/2}(x,y)=y x^{-2/3}/\sqrt{3\nu}$,
\begin{equation}\label{psi05_xy}
\psi_{1/2}(x,y,a)=a \sqrt{3\nu} \,x^{1/3} \tanh{\left(a y x^{-2/3}/
(2 \sqrt{3\nu})\right)};
\end{equation}
\begin{equation}\label{u05_xy}
u_{1/2}(x,y,a)=(a^2/2) \,x^{-1/3}\Big/\cosh^2{\left(a y x^{-2/3}/
(2 \,\sqrt{3\nu}) \right)},
\end{equation}
\begin{equation}\label{v05_xy}
\begin{array}{c}
v_{1/2}(x,y,a)=a \sqrt{\nu/3} \,x^{-2/3}\Bigl[\left(a/ \sqrt{3\nu}\right)
y x^{-2/3}\Big/\cosh^2{\left(a y x^{-2/3}/(2\sqrt{3\nu}) \right)}-
\\ \\
-\tanh{\left(a y x^{-2/3}/ (2 \sqrt{3\nu}) \right)}\Bigr].
\end{array}
\end{equation}

Figs.1.1, 1.2 and Table 1 show the numerical results for \,$a=\nu=1$. Authors are not aware whether
such numerical results were obtained before. E.g., in \cite{shlich}--\cite{loit_73} only illustrations
of qualitative character are given for horizontal velocity profiles of the flows (cf. Fig.1.1 and Figs.1.3--1.5),
and vertical velocity profiles are not presented in any of these monographs.
\begin{center}
\includegraphics[width=10cm, height=6cm]{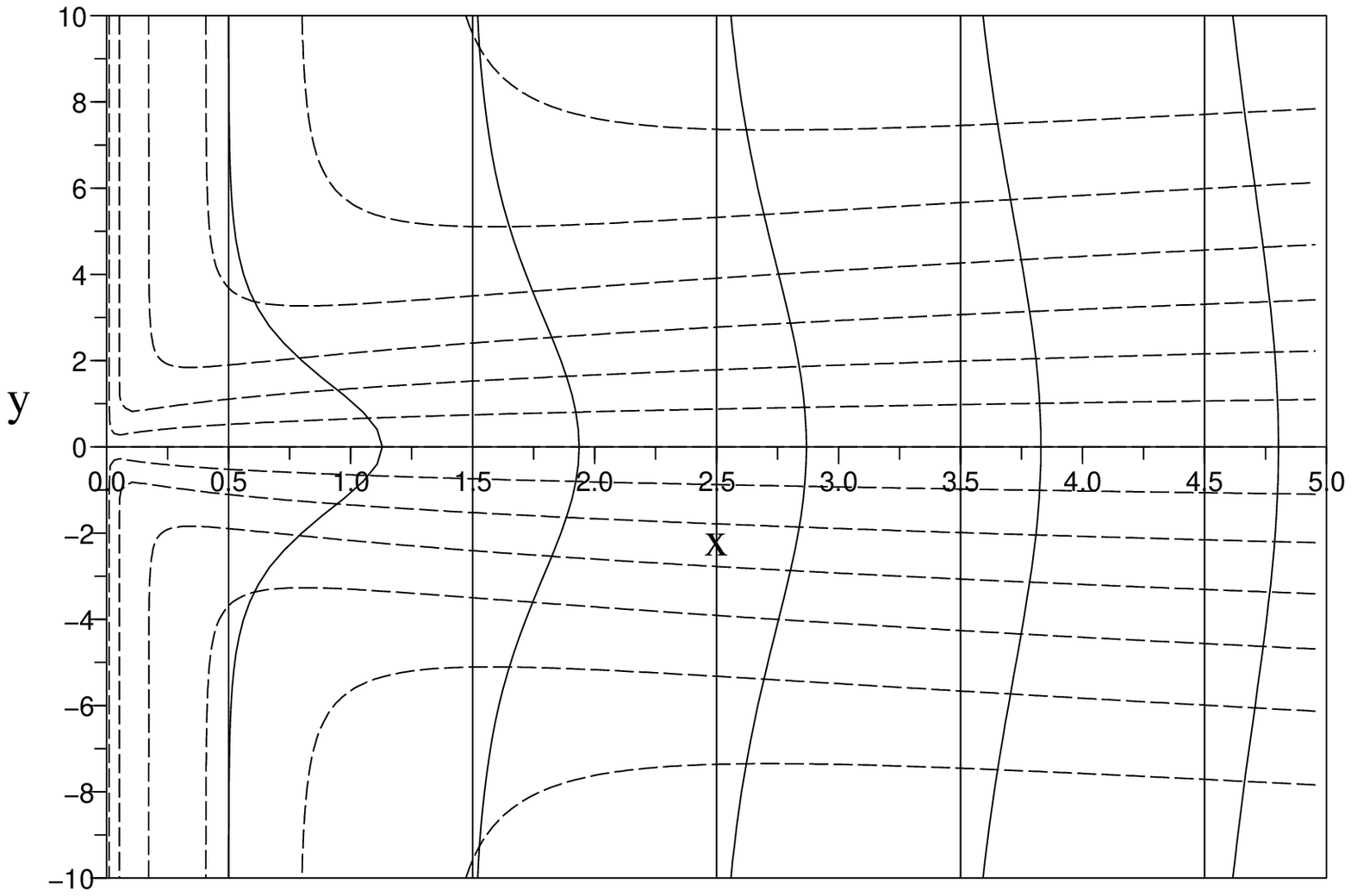}
\end{center}
\begin{center}
Fig.1.1. The curves \,$\psi_{1/2}(x,y)=$const (dotted lines), and the profiles
of the horizontal velocity component \,$\widetilde{u}_{1/2}(y)=u_{1/2}(x,y)|_{x={\rm const}}$.
\end{center}
\begin{center}
\includegraphics[width=10cm, height=6cm]{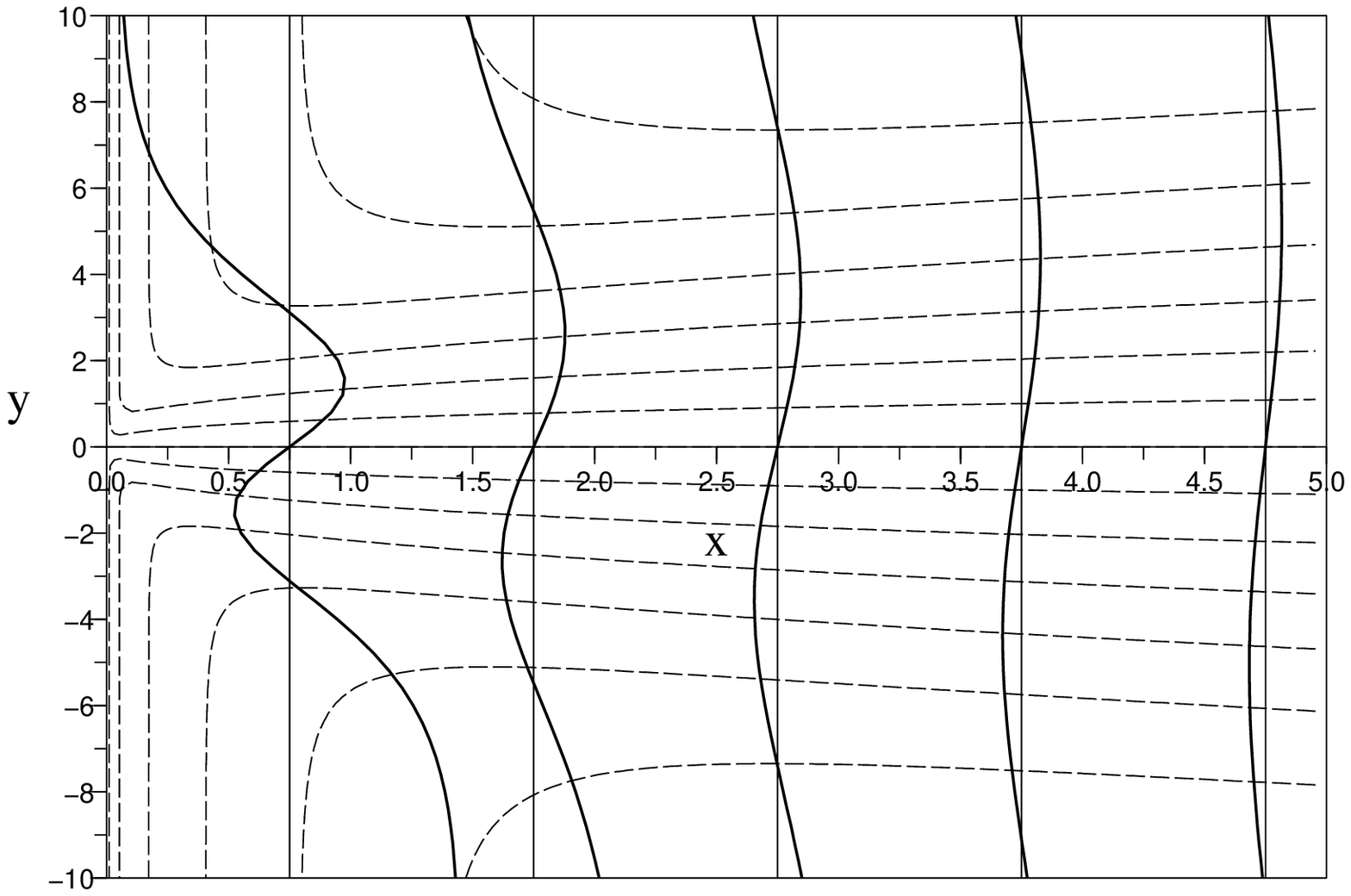}
\end{center}
\begin{center}
Fig.1.2. The curves \,$\psi_{1/2}(x,y)=$const (dotted lines), and the profiles
of the vertical velocity component \,$\widetilde{v}_{1/2}(y)=v_{1/2}(x,y)|_{x={\rm const}}$.
\end{center}

\begin{center}
Table 1 (to Fig.1.2)

\vspace{2mm}

\begin{tabular}{|l|l|l|l|l|l|}
\hline
$x$ &  $y_0(x)$ &  $y_{\rm max}(x)$ &  $v_{1/2, {\rm max}}(x)$ &  $v_{1/2, {\rm lim}}(x)$ &
$\widetilde{v}_{1/2, {\rm lim}}(x)$ \\
\hline
$0.75$ & $3.11308$ & $1.49215$ & $0.227294$ & $-0.69941$  & $-0.69941$ \\
$1.75$ & $5.47656$ & $2.62501$ & $0.129202$ & $-0.39757$  & $-0.39757$ \\
$2.75$ & $7.40238$ & $3.54809$ & $0.0955887$ & $-0.294138$  & $-0.294138$ \\
$3.75$ & $9.1027$ & $4.36308$ & $0.0777334$ & $-0.239195$  & $-0.239194$ \\
$4.75$ & $10.6564$ & $5.10781$ & $0.0663997$ & $-0.20432$  & $-0.204312$ \\
\hline
\end{tabular}
\end{center}

\vspace{2mm}

In Table 1, we use the following notation: $y_0(x): \,v_{1/2}(x,y_0)=0$; \quad
\begin{equation*}
y_{\rm max}(x), v_{1/2,{\rm max}}(x): v_{1/2,{\rm max}}(x)=
v_{1/2}(x, y_{\rm max}(x))=\max_{y\in [0, y_0]} {v_{1/2}(x, y)};
\end{equation*}
\begin{equation*}
v_{1/2,{\rm lim}}(x)=\lim_{y\to+\infty} {v_{1/2}(x, y)}, \quad \widetilde{v}_{1/2,{\rm lim}}(x)=v_{1/2}(x, 70).
\end{equation*}
\begin{center}
\includegraphics[width=10cm, height=6cm]{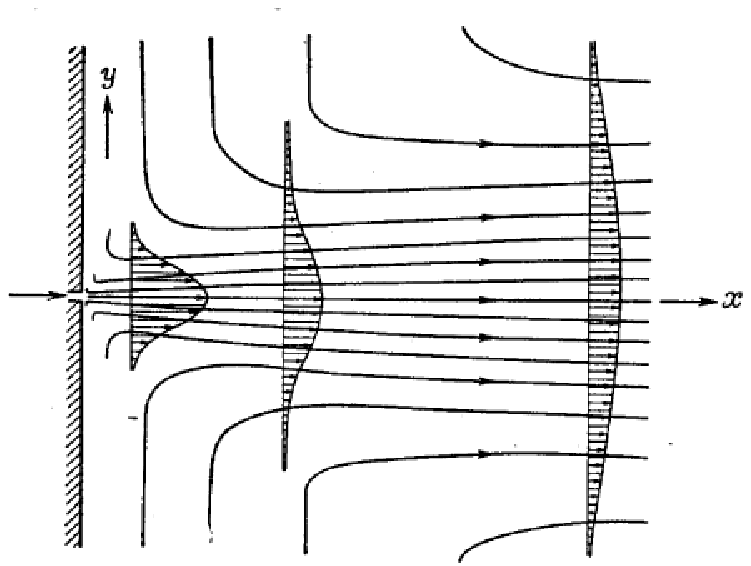}
\end{center}
\begin{center}
Fig.1.3 (\cite{shlich}, p.177)
\end{center}
\begin{center}
\includegraphics[width=6cm, height=4cm]{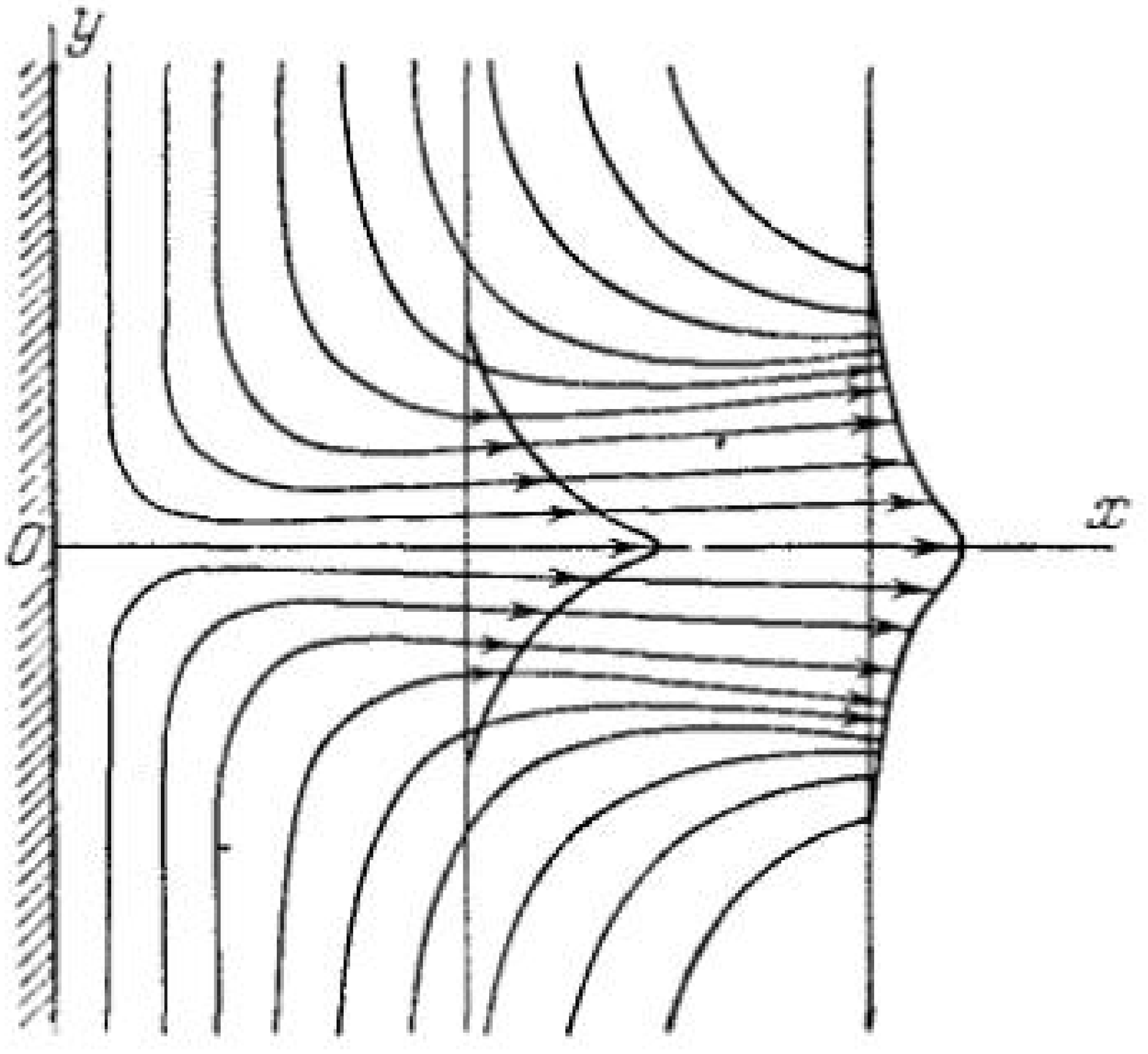}
\end{center}
\begin{center}
Fig.1.4 (\cite{slez_55}, p.287)
\end{center}
\begin{center}
\includegraphics[width=6cm, height=4cm]{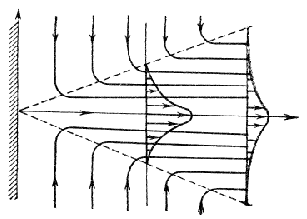}
\end{center}
\begin{center}
Fig.1.5 (\cite{loit_73}, p.532)
\end{center}
\subsubsection{Problem 2: it is connected with a non-stationary laminar
boundary layer separation ($m\to\infty$)}
For \,$m\to\infty$, the exact solution \,$\Phi_\infty(\tau, a)$ to the
singular nonlinear BVP \eqref{eq2_phi}--\eqref{eq2_0} exists globally on \,$\mathbb{R}$ and is
given by the formula
\begin{equation}\label{phi_inf_d}
\Phi_\infty (\tau,a)=a \left[\exp{\left(a\tau\right)}-1\right], \qquad \tau\in\mathbb{R}.
\end{equation}
For data at the point \,$\tau=0$ and for the parameter \,$d=d_{\infty}(a)$ in the Lyapunov series
\eqref{l_ser}, \eqref{hl}, we have
\begin{equation}\label{phiinf_0}
\Phi_\infty (0)=0, \qquad \Phi_\infty^\prime(0)=a^2, \qquad
\Phi_\infty^{\prime\prime}(0)=a^3, \qquad d_{\infty}(a)=a.
\end{equation}

For this limit case with \,$a=1$, according to the remark in \cite{dies2_86} with the reference to \cite{sych_83},
singular nonlinear BVP \eqref{eq2_phi}--\eqref{eq2_0} is connected with the problem concerning
\textbf{a non-stationary laminar boundary layer separation}.

For the variables \,$x,y$:  \,$x\ge 0$, \,$y\in\mathbb{R}$,  we have:
\begin{equation}\label{tpsiinf_xy}
\tau_{\infty}(x,y)\equiv y/\sqrt{\nu}, \quad \psi_{\infty}(x,y,a)=a \sqrt{\nu}\, x\,
\left[\exp{\left(a\,y/\sqrt{\nu}\right)}-1\right];
\end{equation}
\begin{equation}\label{uvinf_xy}
u_{\infty}(x,y,a)=a^2 x \exp{\left(a y/\sqrt{\nu}\right)},
\quad v_{\infty}(x,y,a)\equiv a \sqrt{\nu} \left[1-\exp{\left(a y/\sqrt{\nu}\right)}\right],
\end{equation}
and \,$\psi_{\infty}(x,y,a)=$const implies the formula
\begin{equation}\label{y_inf_x}
y(x,a)=(\sqrt{\nu}/a) \ln{\Bigl(1 + \mbox{const}/(a x)\Bigr)}, \quad x>0.
\end{equation}

Then, according to the Prandtl definition,  we may consider formally the value \,${x=x_0=0}$
as a point of a boundary layer separation because we obtain
\begin{equation}\label{separ_inf}
u_\infty(x_0,0)=0, \qquad \frac{\partial u_\infty}{\partial y}(x_0,0)=0.
\end{equation}

Figs.2.1--2.3 show the numerical results with \,$a=\nu=1$.
We do not aware of any other computational results of this kind
and/or of their more accurate physical interpretation.
\begin{center}
\includegraphics[width=10cm,height=6cm]{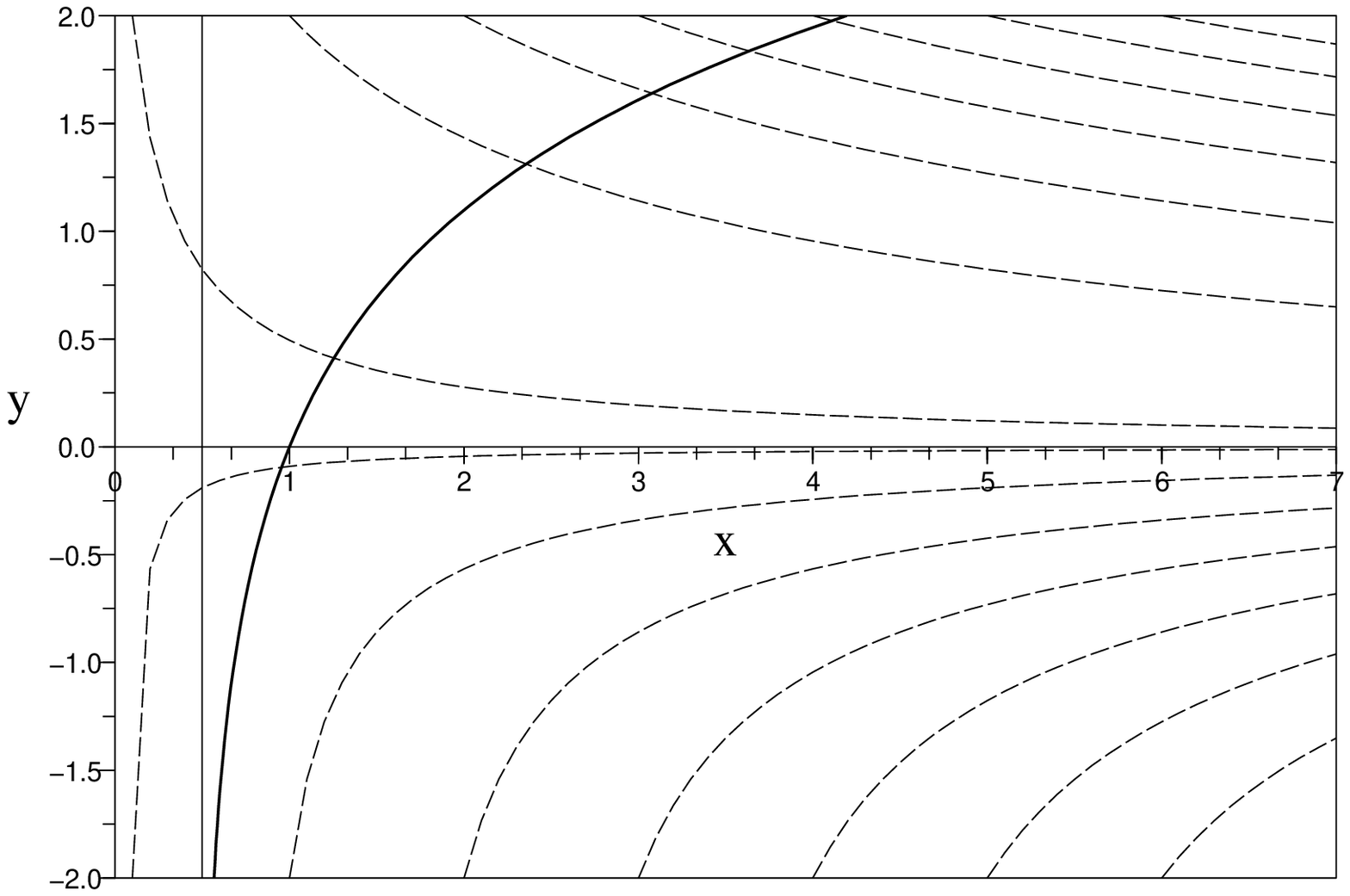}
\end{center}
\begin{center}
Fig.2.1. The curves $\psi_{\infty}(x,y)=$const (dotted lines), and one profile
of the horizontal velocity component $\widetilde{u}_{\infty}(y)=u_{\infty}(0.5, y)$.
\end{center}
\begin{center}
\includegraphics[width=10cm,height=6cm]{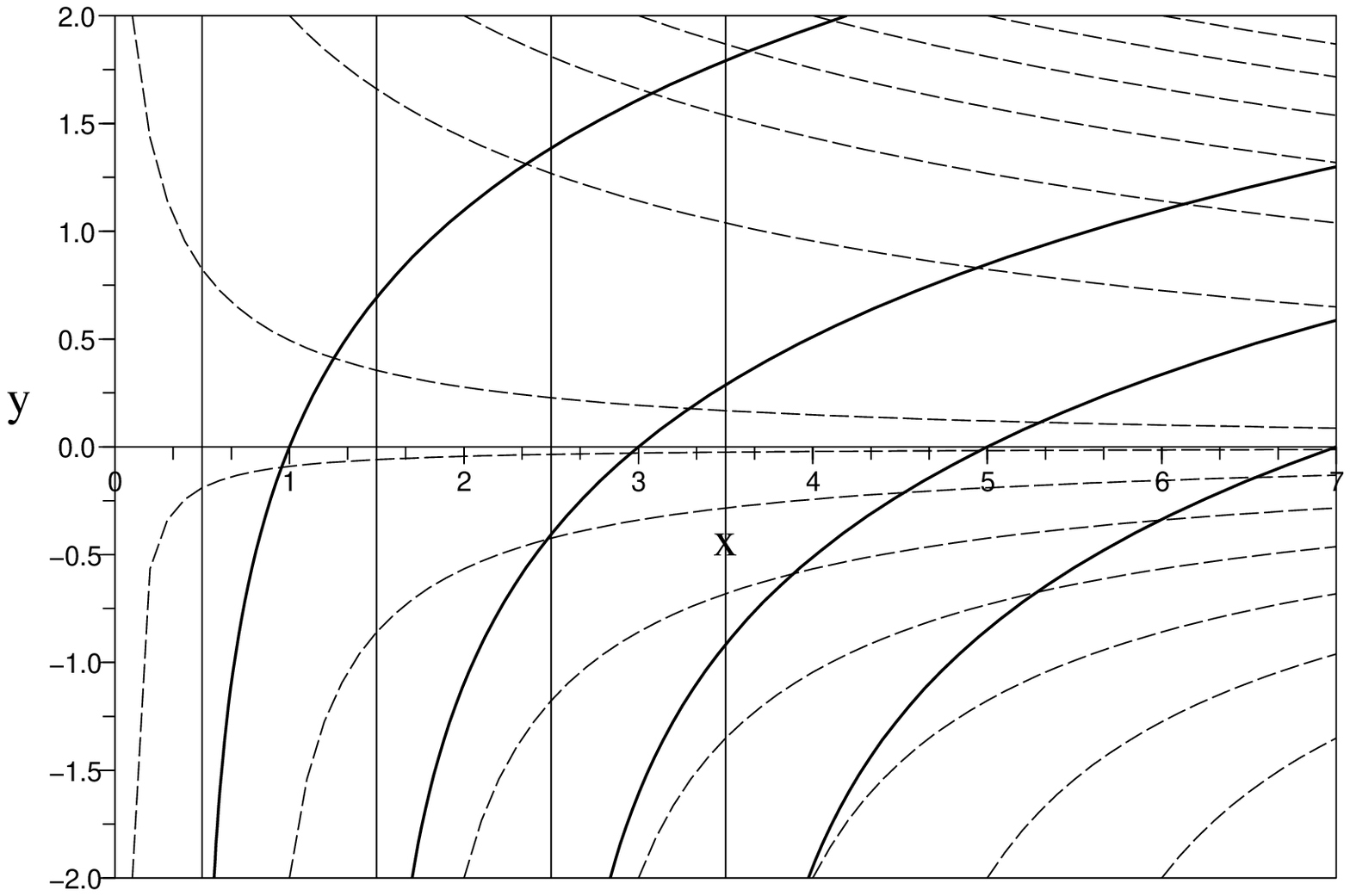}
\end{center}
\begin{center}
Fig.2.2. The curves \,$\psi_{\infty}(x,y)=$const (dotted lines), and the profiles
of the horizontal velocity component $\widetilde{u}_{\infty}(y)=u_{\infty}(x, y)|_{x={\rm const}}$.
\end{center}
\begin{center}
\includegraphics[width=10cm,height=6cm]{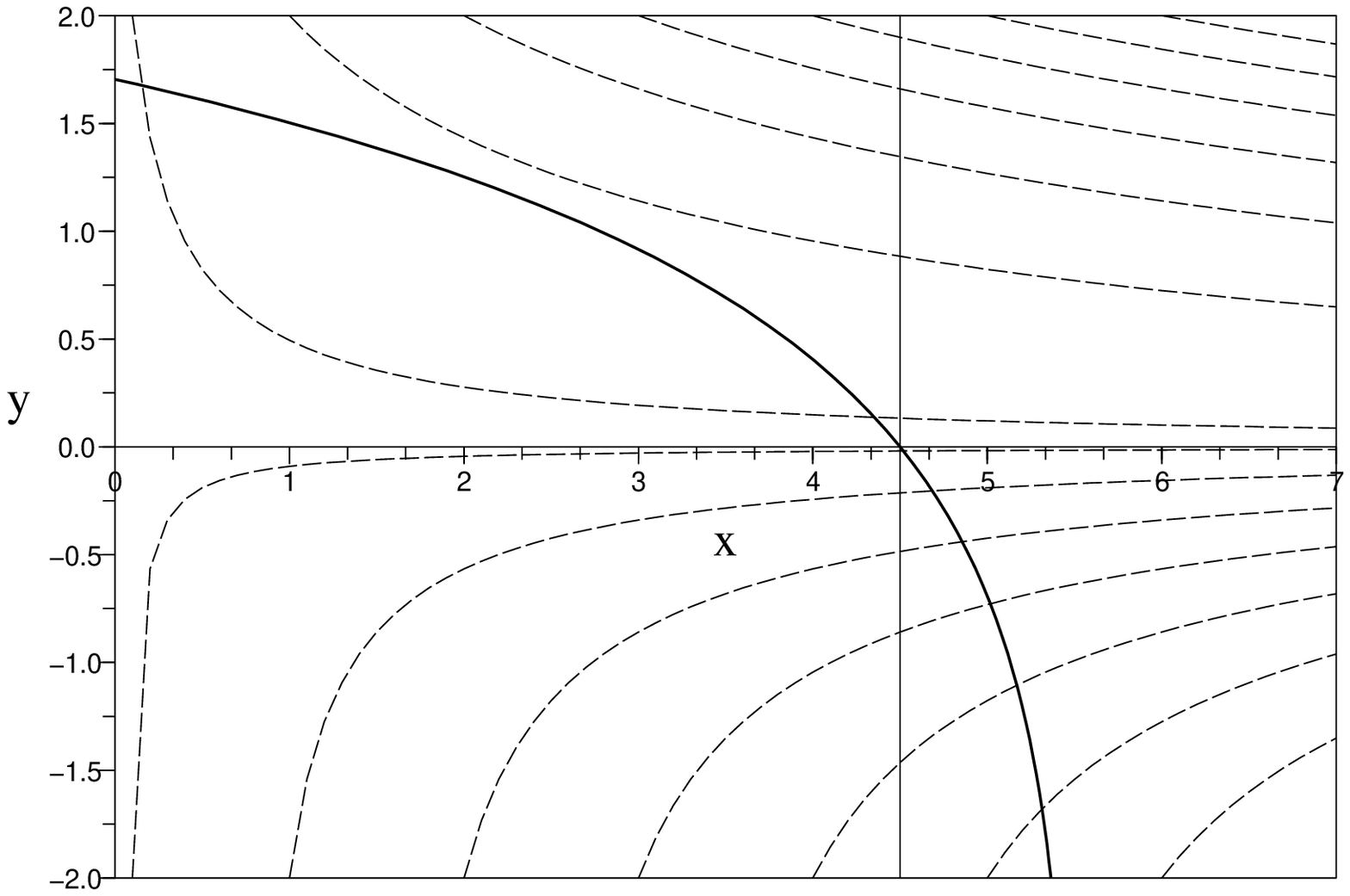}
\end{center}
\begin{center}
Fig.2.3. The curves \,$\psi_{\infty}(x,y)=$const (dotted lines), and the profile
of the vertical velocity component independent of \,$x$: \,$\widetilde{v}_{\infty}(y)=v_{\infty}(4.5,y)
\equiv v_{\infty}(x,y)$.
\end{center}
\subsubsection{Problem 3: a plane laminar "near-wall jet" \,($m=1/3$)}
For \,$m=1/3$, the exact solution \,$\Phi_{1/3}(\tau,a)$ to the singular nonlinear BVP
\eqref{eq2_phi}--\eqref{eq2_0} doesn't exist globally on \,$\mathbb{R}$ (for some \,$\tau=\tau_p>0$, it has
a singularity) and is given by the implicit formulas (in detail, see \cite{dks_07}):
\begin{equation*}
\tau (\Phi_{1/3})=\frac{\sqrt{3} \pi}{6a}-\frac{1}{2a}\ln{\left(
\frac{a+\sqrt{-a\Phi_{1/3}}-\Phi_{1/3}} {[\sqrt{a}-\sqrt{-\Phi_{1/3}}]^2} \right)}-
\end{equation*}
\begin{equation}\label{tau1/3}
-\frac{\sqrt{3}}{a} \arctan{\left(\frac{\sqrt{a}+2\sqrt{-\Phi _{1/3}}}{\sqrt{3a}}\right)}, \qquad
 -\infty <\tau \le 0, \qquad -a <\Phi _{1/3} (\tau) \le 0;
\end{equation}
\begin{equation*}
\tau (\Phi_{1/3})=\frac{\sqrt{3} \pi}{6a}-\frac{1}{2a}\ln{\left(\frac{a-\sqrt{-a\Phi_{1/3}}-\Phi_{1/3}}
{[\sqrt{a}+\sqrt{-\Phi_{1/3}}]^2}\right)} -
\end{equation*}
\begin{equation}\label{tau_1/3}
-\frac{\sqrt{3}}{a} \arctan \left( \frac{\sqrt{a}-2\sqrt{-\Phi_{1/3}}}{\sqrt{3a}}\right), \quad
0<\tau<\tau_p=2\pi \sqrt{3}/(3a)\approx 3.6275987/a,
\end{equation}
\begin{equation}\label{lim_taup}
\Phi_{1/3}(\tau)<0, \qquad \lim_{\tau\to\tau_p-0}{\Phi _{1/3}(\tau)}=-\infty.
\end{equation}

For the point \,$\tau=0$ and the inflection point \,$\tau=\tau_{\rm in}(a)\in\mathbb{R}_-$, we obtain the relations:
\begin{equation}\label{phi1/3_0}
\Phi_{1/3}(0,a)=0, \quad \Phi^\prime_{1/3}(0,a)=a^2/2+\int\limits_{-\infty}^0 {s\,[\Phi_{1/3}^\prime(s,a)]^2\,ds}=0,
\end{equation}
\begin{equation}\label{dphi1/3_0}
\Phi^{\prime\prime}_{1/3}(0,a)=-\int\limits_{-\infty}^0 [\Phi_{1/3}^\prime(s,a)]^2\,ds=-2a^3/9,
\end{equation}
\begin{equation}\label{phi1/3_in}
\Phi_{1/3}(\tau_{\rm in},a)\Phi^\prime_{1/3}(\tau_{\rm in},a)=-\int\limits_{-\infty}^{\tau_{\rm in}}{[\Phi_{1/3}^\prime(s,a)]^2\,ds}, \quad \tau_{\rm in}\in\mathbb{R}_-,
\end{equation}
and for the parameter \,$d=d_{1/3}(a)$ in the Lyapunov series \eqref{l_ser}, \eqref{hl}, there exists
the exact formula
\begin{equation}\label{d1/3}
d_{1/3}(a)=2\sqrt{3}\,a\,\exp{\left(\sqrt{3}\,\pi/6\right)}\approx 8.579306\,a.
\end{equation}

For the variables \,$x,\,y$: \,$x>0$, \,$y\in\mathbb{R}$, and \,$\tau$: \,$-\infty<\tau<\tau_p$, we obtain:
\begin{equation}\label{tau_psi1/3}
\tau_{1/3}(x,y)=y x^{-3/4}/(2\sqrt{\nu}), \qquad \psi_{1/3}(x,y,a)=2 \sqrt{\nu}\,x^{1/4}\,\Phi _{1/3}(\tau,a);
\end{equation}
\begin{equation}\label{u1/3_xy}
u_{1/3}(x,y,a)=x^{-1/2}\, \Phi^{\prime}_{1/3}(\tau, a),
\end{equation}
\begin{equation}\label{v1/3_xy}
v_{1/3}(x,y,a)=(3\sqrt{\nu}/2)\,x^{-3/4}\,\left[\tau \Phi^{\prime}_{1/3}(\tau,a)-\Phi_{1/3}(\tau,a)/3 \right].
\end{equation}

From  BVP \eqref{eq2_phi}--\eqref{eq2_0} with \,$m=1/3$, replacing \,$\tau$ by \,$-\tau$
and \,$\Phi$ by $-\Phi$, including formula \eqref{tau1/3}, we obtain the
problem concerning \textbf{unbounded jet near a wall} (see \cite{loit_73}).
A shift parameter is fixed by condition \eqref{eq2_0}, and the parameter value
\,$a>0$ is defined when a value of the integral
\begin{equation}\label{int1/3}
I_{1/3}(a) = \int_0^{\infty}{\Phi_{1/3}(\tau, a)[\Phi^\prime_{1/3}(\tau, a)]^2\,d\tau}
\end{equation}
is given.

This model is described in \cite{loit_73}, ננ.541--543, including the main formula \eqref{tau1/3}
(with different notation and slightly transformed non-logarithmic part what we in \cite{kss_09} considered as a partial mismatch). Formulas \eqref{tau_1/3} are not given in \cite{loit_73}. We do not aware of any other articles or monographs
studying this problem.

Figs.3.1, 3.2 show the numerical results with \,$a=\nu=1$, where the above change of variables
is taken into account. We do not know of any computations of this kind. In \cite{loit_73}, only an illustration picture
(see Fig.3.3 in this paper) is given which does not express the true behavior of the flow
(cf. Figs.3.1, 3.2).
\begin{center}
\includegraphics[width = 10cm,height=6cm]{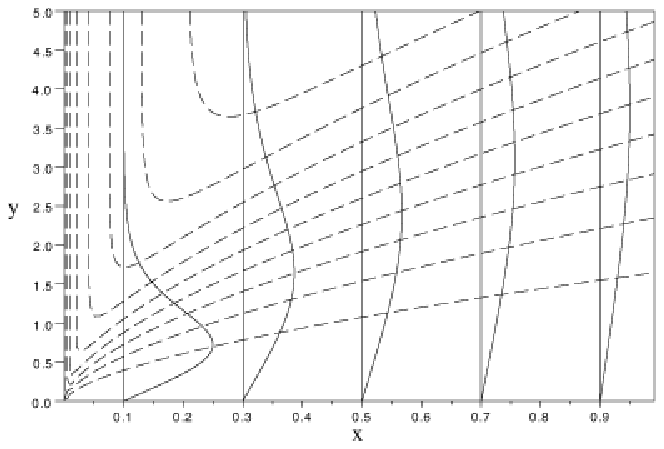}
\end{center}
\begin{center}
Fig.3.1. The curves \,$\psi_{1/3}(x,y)=$const (dotted lines), and the profiles
of the horizontal velocity component with the scale factor: \,$\widetilde{u}_{1/3}(y)=0.15\,
u_{1/3}(x, y)|_{x={\rm const}}$.
\end{center}
\begin{center}
\includegraphics[width = 10cm,height=6cm]{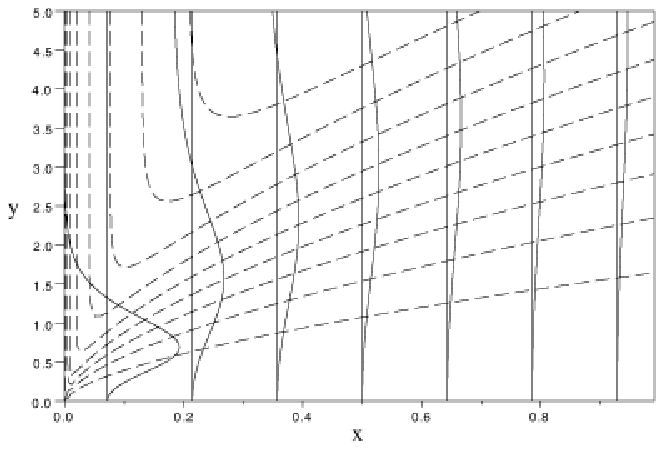}
\end{center}
\begin{center}
Fig.3.2. The curves  \,$\psi_{1/3}(x,y)=$const (dotted lines), and the profiles
of the vertical velocity component with the scale factor: $\widetilde{v}_{1/3}(y)=0.02\,
v_{1/3}(x, y)|_{x={\rm const}}$.
\end{center}
\begin{center}
\includegraphics[width = 10cm,height=6cm]{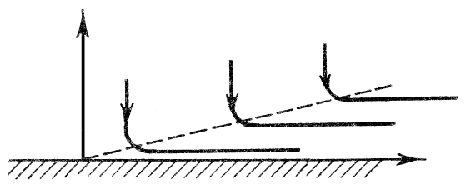}
\end{center}
\begin{center}
Fig.3.3. (\cite{loit_73}, p.541)
\end{center}
\subsection{The Existence and Uniqueness Theorems and Two-Sided Estimates
for Strictly Increasing Solutions}
\subsubsection{The case \,$m\in\{1/3; 1/2; \infty\}$: previous remarks on the exact solutions}
First, we need to summarize some results of previous Subsection 3.2.
\begin{proposition}\label{p_3}.
For fixed \,$a>0$ and \,$m\in\{1/3; 1/2; \infty\}$, singular nonlinear BVP  \eqref{eq2_phi}--\eqref{eq2_0}
has a unique solution \,$\Phi_m(\tau, a)$; each from these solutions is a strictly increasing on \,$\mathbb{R}_-$
function:

(i) solution  \,$\Phi_{1/3}(\tau, a)$ is represented by exact implicit formula \eqref{tau1/3},
where the data at the point \,$\tau=0$ and the parameter \,$d=d_{1/3}(a)$ in the Laypunov series
are defined by \eqref{phi1/3_0}, \eqref{dphi1/3_0} and \eqref{d1/3} respectively;
this solution has the inflection point \,$\tau=\tau_{\rm in}<0$ defined by the relation \eqref{phi1/3_in};

(ii) solution \,$\Phi_{1/2}(\tau, a)$ is represented by exact formula \eqref{phi05}, where the data at
the point \,$\tau=0$ and the parameter \,$d=d_{1/2}(a)$ in the Laypunov series  are defined by
\eqref{phi05_0}; it is a convex  function on \,$\mathbb{R}_-$;

(iii) solution \,$\Phi_{\infty}(\tau, a)$ is represented by exact formula \eqref{phi_inf_d},
where the data at the point \,$\tau=0$ and the parameter \,$d=d_{\infty}(a)$ in the Laypunov series
are defined by \eqref{phiinf_0}; it is a convex function on $\mathbb{R}_-$.
\end{proposition}

\begin{corollary}\label{c_6}.
For any fixed \,$a>0$ and \,$m\in \{1/3; 1/2; \infty \}$, singular nonlinear IBVP \eqref{eq2_phi}--\eqref{eq2_0}, \eqref{con3f_phi} defined on the entire real line has no solutions, for any fixed \,$b\neq 0$. Moreover,
for \,$m=1/3$ the exact solution  doesn't exist globally on \,$\mathbb{R}$, it blows up at finite point \,$\tau =\tau _p=2\pi\sqrt{3}/(3a)>0$.
\end{corollary}

The graphs of the above exact solutions to singular nonlinear BVP  \eqref{eq2_phi}--\eqref{eq2_0} continued to the right are demonstrated on Fig.1, where \,$a=1$.
\begin{center}
\includegraphics[width = 10cm,height=6cm]{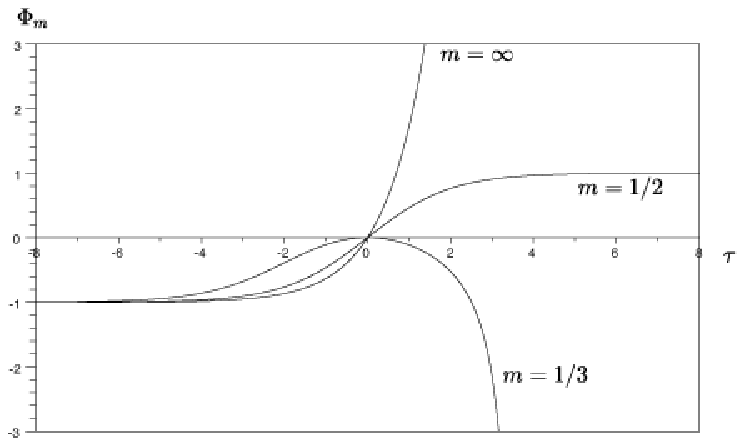}
\end{center}
\begin{center}
Fig.1
\end{center}

These solutions don't satisfy original singular nonlinear IBVP \eqref{eq2_phi}--\eqref{eq2_0}, \eqref{con3f_phi},
but they define the existence regions for the solutions of ODE \eqref{eq2_phi} with different behavior of patterns
on \,$\mathbb{R}_+$ (values \,$m\in\{1/3; 1/2; \infty\}$ are critical for \,$m$ as the bifurcation parameter).
\subsubsection{The case \,$m \ge 1/2$}
For fixed \,$m\ge 1/2$ and any  \,$d>0$, by virtue of \eqref{intdif1} and \eqref{intdif2},
a solution \,$\Phi_m(\tau,a,d)$ from series \eqref{l_ser}, \eqref{hl} strictly increases and is convex
at least while it remains non-positive (see Corollary \ref{c_5}). Then, by continuity, there is \,$d=d(m, a)>0$
such that condition \eqref{eq2_0} is satisfied (a variation in \,$d$ is equivalent to a shift in \,$\tau$).
\begin{theorem}\label{t_1}.
For any fixed  \,$m \ge 1/2$  and \,$a>0$, the singular nonlinear BVP \eqref{eq2_phi}--\eqref{eq2_0},
defined on \,$\mathbb{R}_-$, has a unique solution  \,$\Phi _m(\tau, a)$;  it is convex  strictly
increasing  function, belonging to the family \eqref{l_ser}, \eqref{hl} with some \,$d=d(m,a)>0$
and satisfying relations \eqref{intdif1}--\eqref{intd02}. Moreover, the following two-sided estimates
are valid:
\begin{equation}\label{rest}
\Phi_\infty(\tau,a)\le\Phi_m(\tau,a)\le\Phi_{1/2}(\tau,a) \qquad \forall m\ge 1/2,
\qquad \tau\in\mathbb {R}_-;
\end{equation}
here the "super-solution" \,$\Phi_{1/2}(\tau, a)$ and the "sub-solution" \,$\Phi_\infty (\tau,a)$
are defined by \eqref{phi05} and \eqref{phi_inf_d} respectively.
\end{theorem}

For \,$a=1$ and several values of \,$m$, the graphs of the solutions to the singular nonlinear BVP
\eqref{eq2_phi}--\eqref{eq2_0} are represented  on Fig.2, where the graphs of the super-solution ($m=1/2$) and the sub-solution ($m=\infty$) are indicated by dotted lines.
\begin{center}
\includegraphics[width = 10cm,height=6cm]{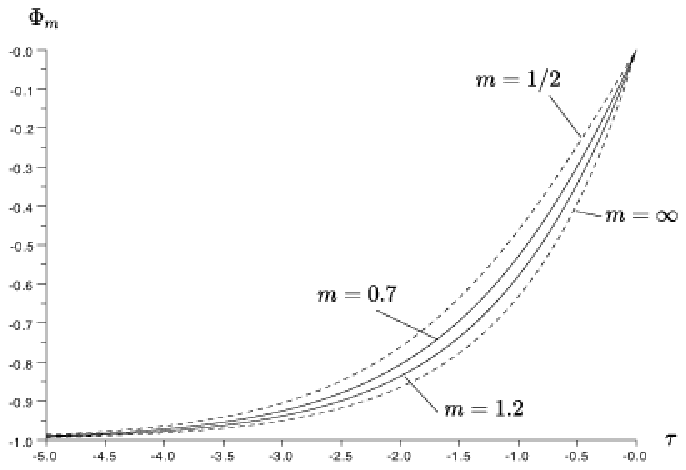}
\end{center}
\begin{center}
Fig.2
\end{center}
\subsubsection{The case \,$0< m< 1/2$}
In what follows, we take into account: \,1) for \,$m=1/2$
and \,$m=1/3$, the above analytical formulas for the exact solutions;
\,2) the behavior of solutions from family \eqref{l_ser}, \eqref{hl};
\,3) the continuity of solutions with respect to \,$m$, \,$1/3<m<1/2$ (see Remark \ref{r_7}).
\begin{theorem}\label{t_2}.
For any fixed  \,$m: \,1/3 \le m \le 1/2,$ and \,$a>0$, singular nonlinear BVP  \eqref{eq2_phi}--\eqref{eq2_0}
has a unique solution  \,$\Phi_m(\tau,a)$; it is a strictly increasing function, belonging  to the family \eqref{l_ser}, \eqref{hl}, with some  \,$d=d(a,m)>0$,  and satisfying estimates
\begin{equation}\label{restt}
\Phi_{1/3}(\tau,a)\le\Phi_m(\tau,a)\le\Phi_{1/2}(\tau,a) \qquad
\forall m: \,1/3 \le m \le 1/2, \qquad \tau\in \mathbb{R}_-.
\end{equation}

For fixed \,$m:\,1/3\le m<1/2$, the solution  has a point of inflection \,$\tau=\tau_{\rm in}\in\mathbb{R}_-$
defined by the relation
\begin{equation}\label{rel_in}
\Phi_{m}(\tau_{\rm in},a)\Phi^\prime_m(\tau_{\rm in},a)=-\int\limits_{-\infty}^{\tau_{\rm in}}
{\left[\Phi_{m}^\prime(s,a)\right]^2\,ds}, \quad \tau_{\rm in}\in\mathbb{R}_-,
\end{equation}
and doesn't exist globally on \,$\mathbb{R}$: it has a simple pole singularity at a finite point
\,${\tau=\tau_p(a,m)>0}$,  where  \,$\tau_p(a,1/3)=2\pi\sqrt{3}/(3a)$, and \,$\tau_p(a,m)>\tau _p(a,1/3)$
\,${\forall m:\,1/3<m<1/2}$.

For any fixed  \,$m:\,0<m<1/3$  and  \,$a>0$, singular nonlinear  BVP  \eqref{eq2_phi}--\eqref{eq2_0} has
no solutions.
\end{theorem}

The illustrations to this case are given on Figs.3, where \,$a=1$ (for the flows, see further).
\begin{center}
\includegraphics[width = 10cm,height=6cm]{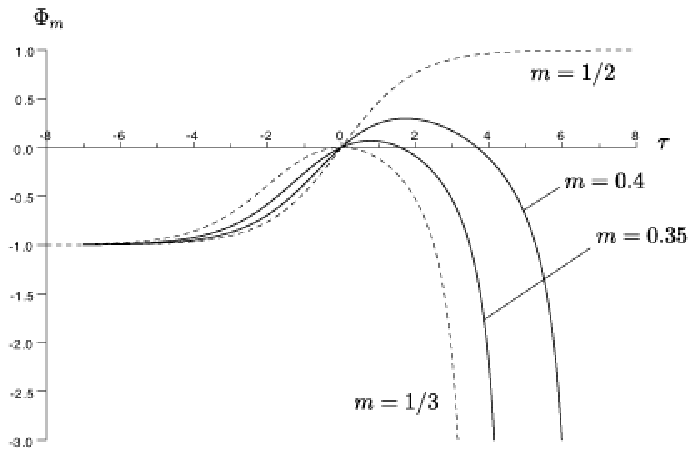}
\end{center}
\begin{center}
Fig.3a
\end{center}
\begin{center}
\includegraphics[width=7cm,height=5cm]{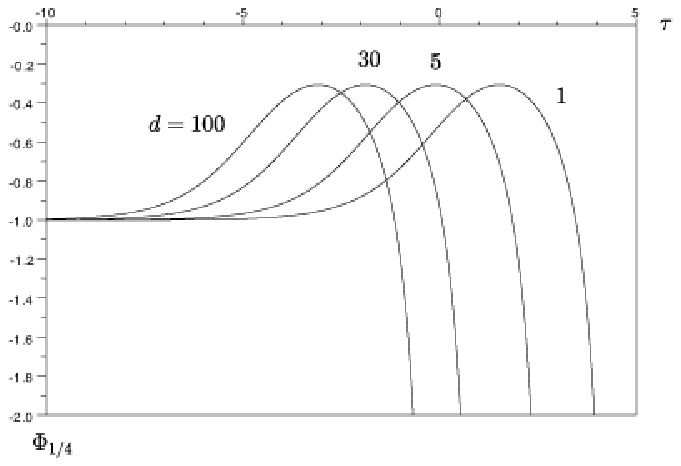}
\end{center}
\begin{center}
Fig.3b
\end{center}
\begin{center}
\includegraphics[width=7cm,height=5cm]{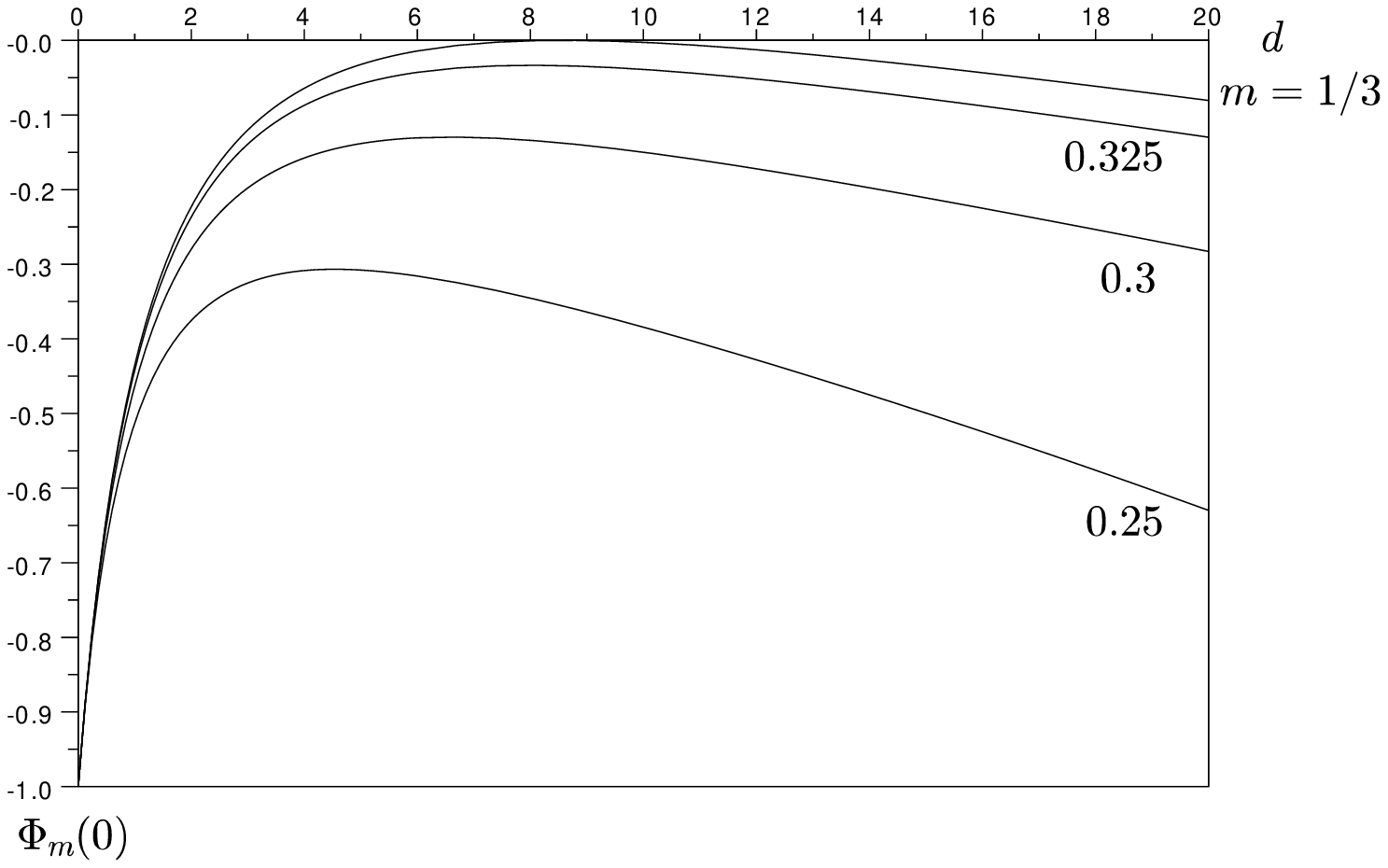}
\end{center}
\begin{center}
Fig.3c
\end{center}
\subsubsection{Calculation of the parameter \,$d=d_m(a)$ \,($m\ge 1/3$) in the Lyapunov series,
for a solution to BVP  \eqref{eq2_phi}--\eqref{eq2_0}}
Let us recall the exact values for the parameters \,$d=d_m(a)$ in the Lyapunov series, corresponding
to the exact solutions \,$\Phi_m(\tau, a)$ of BVP  \eqref{eq2_phi}--\eqref{eq2_0}:
\begin{equation*}
d_\infty(a)=a, \quad d_{1/2}(a)=2a, \quad d_{1/3}(a)=2\sqrt{3} a \exp{\left(\sqrt{3} \pi/6\right)}
\approx8.579306\,a.
\end{equation*}

Moreover, the following formulas are valid due to the scaling
transformations (see \eqref{phim_dm}):
\begin{equation*}
 d_m(a)=a\,d_m(1), \qquad \Phi_m(\tau, a)=a\,\Phi_m(a\tau, 1).
\end{equation*}

The values \,$d_m(1)$ are indicated in Table 2 and the graph \,$d=d_m(1)$ is represented on Fig.4. Here, in particular,
$d^\prime_m (1)\to - \infty$ as  $m\to 1/3+0$, where derivative is taken on $m$.
\begin{center}
Table 2 \,($\tau_{-\infty}=-T=-7$)

\begin{tabular}[t]{|l|l|}
\hline
$m$ & $d_m(1)$ \\
\hline
$1/3$ & $8.579306$   \\
$0.33334$ & $8.4470$   \\
$0.3334$ & $8.1710$   \\
$0.334$ & $7.3800$   \\
$0.335$ & $6.7938$   \\
$0.34$ & $5.5131$   \\
$0.35$ & $4.4467$   \\
$0.36$ & $3.8739$   \\
$0.37$ & $3.4922$   \\
$0.38$ & $3.2125$   \\
$0.39$ & $2.9959$   \\
$0.40$ & $2.8218$   \\
$0.41$ & $2.6781$   \\
$0.42$ & $2.5571$   \\
$0.43$ & $2.4534$   \\
$0.44$ & $2.3635$   \\
\hline
\end{tabular}
\begin{tabular}[t]{|l|l|}
\hline
$m$ & $d_m(1)$ \\
\hline
$0.45$ & $2.2846$   \\
$0.46$ & $2.2148$   \\
$0.47$ & $2.1525$   \\
$0.48$ & $2.0965$   \\
$0.49$ & $2.0460$   \\
$0.50$ & $2.0000$   \\
$0.5001$ & $1.9996$   \\
$0.501$ & $1.9956$   \\
$0.51$ & $1.9580$   \\
$0.52$ & $1.9196$   \\
$0.53$ & $1.8841$   \\
$0.54$ & $1.8514$   \\
$0.55$ & $1.8211$   \\
$0.56$ & $1.7929$   \\
$0.57$ & $1.7666$   \\
$0.58$ & $1.7421$   \\
\hline
\end{tabular}
\begin{tabular}[t]{|l|l|}
\hline
$m$ & $d_m(1)$ \\
\hline
$0.59$ & $1.7191$   \\
$0.60$ & $1.6975$   \\
$0.70$ & $1.5370$   \\
$0.80$ & $1.4370$   \\
$0.90$ & $1.3686$   \\
$0.99$ & $1.3232$   \\
$0.999$ & $1.3192$   \\
$1.00$ & $1.3188$   \\
$1.01$ & $1.3146$   \\
$1.02$ & $1.3104$   \\
$1.04$ & $1.3025$   \\
$1.06$ & $1.2949$   \\
$1.08$ & $1.2877$   \\
$1.10$ & $1.2809$   \\
$1.20$ & $1.2511$     \\
$1.50$ & $1.1904$   \\
\hline
\end{tabular}
\begin{tabular}[t]{|l|l|}
\hline
$m$ & $d_m(1)$ \\
\hline
$1.70$ & $1.1641$      \\
$1.90$ & $1.1441$       \\
$2.00$ & $1.1358$   \\
$2.50$ & $1.1056$   \\
$3.00$ & $1.0864$    \\
$3.50$ & $1.0731$   \\
$4.00$ & $1.0633$   \\
$4.50$ & $1.0558$   \\
$5.00$ & $1.0500$   \\
$10.00$ & $1.0243$   \\
$15.00$ & $1.0161$   \\
$20.00$ & $1.0120$   \\
$25.00$ & $1.0096$   \\
$50.00$ & $1.0048$   \\
$100.00$ & $1.0024$   \\
$\infty$ & $1.0000$  \\
\hline
\end{tabular}
\end{center}
\begin{center}
\includegraphics[width = 10cm,height=6cm]{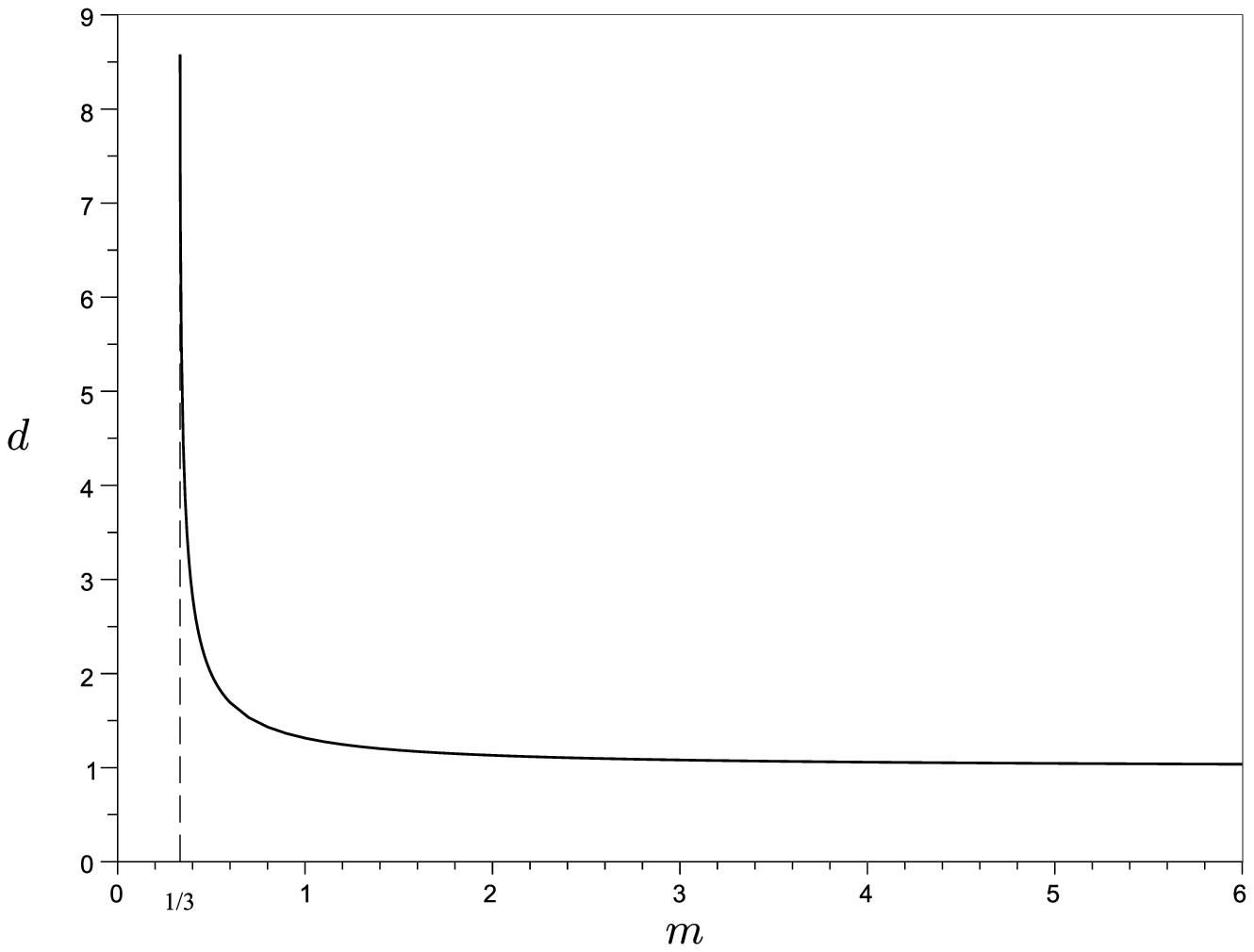}
\end{center}
\begin{center}
Fig.4 ($a=1$)
\end{center}
\begin{remark}\label{r_8}.
From singular nonlinear BVP \eqref{eq2_phi}--\eqref{eq2_0}, using Corollary \ref{c_1}, we obtain
the equivalent nonsingular BVP on the interval \,$[-T,0]$. After linearizing the obtained nonsingular BVP
about stationary solution \eqref{s_p} and taking into account Remark \ref{r_6}, we have the independent
of \,$m$ linear BVP of the form
\begin{equation*}
\Phi^{\prime\prime\prime} - a\Phi ^{\prime\prime} = 0,
\qquad -T\le \tau \le 0,
\end{equation*}
\begin{equation*}
\Phi (-T) + a - \Phi^{\prime\prime} (-T)/a^2 = 0, \quad
\Phi^\prime (-T) - \Phi^{\prime\prime} (-T)/a = 0,
\end{equation*}
\begin{equation*}
\Phi (0)=0.
\end{equation*}
This  BVP  has exact solution  \eqref{phi_inf_d} which coincides with the exact solution to the
original singular nonlinear BVP \eqref{eq2_phi}--\eqref{eq2_0} as \,$m\to \infty$.

As a result, applying a quasilinearization method for solving the nonsingular nonlinear BVP on
the interval \,$[-T,0]$, we obtain function \eqref{phi_inf_d} as an initial approximation
for any  \,$m\not=0$. Computations show that for different fixed \,$m\ge 1/3$,  and
particularly for \,$m\ge 1/2$,  the solutions to the singular nonlinear BVP
\eqref{eq2_phi}--\eqref{eq2_0} are fairly similar indeed (see Figs.1,2).
\end{remark}
\subsubsection{Additional remarks: numerical results for the flows corresponding \,$m:\,1/3<m<1/2$\,
("singular flows")}
The numerical results we present here demonstrates, e.g., the formal pass from \,$m=1/3$ to \,$m=1/2$ case.
The computations were made up to a certain neighborhood of the pole depended on \,$m$.
We don't know if any physical interpretation of these illustrations is possible.

Here, analogously to the case \,$m=1/3$, we consider BVP \eqref{eq2_phi}--\eqref{eq2_0}
with the change of variables \,$\tau$ by \,$-\tau$ and \,$\Phi$ by $-\Phi$ (the same is
done on the graphs).

Figs.5--7 show the results of the computations with \,$a=\nu=1$.
\begin{center}
\includegraphics[width=7cm,height=5cm]{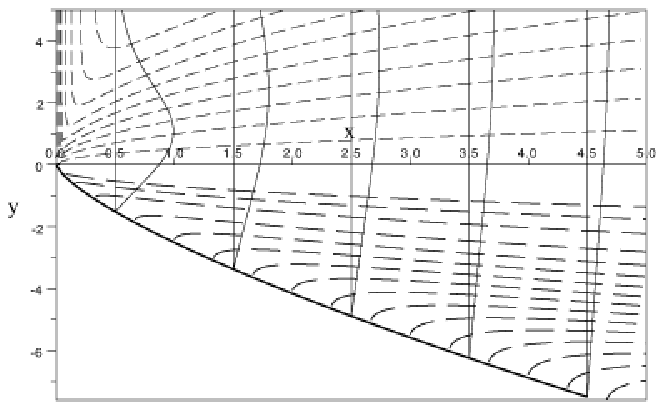}
\end{center}
\begin{center}
Fig.5a
\end{center}
\begin{center}
\includegraphics[width=7cm,height=5cm]{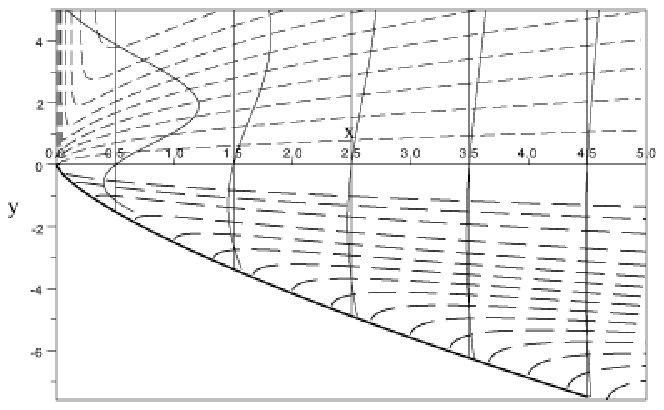}
\end{center}
\begin{center}
Fig.5b
\end{center}
The curves \,$\psi_{3/8}(x,y)=$const (dotted lines), the profiles of the horizontal
velocity component \,$\widetilde{u}_{3/8}(y)=u_{3/8}(x, y)|_{x={\rm const}}$ (Fig.5a)
and the profiles of the vertical velocity component
\,$\widetilde{v}_{3/8}(y)=v_{3/8}(x, y)|_{x={\rm const}}$ (Fig.5b).
\begin{center}
\includegraphics[width=7cm,height=5cm]{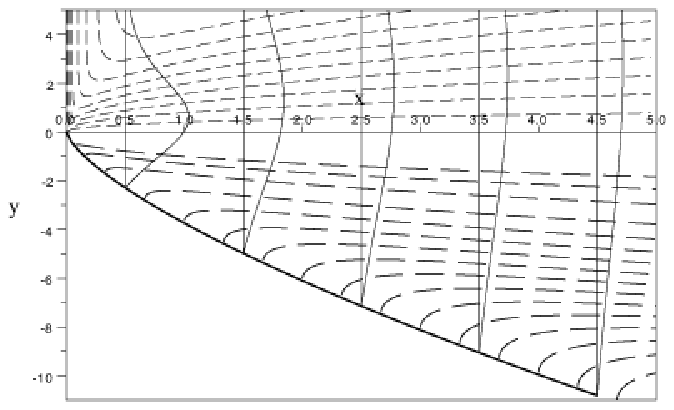}
\end{center}
\begin{center}
Fig.6a
\end{center}
\begin{center}
\includegraphics[width=7cm,height=5cm]{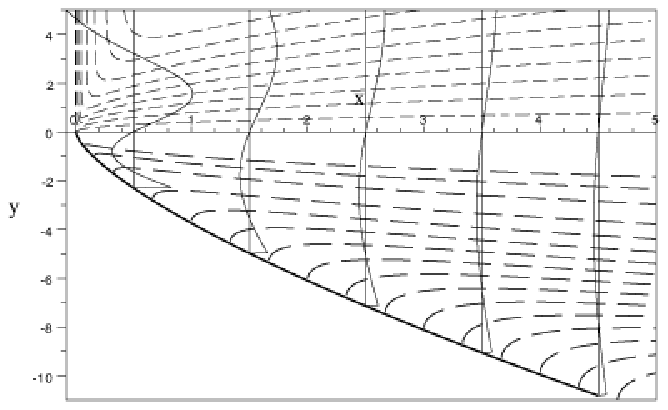}
 \end{center}
\begin{center}
Fig.6b
\end{center}
The curves \,$\psi_{5/12}(x,y)=$const (dotted lines), the profiles of the horizontal
velocity component \,$\widetilde{u}_{5/12}(y)=u_{5/12}(x, y)|_{x={\rm const}}$ (Fig.6a) and
\,the profiles of the vertical velocity component
\,$\widetilde{v}_{5/12}(y)=v_{5/12}(x, y)|_{x={\rm const}}$ (Fig.6b).
\begin{center}
\includegraphics[width=7cm,height=5cm]{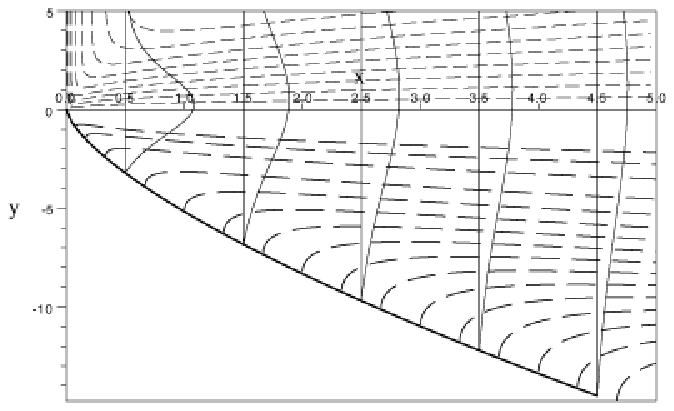}
\end{center}
\begin{center}
Fig.7a
\end{center}
\begin{center}
\includegraphics[width=7cm,height=5cm]{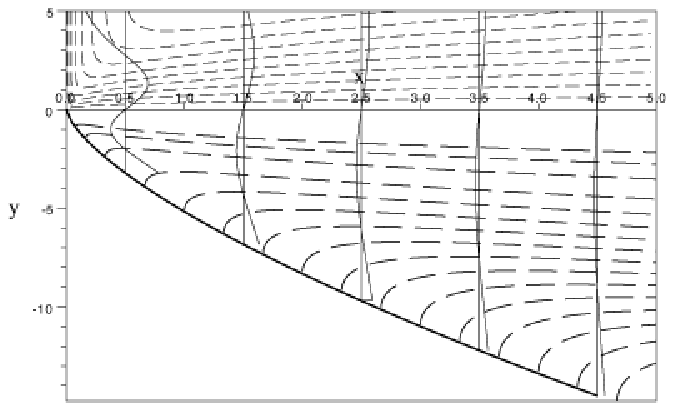}
\end{center}
\begin{center}
Fig.7b
\end{center}
The curves \,$\psi_{11/24}(x,y)=$const (dotted lines), the profiles of the horizontal
velocity component \,$\widetilde{u}_{11/24}(y)=u_{11/24}(x, y)|_{x={\rm const}}$ (Fig.7a)
and the profiles of the vertical velocity component with the scale factor
\,$\widetilde{v}_{11/24}(y)=0.5 v_{11/24}(x, y)|_{x={\rm const}}$ (Fig.7b).

The lowermost curve (solid line) on each of the figures shows approximately the line
\,${y(x)=-\tau_p(m)\,x^{1/(m+1)}\,\sqrt{(m+1)/m}}$\, ($x\ge 0$), where \,$\tau_p(m)$ is the
pole point for the solution \,$\Phi_m(\tau)$ of BVP \eqref{eq2_phi}--\eqref{eq2_0} extended to the right
(note that the exact position is known for pole point when \,$m=1/3$ (see \eqref{tau_1/3})).

Slightly above that line the other singular line $y(x)=-\tau_{\rm max}(m)x^{1/(m+1)}\times$ $\times\sqrt{(m+1)/m}$ \,($x\ge 0$)
goes, where \,$\tau_{\rm max}(m)$: \,$\Phi^\prime_m(\tau_{\rm max})=0$ (for \,$m=1/3$ we have \,$\tau_{\rm max}(1/3)=0$).
In the points of such line, the horizontal velocity component is equal to zero (i.e., there is a "curvilinear  wall"\,;
such line is absent on figures in view of difficulties to calculate and construct it graphically).

Once again, we should note that we do not aware if these calculated data has a physical meaning.
Anyway, it makes possible to observe the transfer from the "near-wall flow" \,$m=1/3$ to the
"flooded jet" \,$m=1/2$.
\section{Existence, Two-Sided Estimates and Behavior of the Solutions to
the Original Singular Nonlinear IBVP}
\subsection{The Necessary and Sufficient Conditions for the Solution of the Accompanying
Singular Nonlinear BVP to Be Continued on the Positive Half-Line}
We analyze the behavior of solutions to the singular nonlinear BVP
\eqref{eq2_phi}--\eqref{eq2_0} when they are indefinitely extended rightward, which
is possible if \,$m\ge 1/2$. There are two limiting cases
corresponding to \,$m=1/2$ and \,$m\to\infty$:
\begin{equation}\label{m_1/2_inf}
\lim _{\tau \to\infty}\Phi _{1/2}(\tau,a)=\lim _{\tau \to \infty}[\Phi
_{\infty} (\tau,a)/\exp (a\tau)] = a.
\end{equation}
Thus we consider the case \,$1/2 < m < \infty$ and look for the needed
solutions in the form
\begin{equation}\label{Y}
\Phi (\tau) = \tau ^m [b + Y(\tau)], \quad \tau >0,
\end{equation}
We use the change of variables:
\begin{equation}\label{xi_Y}
\xi = \tau ^{m+1}/(m+1), \qquad  Y(\tau (\xi)) = v (\xi),
\end{equation}
where \,$b$ is a parameter ($b \not = 0$) and
\,$\lim_{\tau\to\infty}{Y(\tau)} = \lim_{\xi\to\infty}{v(\xi)} = 0$.

For \,$v (\xi)$, we  obtain the singular nonlinear CP at infinity
(the dotted letters denote derivatives with respect to \,$\xi$):
\begin{equation*}
\dddot{v} + \left [\frac {6m}{(m+1)\xi} + b \right ]\ddot{v} + \left [\frac {(7m-4)m}{(m+1)^2\xi ^2} +
\frac {(m+2)b}{(m+1)\xi }\right ]\dot{v} +
\frac{m(m-1)(m-2)}{(m+1)^3 \xi ^3} \,v  +
\end{equation*}
\begin{equation}\label{xi_Y}
 + \frac {m(m-1)(m-2)b}{(m+1)^3\xi ^3} = -v \ddot{v} -
\frac{m+2}{(m+1)\xi}\, v\dot{v} + \frac {m-1}{m}\,\dot{v}^2, \quad 0< \xi < \infty,
\end{equation}
\begin{equation}\label{xi_inf}
\lim_{\xi\to\infty}{v(\xi)} = \lim_{\xi\to\infty}\dot{v}(\xi)
= \lim_{\xi\to\infty}\ddot{v}(\xi) = 0.
\end{equation}

Nonlinear ODE \eqref{xi_Y} has an irregular singularity  of a rank \,$1$ as \,$\xi\to\infty$.
Next two propositions follow from the general theory of the ODE systems with irregular singular points,
including certain classes of nonlinear ODEs (see, e.g., \cite{was}).

In the first place, let us remark that the following formal series satisfies Eq.\eqref{xi_Y} with
\,$b\ne 0$, \,$m>0$:
\begin{equation}\label{xi_par}
v_{\rm form} (\xi, b)= \sum_{k=1}^{\infty} v _k/\xi ^k,
\end{equation}
where
\begin{equation}\label{xi_1}
v_1=-(m-1)(m-2)/(m+1)^2,
\end{equation}
\begin{equation*}
v_{k+1}=\Bigl\{\Bigl[k(k+1)(k+2) - \frac {6m}{m+1}k(k+1) +\frac{(7m-4)m}{(m+1)^2} k-
\end{equation*}
\begin{equation*}
-\frac {m(m-1)(m-2)}{(m+1)^3}\Bigr] v_k - \sum_{l=1}^kl\left [l+1 - \frac {m+2}{m+1} -
\frac {m-1}{m} (k-l+1)\right ]v_lv_{k-l+1} \Bigr \}\Big /
\end{equation*}
\begin{equation}\label{xi_k+1}
\Big / \left [b(k+1)\left (k+2 - \frac {m+2}{m+1}\right )\right ], \qquad k=1, 2, \ldots,
\end{equation}
\begin{proposition}\label{p_4} For any fixed \,$b\ne 0$ and \,$m>0$, the
singular nonlinear CP \eqref{xi_Y}, \eqref{xi_inf} has a particular solution
\,$v_{\rm par} (\xi)$ which has the series \eqref{xi_par}--\eqref{xi_k+1} as its
asymptotic expansion for large \,$\xi$.
\end{proposition}
\begin{proposition}\label{p_5} For every fixed \,$m>0$, nonlinear  ODE \eqref{eq2_phi}
has a three-parameter family of solutions \,$\Phi_m (\tau + \tau _s,b,D)$ that can be represented
for large positive \,$\tau $ in the principal approximation as
\begin{equation*}
\Phi_m(\tau+\tau_s,b,D)=(\tau+\tau_s)^m \Bigl\{b+v_{\rm par}
\Bigl((\tau+\tau_s)^{m+1}/(m+1),\,b\Bigr)+
\end{equation*}
\begin{equation}\label{phim_inf}
+D(\tau+\tau_s)^\gamma\exp{\Bigl(-b(\tau+\tau_s)^{m+1}/(m+1)\Bigr)}
\Bigl[1+o(1)\Bigr] \Bigr\}, \qquad \tau\to\infty,
\end{equation}
where \,$\tau_s$, \,$b$, and \,$D$ are parameters, \,$b>0$,
\begin{equation}\label{gam}
\gamma=(-4m^2-6m+4)/(m+1),
\end{equation}
and \,$v_{\rm par}(\xi,b)$ is defined by Proposition \ref{p_4}.
\end{proposition}
\begin{remark}\label{r_9} For \,$m=1$ and \,$m=2$, it follows from \eqref{xi_par}--\eqref{xi_k+1}
that $v_{\rm par} (\xi)\equiv 0$. Then from \eqref{phim_inf}, \eqref{gam} we have, for large \,$\tau>0$,
\begin{equation}\label{phi1_inf}
\Phi_1(\tau+\tau_s,b,D)=(\tau+\tau_s)\Bigl\{b+D(\tau +\tau_s)^{-3}\exp {\Bigl(-b(\tau +\tau_s)^{2}/2\Bigr)}\Bigl[1+o(1)\Bigr] \Bigr\},
\end{equation}
\begin{equation}\label{phi2_inf}
\Phi_2(\tau+\tau_s,b,D)=(\tau+\tau_s)^2\Bigl\{b+D(\tau +\tau _s)^{-8}\exp{\Bigl(-b (\tau+\tau_s)^{3}/3\Bigr)}\Bigl[1+o(1)\Bigr] \Bigr\},
\end{equation}
and the two-parameter set of the exact solutions \,$\Phi_1(\tau+\tau_s,b)=b(\tau+\tau_s)$ {($\Phi_2(\tau+\tau_s,b)=b(\tau+\tau_s)^2$)} belongs to family \eqref{phi1_inf}  (resp., to family \eqref{phi2_inf})
at \,$D=0$ and \,$b\in\mathbb{R}$.
\end{remark}

For a solution \,$\Phi (\tau)$ to ODE \eqref{eq2_phi} that satisfies condition \eqref{eq2_0}, we
infer  that
\begin{equation}\label{dphi_rel}
[\Phi^\prime(\tau) + \Phi^2 (\tau)/2]^\prime=\Phi^{\prime\prime}(0) +
[(2m-1)/m] \int\limits_0^\tau [\Phi^\prime (t)]^2\,dt,
\end{equation}
and that \,$\Phi^\prime (0)>0$ and \,$\Phi ^{\prime\prime}(0)>0$ for
\,$m>1/2$ (see \eqref{intd01}, \eqref{intd02}).  It then follows from \eqref{dphi_rel} that, for
\,$m>1/2$, the value \,$\Phi ^\prime (\tau) + \Phi ^2 (\tau)/2$
increases on \,$\tau$ for \,$\tau>0$ and the following assertion holds.
\begin{proposition}\label{p_6} For any fixed $a>0$ and $m>1/2$, the solution
$\Phi_m (\tau, a)$ to the singular nonlinear BVP \eqref{eq2_phi}--\eqref{eq2_0} can
be indefinitely extended rightward and, for large $\tau >0 $, has
a representation of form \eqref{phim_inf} with some $b=b(m,a)>0$, $D=D(m,a)$,
and $\tau _s=\tau _s(m,a)$.
\end{proposition}
\begin{corollary}\label{c_7} For any fixed \,$m>1/2$ and \,$b>0$, singular
nonlinear  IBVP \eqref{eq2_phi}--\eqref{eq2_0}, \eqref{con3f_phi} defined on the entire real
line has a unique solution \,$\Phi_m (\tau, a, b)$; it belongs to
the Lyapunov series family \eqref{l_ser}, \eqref{hl} for some \,$a=a(m,b)>0$ and
\,$d=d(m,b)>0$.
\end{corollary}

Note once again, that, solving the  singular nonlinear BVP \eqref{eq2_phi}--\eqref{eq2_0}, we obtain
the Cauchy data  \,$\Phi(0)=0$,  \,$\Phi^\prime (0,a)>0$, and  \,$\Phi^ {\prime\prime}(0,a)>0$
which determine three parameters in \eqref{phim_inf}, including \,$b=b_m(a)$.  A remarkable
circumstance is that it suffices to solve this BVP  for  \,$a=1$  and to find the corresponding value
of  \,$b=b_m(1)$  because formulas \eqref{b_m}, \eqref{phim_dm} hold.

Fig.8 shows the solutions to the singular nonlinear BVP \eqref{eq2_phi}--\eqref{eq2_0}
for  \,$a=1$ extended to positive \,$\tau$ for different values of \,$m$.
\begin{center}
\includegraphics[width = 10cm,height=6cm]{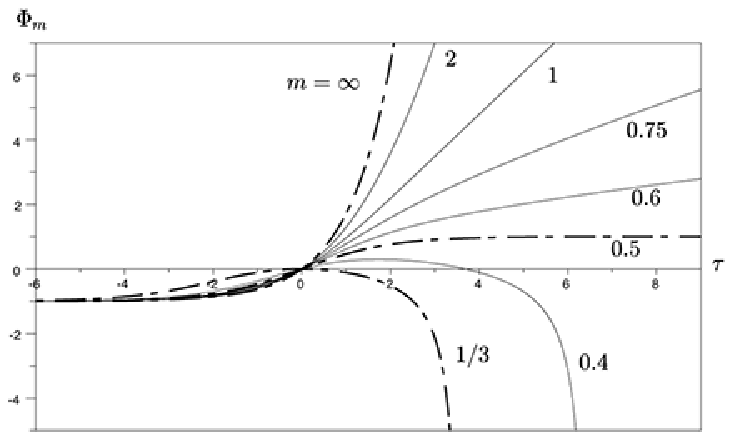}
\end{center}
\begin{center}
Fig.8
\end{center}

Let us recall once more that, due to the scaling transformations,
the following formulas are valid (see \eqref{b_m}, \eqref{phim_dm}):
\,$b_m(a)=a^{m+1} b_m(1)$, \,$\Phi_m(\tau, a)=a\,\Phi_m(a\tau, 1)$.

The values \,$b_m(1)$, for various \,$m>1/2$, are given in Table 3 and the graph \,$b=b_m (1)$
is represented on Fig.9. Here,  \,$b_m(1)\to 0$  as  \,$m\to 1/2+0$  or \,$m\to\infty $ (for these \,$m$,
the asymptotic behavior of the solution to the singular nonlinear BVP \eqref{eq2_phi}--\eqref{eq2_0}
is changing for large \,$\tau$), and \,$b^\prime_m(1)\to\infty$ as \,$m\to 1/2+0$, where derivative is
taken on \,$m$ \,($m=1/2$ is the branch point of the function \,$b_m(1)$).

\begin{center}

Table 3 $\left(\tau_{-\infty} = -T = -7\right)$

\begin{tabular}[t]{|l|l|}
\hline
$m$ & $b_m(1)$ \\
\hline
$0.5001$ & $0.062134$   \\
$0.501$ & $0.092668$   \\
$0.505$ & $0.17106$   \\
$0.51$ & $0.23382$   \\
$0.52$ & $0.32389$   \\
$0.53$ & $0.39361$   \\
$0.54$ & $0.45276$   \\
$0.55$ & $0.50516$   \\
$0.56$ & $0.55272$   \\
$0.57$ & $0.59657$   \\
$0.58$ & $0.6374$   \\
$0.59$ & $0.67571$   \\
$0.60$ & $0.71184$   \\
$0.70$ & $0.98975$   \\
\hline
\end{tabular}
\begin{tabular}[t]{|l|l|}
\hline
$m$ & $b_m(1)$ \\
\hline
$0.80$ & $1.1634$   \\
$0.85$ & $1.221$   \\
$0.90$ & $1.2621$   \\
$0.92$ & $1.2744$   \\
$0.94$ & $1.2845$   \\
$0.95$ & $1.2887$     \\
$0.96$ & $1.2925$   \\
$0.98$ & $1.2984$   \\
$0.99$ & $1.3007$   \\
$0.999$ & $1.3023$   \\
$1.00$ & $1.3025$   \\
$1.01$ & $1.3038$   \\
$1.02$ & $1.3047$   \\
$1.03$ & $1.3052$   \\
\hline
\end{tabular}
\begin{tabular}[t]{|l|l|}
\hline
$m$ & $b_m(1)$ \\
\hline
$1.04$ & $1.3053$   \\
$1.05$ & $1.3049$   \\
$1.06$ & $1.3042$   \\
$1.07$ & $1.3031$   \\
$1.08$ & $1.3016$   \\
$1.09$ & $1.2997$   \\
$1.10$ & $1.2975$   \\
$1.12$ & $1.2975$   \\
$1.14$ & $1.2854$   \\
$1.16$ & $1.2775$   \\
$1.18$ & $1.2684$   \\
$1.20$ & $1.2584$   \\
$1.35$ & $1.1559$   \\
$1.4$ & $1.1141$      \\
\hline
\end{tabular}
\begin{tabular}[t]{|l|l|}
\hline
$m$ & $b_m(1)$ \\
\hline
$1.50$ & $1.0236$   \\
$1.60$ & $0.92833$   \\
$1.70$ & $0.8324$   \\
$1.80$ & $0.73884$   \\
$1.90$ & $0.64985$   \\
$2.00$ & $0.56684$   \\
$2.50$ & $0.25817$   \\
$3.00$ & $0.10274$   \\
$3.50$ & $0.036882$   \\
$4.00$ & $0.012181$   \\
$4.50$ & $0.0122$   \\
$5.00$ & $0.0010872$   \\
$6.00$ & $7.832e-005$   \\
$\infty$ & 0.00000 \\
\hline
\end{tabular}
\end{center}
\begin{center}
\includegraphics[width = 10cm,height=6cm]{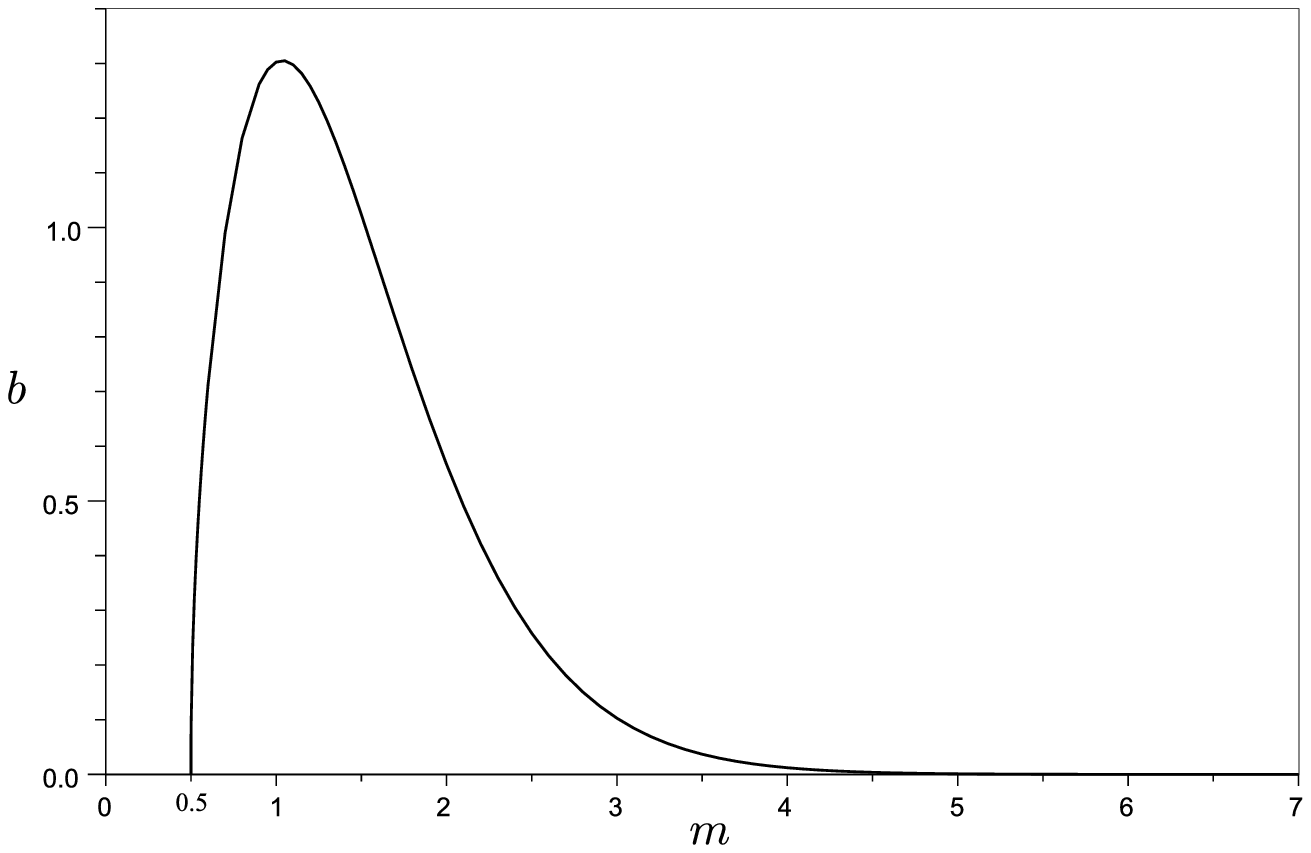}
\end{center}
\begin{center}
Fig.9 ($a=1$)
\end{center}
\subsection{Main Result for the Singular Nonlinear IBVP}
\begin{theorem}\label{t_3} For any fixed  \,$m>1/2$  and  \,$b>0$,
the singular nonlinear  IBVP \eqref{eq2_phi}--\eqref{eq2_0}, \eqref{con3f_phi}
defined on the entire real line has a unique solution \,$\Phi_m(\tau,a,b)$,
where \,$a=a(m,b)>0$, and the following assertions are valid:

(i) \,$\Phi_m(\tau,a,b)$ is convex  monotone increasing
function, belonging to the family \eqref{l_ser}, \eqref{hl} for some
\,$d=d(m,a,b)>0$ and satisfying the restrictions
\begin{equation}\label{subsup}
a[\exp{(a\tau)}-1]\le\Phi_m(\tau,a,b)\le a\tanh{(a\tau/2)}, \qquad -\infty<\tau\le 0,
\end{equation}
\begin{equation}\label{subsup1}
a\tanh{(a\tau/2)}<\Phi_m(\tau,a,b)<a[\exp{(a\tau)}-1], \qquad \tau>0;
\end{equation}

(ii) for large \,$\tau>0$, \,$\Phi_m(\tau,a,b)$ has the
representation of form \eqref{phim_inf} with certain \,$D=D(m,a,b)$ and
\,$\tau_s=\tau_s(m,a,b)$;

(iii) the solution \,$\Phi_m(\tau,a,b)$ may be obtained as
follows: fix \,$a=1$ and define the solution \,$\Phi_m(\tau,1)$ of singular
nonlinear BVP \eqref{eq2_phi}--\eqref{eq2_0} which exists, is unique
and belongs to the Lyapunov series family \eqref{l_ser}, \eqref{hl} with some
\,$d=d_m(1)$ (according to Theorem \ref{t_1}); being extended to the
right, this solution satisfies the limit condition
\begin{equation}\label{phi_m_inf}
\lim_{\tau\to\infty}\Phi_m (\tau,1)/\tau^{m}=b_m(1)>0;
\end{equation}
due to the scaling transformations, the needed solution is defined as
\begin{equation}\label{phi_ab}
\Phi_m(\tau,a,b)=a\Phi_m (a\tau,1), \qquad \tau\in\mathbb{R},
\end{equation}
where
\begin{equation}\label{adbm}
a=a(m, b)=[b/b_m(1)]^{1/(m+1)},\quad d=d_m(a(m, b))=a(m, b)d_m(1).
\end{equation}
\end{theorem}
\subsection{Numerical Results for Different Values of the Self-Similarity
Parameter $m$}
For all numerical results,  we set \,$a=1$ and \,$\nu=1$.
\subsubsection{The case \,$1/2<m<1$ \,($m\in\{11/20, 3/5, 4/5\}$)}
\begin{center}
\includegraphics[width = 10cm,height=6cm]{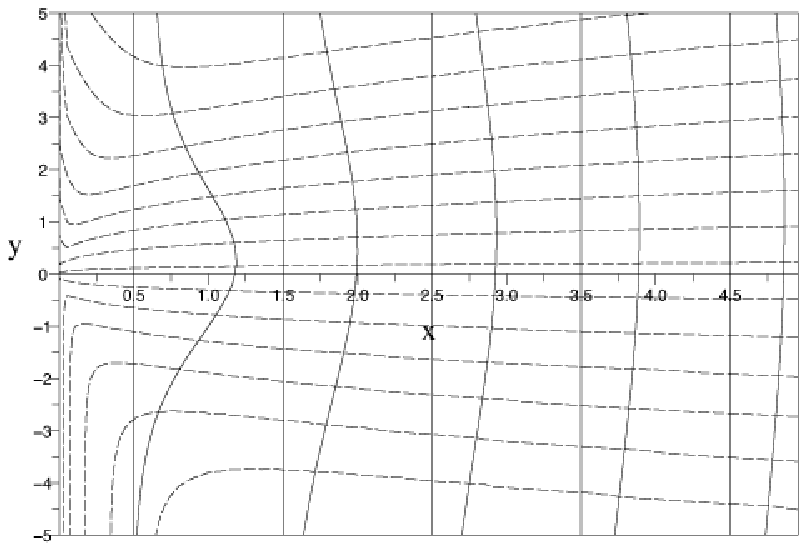}
\end{center}
\begin{center}
Fig.10a. The curves \,$\psi_{11/20}(x,y)=$const (dotted lines) and the profiles
of the horizontal velocity component \,$\widetilde u_{11/20}(y)=u_{11/20}(x, y)|_{x={\rm const}}$.
\end{center}
\begin{center}
\includegraphics[width = 10cm,height=6cm]{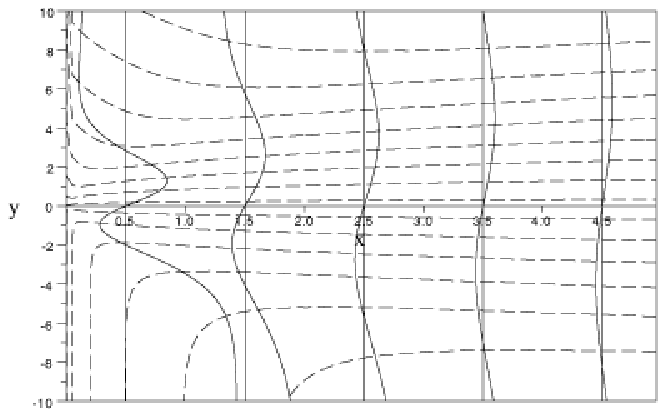}
\end{center}
\begin{center}
Fig.10b. The curves \,$\psi_{11/20}(x,y)=$const (dotted lines) and the profiles
of the vertical  velocity component  \,$\widetilde v_{11/20}(y)=v_{11/20}(x, y)|_{x={\rm const}}$.
\end{center}
\begin{center}
\includegraphics[width = 10cm,height=6cm]{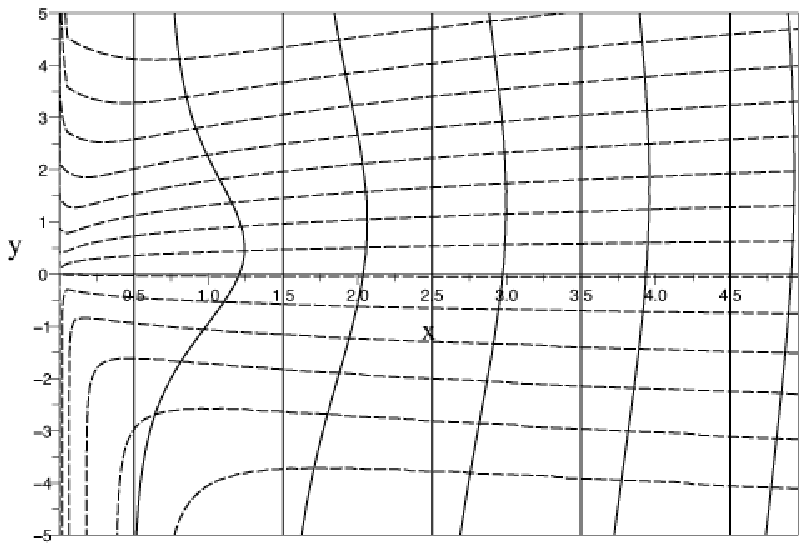}
\end{center}
\begin{center}
Fig.11a. The curves \,$\psi_{3/5}(x,y)=$const (dotted lines) and the profiles
of the horizontal velocity component \,$\widetilde u_{3/5}(y)=u_{3/5}(x, y)|_{x={\rm const}}$.
\end{center}
\begin{center}
\includegraphics[width = 10cm,height=6cm]{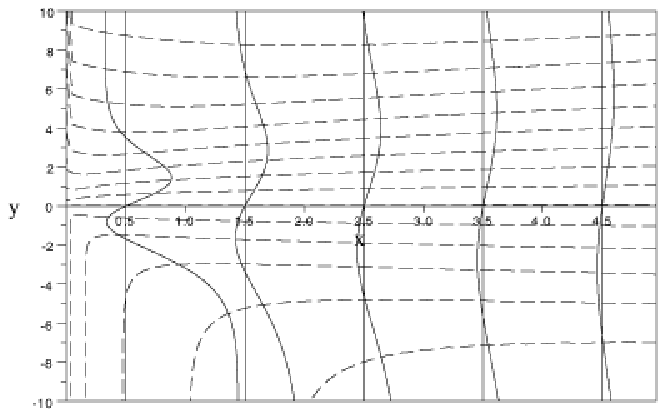}
\end{center}
\begin{center}
Fig.11b. The curves \,$\psi_{3/5}(x,y)=$const (dotted lines) and the profiles
of the vertical  velocity component \,$\widetilde  v_{3/5}(y)=v_{3/5}(x,y)|_{x={\rm const}}$.
\end{center}
\begin{center}
\includegraphics[width = 10cm,height=6cm]{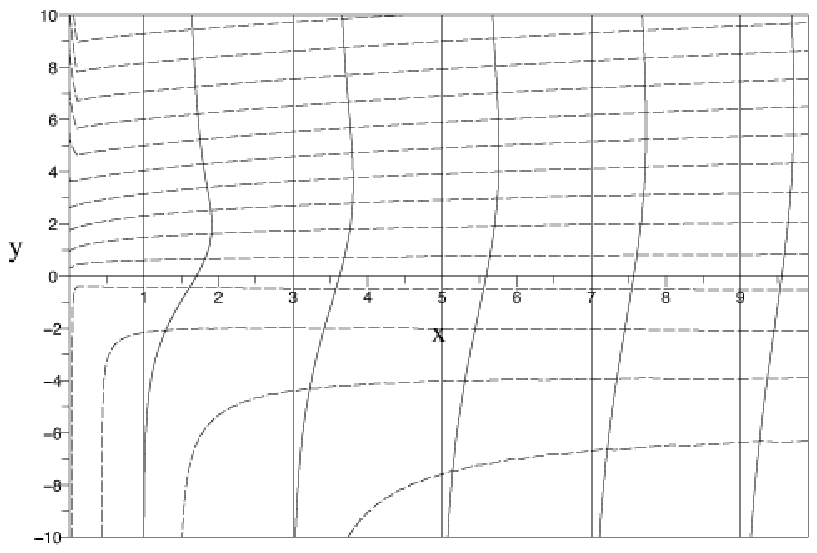}
\end{center}
\begin{center}
Fig.12a. The curves \,$\psi_{4/5}(x,y)=$const (dotted lines) and the profiles
of the horizontal velocity component \,$\widetilde u_{4/5}(y)=u_{4/5}(x, y)|_{x={\rm const}}$.
\end{center}
\begin{center}
\includegraphics[width = 10cm,height=6cm]{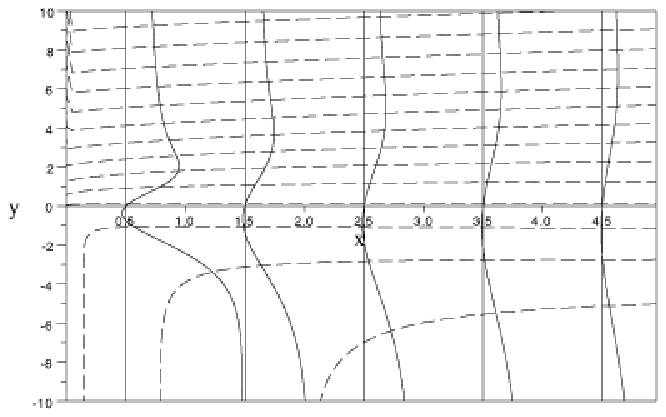}
\end{center}
\begin{center}
Fig.12b. The curves \,$\psi_{4/5}(x,y)=$const (dotted lines) and the profiles
of the vertical  velocity component \,$\widetilde v_{4/5}(y)=v_{4/5}(x,y)|_{x={\rm const}}$.
\end{center}
\subsubsection{The case \,$m=1$}
This case is known as a problem on a plane "semi-jet". More exactly, the following problem is
formulated in \cite{shlich}, pp.180-181, as a problem for laminar layer on an interface between two flows
(see also references to this problem in \cite{shlich}):
\begin{equation}\label{eqf_m1}
f^{\prime\prime\prime}(\eta) + f(\eta)\,f^{\prime\prime}(\eta)/2=0, \qquad -\infty<\eta<\infty,
\end{equation}
\begin{equation}\label{con1f_m1}
\lim_{\eta\to -\infty}f^{\prime}(\eta)=U_2/U_1=\lambda,
\end{equation}
\begin{equation}\label{con2f_m1}
f(0)=0,
\end{equation}
\begin{equation}\label{con3f_m1}
\lim_{\eta\to \infty}f^{\prime}(\eta)=1,
\end{equation}
where \,$\eta$ is the self-similar variable,
\begin{equation}\label{eta_m1}
\eta(x, y) = y\sqrt{U_1/(\nu x)},
\end{equation}
\,$U_1$ and \,$U_2$ are the constant velocities of the upper and lower
flows respectively (in our case, $U_2=0$, i.e., \,$\lambda=0$).

For the stream function \,$\psi (x, y, U_1)$ and velocity  \,$x$-component \,$u(\eta, U_1)$, we have
\begin{equation}\label{psi_fm1}
\psi(x, y, U_1)=\sqrt{\nu U_1 x}\, f(\eta), \qquad u(\eta, U_1)=U_1\, f^\prime (\eta).
\end{equation}

Putting  $U_1 = 1/2$, we obtain $\eta(x,y) = \tau(x,y)$,
\begin{equation}\label{phi_fm1}
\Phi_1(\tau)= f(\eta)/2, \quad  \widetilde u(\eta) = f^\prime(\eta) =
2 \Phi^\prime_1(\tau, a) = 2 u_1(\tau,a),
\end{equation}
where $a: b_1(a)=1/2$.

From our tables, we have $d_1(1) \approx 1.3188$,\, $b_1(1)
\approx 1.3025$. Then, using the relations $b_1(a)= b_1(1) a^2 = 1/2$,\,
$d_1(a)=d_1(1)a$,  we obtain $a\approx 0.61958$, $d(a)\approx
0.8171$.

On Fig.13a, there are the graphs of the function
\begin{equation}\label{ueta_m1}
\widetilde u(\eta) = u(\eta, U_1)/U_1 =  f^\prime(\eta),
\end{equation}
for the values \,$\lambda=0$ and \,$\lambda=0.5$.

On Fig.13b, there is the graph of \,$\widetilde u(\eta) = 2\Phi^\prime_1(\tau, a)$, to be compared with the
\,$\lambda = 0$ graph on Fig.13a. (In case \,$\lambda\not=0$ the singular problem is mathematically different.)
\begin{center}
\includegraphics[width = 7cm,height=6cm]{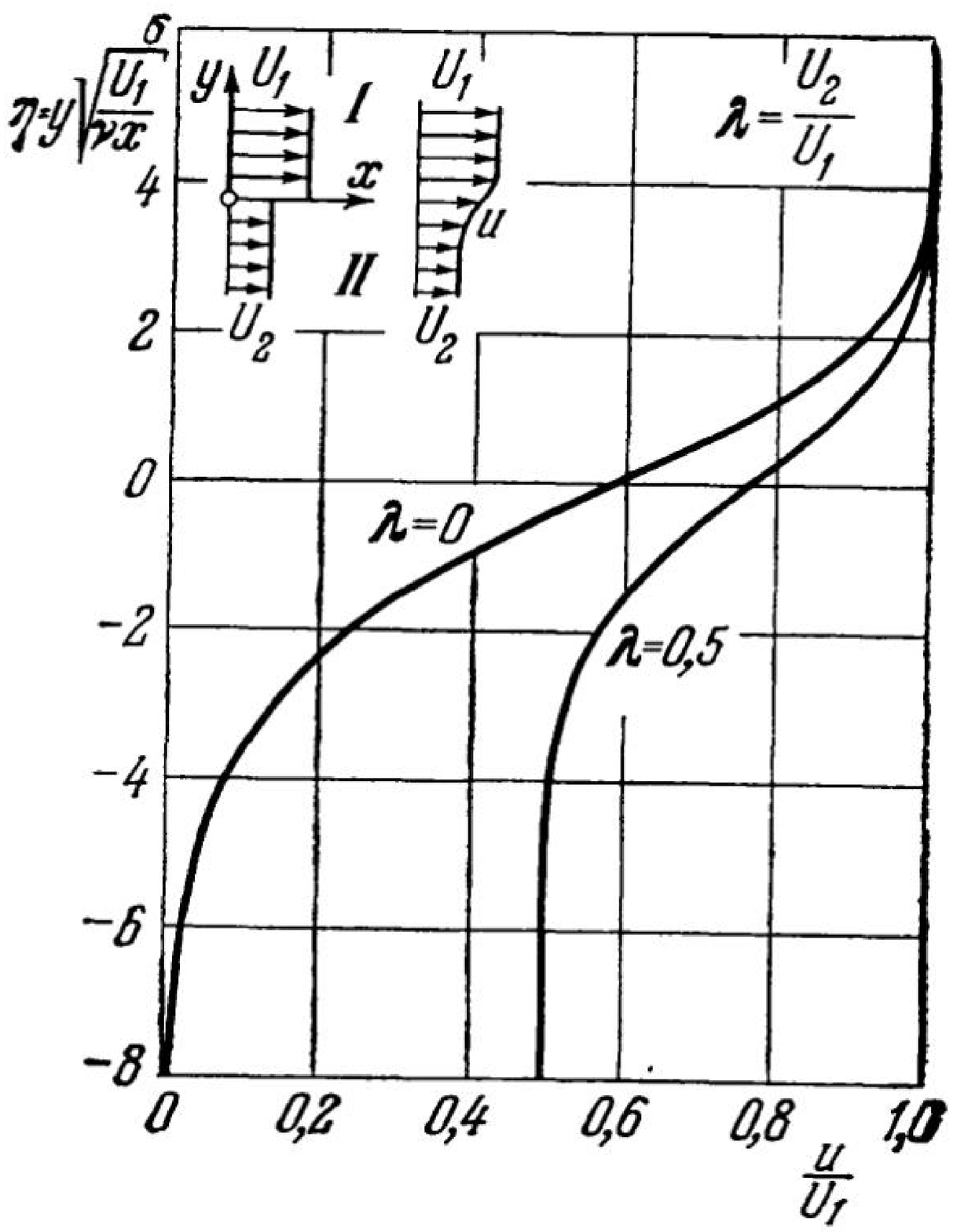}
\end{center}
\begin{center}
Fig.13a (\cite{shlich}, p.181)
\end{center}
\begin{center}
\includegraphics[width = 7cm,height=6cm]{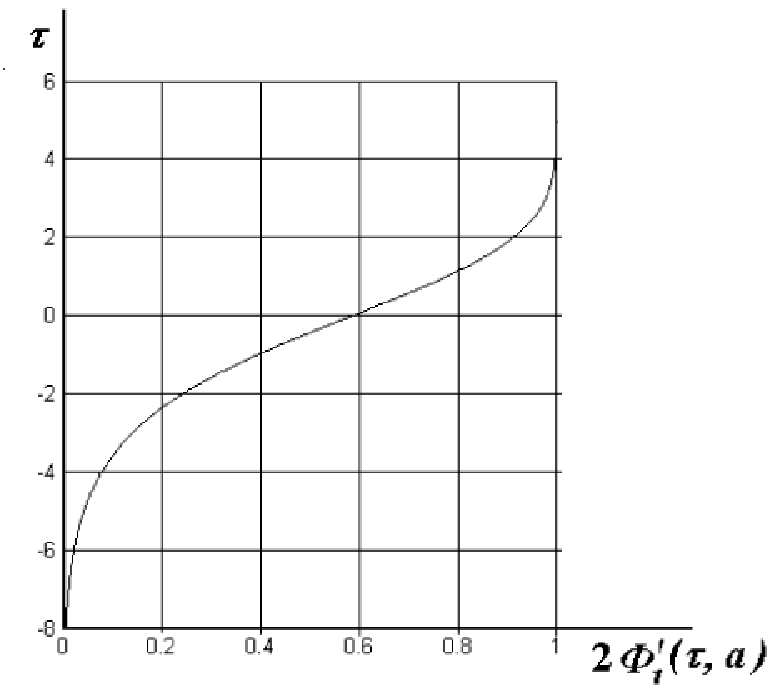}
\end{center}
\begin{center}
Fig.13b ($m=1$, $a=0.619583$)
\end{center}

On Figs.13c,d, the results of our calculations, for \,$m=1$, are presented.
\begin{center}
\includegraphics[width = 10cm,height=6cm]{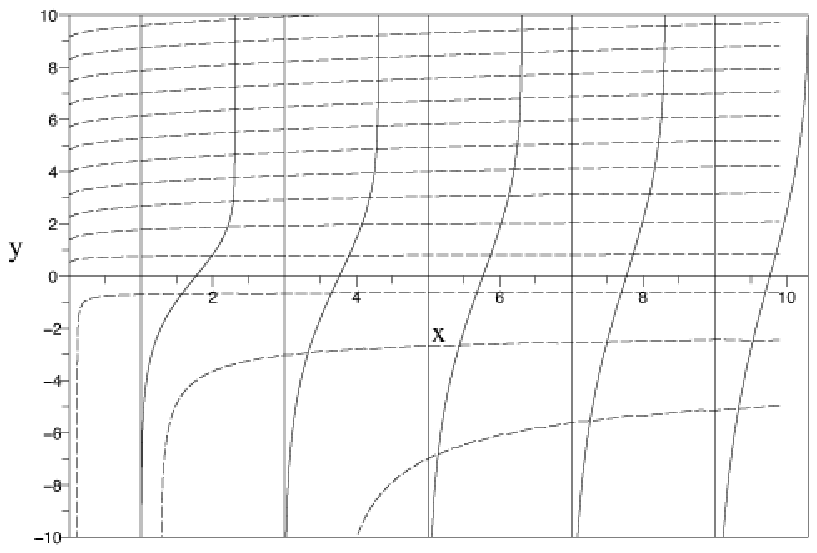}
\end{center}
\begin{center}
Fig.13c. The curves \,$\psi_{1}(x,y)=$const (dotted lines) and the profiles
of the horizontal velocity component  \,$\widetilde u_{1}(y)=u_{1}(x,y)|_{x={\rm const}}$.
\end{center}
\begin{center}
\includegraphics[width = 10cm,height=6cm]{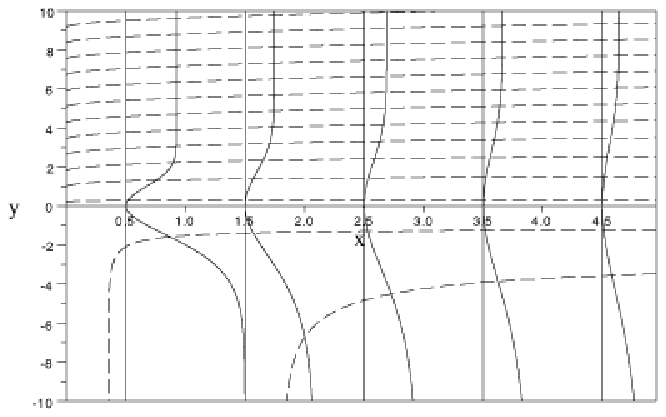}
\end{center}
\begin{center}
Fig.13d. The curves \,$\psi_{1}(x,y)=$const (dotted lines) and the profiles
of the vertical velocity component \,$\widetilde v_{1}(y)=v_{1}(x,y)|_{x={\rm const}}$.
\end{center}
\subsubsection{The case \,$m>1$ \,($m=6/5$)}
\begin{center}
\includegraphics[width = 10cm,height=6cm]{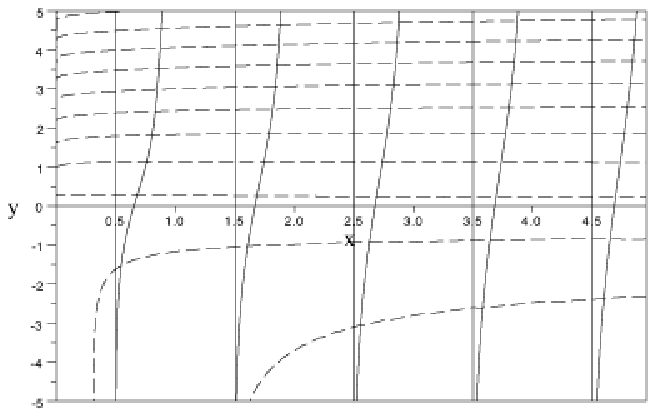}
\end{center}
\begin{center}
Fig.14a. The curves \,$\psi_{6/5}(x,y)=$const (dotted lines) and the profiles
of the horizontal velocity component  with the scale factor:
\,$\widetilde u_{6/5}(y)=0.2\,u_{6/5}(x, y)|_{x={\rm const}}$.
\end{center}
\begin{center}
\includegraphics[width = 10cm,height=6cm]{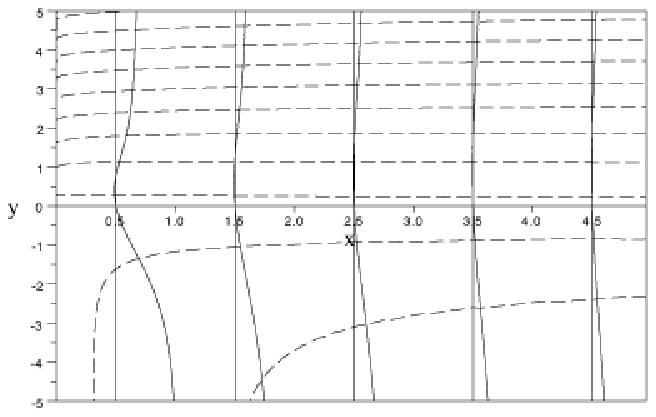}
\end{center}
\begin{center}
Fig.14b. The curves \,$\psi_{6/5}(x,y)=$const and the profiles
of the vertical  velocity component with the scale factor: \,$\widetilde v_{6/5}(y)=0.5\,
v_{6/5}(x, y)|_{x={\rm const}}$.
\end{center}
\section{Appendix A. The Families of Singular (Blow-Up) Solutions to
the Initial Third-Order Nonlinear ODE}
For ODE \eqref{eqf_phi}, let us briefly discuss the behavior of singular
solutions that tend to exact singular solution \eqref{phi1_sing} as
\,$\tau\to\tau _p$. A solution to ODE \eqref{eqf_phi} in the neighborhood of
the pole point \,$\tau=\tau _p$ is sought in the form
\begin{equation}\label{phi_pole}
\Phi (\tau - \tau_p) = \frac{6m}{(m+1)\,(\tau - \tau_p)} \Bigl[1 + Y(\tau -\tau_p)\Bigr],
\quad m>0, \quad \tau_p\in \mathbb{R},
\end{equation}
where
\begin{equation}\label{Y_pole}
\lim_{\tau\to\tau _p} Y(\tau - \tau_p) =0.
\end{equation}
In what follows, it is sufficient without loss of generality to consider the case \,$\tau>\tau_p$.
Then introducing the notation
\begin{equation}\label{x_pole}
x = \tau - \tau_p, \qquad x>0,
\end{equation}
for \,$Y(x)$, we obtain a singular CP for a nonlinear ODE with a
regular singular point at \,$x=0$:
\begin{equation}\label{eqY_pole}
x^3 Y^{\prime \prime \prime} + \frac{3(m-1)}{m+1}\,x^2 Y^{\prime \prime} +
\frac {6(m-1)}{m+1}\,x Y^\prime + 6 Y +
F\left(x, Y, x Y^\prime, x^2 Y^{\prime\prime}\right) = 0, \,\,\, x >0,
 \end{equation}
\begin{equation}\label{Y_con0}
\lim_{x\to +0} Y(x) = \lim_{x \to +0}\Bigl(x\, Y^\prime (x)\Bigr) =
\lim_{x\to +0} \Bigl(x^2 \,Y^{\prime \prime }(x)\Bigr) = 0,
\end{equation}
where the nonlinear part is given by
\begin{equation}\label{F_eqY}
F\left(x,Y,x Y^\prime,x^2 Y^{\prime\prime}\right)=\frac{6}{m+1}\left[
m x^2 Y^{\prime\prime} Y - (m-1)\Bigl(x Y^{\prime}\Bigr)^2 - 2 x Y^{\prime} Y + (m+1)Y^2\right].
\end{equation}
The eigenvalues of the principal linear part of ODE \eqref{eqY_pole} are the
roots of the cubic equation
\begin{equation}\label{cub_eq}
\lambda ^3 - 6\lambda ^2/(m+1)+ \lambda(5m-1)/(m+1) + 6 = 0.
\end{equation}
One root of this equation is \,$\lambda _3 = -1$, and the other two \,$\lambda _{1,2}$ satisfy the equation
\begin{equation}\label{l_1_2_eq}
\lambda ^2 - \lambda (m+7)/(m+1) + 6 = 0,
\end{equation}
whence
\begin{equation}\label{l_1_2_ge}
\lambda _{1,2}(m) = \left(m+7 \pm i\,\sqrt{\varkappa (m)}\right)\Big/\bigl[2(m+1)\bigr],
\qquad m > m_1,
\end{equation}
\begin{equation}\label{l_1_2_le}
\lambda _{1,2}(m) = \left(m+7 \pm \sqrt {-\varkappa (m)}\right)\Big/\bigl[2(m+1)\bigr]>0,
\quad 0 < m \le m_1,
\end{equation}
where
\begin{equation}\label{kappa}
\varkappa (m) = 23m^2 + 34m - 25,
\end{equation}
\begin{equation}\label{m1}
m_1 = \left( -17 + 12\, \sqrt{6}\right)/23 \approx 0.538864213
\end{equation}
is a positive root of the equation \,$\varkappa (m) = 0$ \,($m_1 > 1/2$).

Now, for \,$\lambda_{1,2}$ in \eqref{l_1_2_ge}, denote
\begin{equation}\label{alpha}
\alpha = {\rm Re}\,\lambda_{1,2}(m) = (m+7)/[2(m+1)]>0, \quad m > m_1,
\end{equation}
\begin{equation}\label{beta}
\beta={\rm \Bigl|Im}\,\lambda_{1,2}(m)\Bigr| = \sqrt {\varkappa
(m)}/[2(m+1)]>0, \quad m > m_1.
\end{equation}

A change of variables
\begin{equation}\label{x_Z}
x =\exp{(-t)}, \,\,t >0 \,\,(t= -\ln{x},\,\,0< x < 1), \,\,Z(t)\equiv
Y\left(\exp{(-t)}\right), \,\,t >0,
\end{equation}
leads, for \,$Z(t)$, to a singular CP at infinity for a nonlinear  autonomous ODE:
\begin{equation}\label{eq_Z}
\dddot{Z} + \frac{6}{m+1}\,\ddot{Z} + \frac{5m-1}{m+1}\,\dot{Z} - 6\,Z -
 G\left(Z, \,\dot{Z}, \,\ddot{Z}\right) = 0, \quad t\ge T,
\end{equation}
\begin{equation}\label{Z_con_inf}
\lim_{t\to\infty} Z(t) = \lim_{t\to\infty}\dot{Z}(t) =
\lim_{t\to\infty}\ddot{Z}(t)  = 0,
\end{equation}
where the nonlinear part is given by
\begin{equation}\label{G_Z}
G\left(Z,\dot{Z},\ddot{Z}\right)=\Bigl[6/(m+1)\Bigr]\left[m\,\ddot{Z} Z - (m-1) \dot{Z}^2 + (m+2)
\dot{Z} Z + (m+1) Z^2\right].
\end{equation}

Then using the Lyapunov results for nonlinear autonomous ODEs defined on an infinite interval and
taking into account the substitutions \eqref{x_pole}  and \eqref{x_Z}, we obtain the following
assertion.
\begin{proposition}\label{p_7}
For any fixed \,$m>0$ and \,$\tau _p\in\mathbb{R}$, the singular CP \eqref{eqY_pole}, \eqref{Y_con0}
has a two-parameter family of solutions \,$Y_m \left (\tau - \tau _p, C_1,  C_2\right)$, where
\,$C_1$ and \,$C_2$ are parameters ($C_1, C_2\in\mathbb{R}$). In the neighborhood of the point
\,$\tau=\tau_p$, these solutions can be represented in the principal approximation as:
\begin{equation*}
Y_m\left(\tau-\tau_p,C_1,C_2\right)=|\tau-\tau_p|^{\alpha}\Bigl[C_1 \cos{\bigl(\beta\ln{|\tau-\tau_p|}\bigr)} +
\end{equation*}
\begin{equation}\label{Y_mgem1}
+ C_2 \sin{\bigl(\beta\ln{|\tau-\tau_p|}\bigr)}\Bigr]\Bigl[1+o(1)\Bigr], \quad m>m_1, \quad \tau\to\tau_p,
\end{equation}
where \,$m_1$ \,($m_1>1/2$) is given by \eqref{m1}, \,$\alpha>0$
and \,$\beta>0$ are defined by  \eqref{alpha} and  \eqref{beta} respectively;
\begin{equation*}
Y_{m_1}\left(\tau-\tau_p,C_1,C_2\right)=
\end{equation*}
\begin{equation}\label{Y_meqm1}
 =|\tau-\tau_p|^{\alpha}\Bigl[C_1+C_2\ln{|\tau-\tau_p|}\Bigr]\Bigl[1 + o(1)\Bigr], \quad m=m_1, \quad
\tau\to\tau_p,
\end{equation}
where \,$\alpha=(m_1+7)/[2(m_1+1)]>0$;
\begin{equation*}
Y_m \left(\tau-\tau_p,C_1,C_2\right)=
\end{equation*}
\begin{equation}\label{Y_mlem1}
=\Bigl[C_1|\tau -\tau_p|^{\lambda_1}+C_2|\tau -\tau_p|^{\lambda_2}\Bigr]
\Bigl[1+o(1)\Bigr], \quad 0<m<m_1, \quad \tau\to\tau_p,
\end{equation}
where \,$\lambda_{1,2}=\lambda_{1,2}(m)>0$ are given by \eqref{l_1_2_le}.

More precisely: in the neighborhood of the point \,$\tau = \tau _p$, the two-parameter family
of solutions \,$Y_m \left (\tau - \tau _p, C_1,  C_2\right )$ to the singular CP \eqref{eqY_pole},
\eqref{Y_con0} (where \,$x=|\tau - \tau _p|$) can be represented as the Lyapunov parametric series
in integer powers with principal term given by \eqref{Y_mgem1}, \eqref{Y_meqm1}, or \eqref{Y_mlem1},
respectively (without the \,$o(1)$  term); the  coefficients of this series can be obtained by formal
substitution into ODE \eqref{eqY_pole}.
\end{proposition}

For \,$m = 1/2$, using the exact solutions \eqref{phi2_sing} and taking into
account that \,${\varkappa(1/2)=-9/4<0}$,  \,$\lambda_1(1/2)=3$, \,$\lambda_2(1/2)=2$, and \,$1/2<m_1$,
we also obtain
\begin{corollary}\label{c_8}
For \,$m=1/2$, the singular CP \eqref{eqY_pole}, \eqref{Y_con0} has
the one-parameter set of the exact solutions
\begin{equation*}
Y_{1/2}\bigl(x,a\bigr)=\bigl(ax/2\bigr)\coth{\bigl(ax/2\bigr)}-1, \quad x, a\in\mathbb{R},
\end{equation*}
where \,$a$ is a parameter,  i.e., for any fixed  \,$\tau_p$  ($\tau_p\in\mathbb{R}$), we have
\begin{equation}\label{Y_1/2_taup}
Y_{1/2}\bigl(\tau -\tau_p,a\bigr)=\bigl[a(\tau-\tau_p)/2\bigr]\coth{\Bigl(a(\tau -\tau_p)/2\Bigr)}-1,
\quad \tau,\tau_p,a\in\mathbb{R}.
\end{equation}

This set of the exact solutions belongs to the family \eqref{Y_mlem1} with
\,$C_1=0$,  \,${C_2=C_2(a)=a^2/(3\cdot2^2)}$ and \,$\lambda_2=2$; it has the Lyapunov expansion in the form
\begin{equation*}
Y_{1/2}\bigl(\tau-\tau_p,a\bigr)=C_2(a) (\tau-\tau_p)^2-(1/5) C_2^2(a) (\tau-\tau_p)^4 +
\end{equation*}
\begin{equation}\label{Y_1/2_ls}
 + (2/35) C_2^3(a) (\tau-\tau_p)^6-\ldots, \qquad (\tau-\tau_p)^2<\pi^2.
\end{equation}

In full recording:
\begin{equation}\label{Y_1/2_lsf}
Y_{1/2}\bigl(\tau-\tau_p,a\bigr)=\sum_{n=1}^\infty {D_n\Bigl(C_2(a) (\tau-\tau_p)^2\Bigr)^n},
\qquad (\tau-\tau_p)^2<\pi^2,
\end{equation}
where
\begin{equation}\label{Y_1/2_Dn}
D_n=\bigl[(-1)^{n-1}/(2n!)\bigr] 3^n\cdot 2^{2n} B_n,
\end{equation}
\,$B_n$ are the Bernoulli numbers: \,$B_1=1/6$,  \,$B_2=1/30$, \,$B_3=1/42$, \ldots
(see \rm{\cite{dwight}}).
\end{corollary}
\begin{corollary}\label{c_9}
For any fixed \,$m>0$, nonlinear ODE \eqref{eqf_phi}
has a three-parameter family \,$\Phi_m \left (\tau - \tau_p, C_1, C_2 \right )$ of singular (blow-up) solutions
represented in the form
\begin{equation}\label{Fm_sing1}
\Phi_m\bigl(\tau-\tau_p,C_1,C_2\bigr)=\frac{6m}{(m+1)(\tau-\tau_p)}
\Bigl[1+Y_m\bigl(\tau-\tau_p,C_1,C_2\bigr) \Bigr],
\end{equation}
where \,$\tau _p$,  \,$C_1$,  \,$C_2$  are parameters \,($\tau_p,C_1,C_2\in\mathbb{R}$), and, for fixed
\,$\tau_p\in\mathbb{R}$, the two-parameter set of functions  \,$Y_m\bigl(\tau-\tau_p,C_1,C_2\bigr)$  is
described by Proposition \ref{p_7}.
\end{corollary}
\begin{remark}\label{r_10}
The formulas \eqref{Y_mgem1}--\eqref{Y_mlem1} have been obtained
earlier in {\rm \cite{dies1_86}, \cite{dies2_86}} in an entirely different manner (and, in our
opinion, using an essentially more complex approach). Moreover, the
assertions of Proposition \ref{p_7} and Corollary \ref{c_9} are more complete and
accurate, and the result of Corollary \ref{c_8} is new.
\end{remark}
\section{Appendix B. Some Remarks on the Previous Approach Assuming Complex Analysis of
a Two-Dimensional Dynamical System on "The Poincar\'e Sphere"\, and on Accompanying Singular Problems}
Here we describe very briefly the approach taken in \cite{dies1_86},  \cite{dies2_86}.
We touch in more detail only some auxiliary problems in order to discuss our certain corrections and
remarks and/or to give more complete and exact assertions.
\subsection{Transformation of the Initial Third Order ODE to the First Order
ODE in Nonphysical Variables}
The order of ODE \eqref{eqf_phi} is reduced via treating the desired function \,$\Phi$ as a new
independent variable and introducing a new desired function \,$f(\Phi)$ that is specified along the
trajectory \,$\Phi (\tau)$ of ODE \eqref{eqf_phi} in the form
\begin{equation}\label{f_Phi}
f(\Phi (\tau))=\frac{d\Phi}{d\tau}(\tau).
\end{equation}

For \,$f(\Phi)$ along the trajectory of Eq.\eqref{eqf_phi}, we obtain a second-order ODE (dotted
letters denote derivatives with respect to \,$\Phi$):
\begin{equation}\label{eq_f_Phi}
f\ddot{f} + \dot{f}^2 + \Phi\dot{f} - [(m-1)/m] f = 0.
\end{equation}
In view of the group properties of this ODE, as it is indicated in \cite{dies1_86},  \cite{dies2_86},
new functions \,$F$ and \,$\Psi$ are defined by the formulas
\begin{equation}\label{F_Psi}
\Phi^2F(\Phi)=f(\Phi), \quad \qquad \Psi(\Phi)=\Phi\frac{dF}{d\Phi}(\Phi).
\end{equation}
Differentiating \eqref{F_Psi} with respect to \,$\Phi$ and using \eqref{f_Phi}--\eqref{F_Psi}
yields
\begin{equation}\label{eq1_f}
\frac{df}{d\Phi}=(2F+\Psi)\Phi, \qquad F\frac {d^2f}{d\Phi^2}=-(\Psi+2F)^2-\Psi-\frac {m+1}{m}F;
\end{equation}
\begin{equation}\label{eq2_Phi}
\frac{d^2\Phi}{d\tau^2}=F(2F+\Psi)\Phi^3, \quad \frac{d^3\Phi}{d\tau^3}=
-F\left[\Psi+F(m+1)/m\right]\Phi^4;
\end{equation}
\begin{equation}\label{eq3_F_Psi}
\Psi F\frac{d\Psi}{dF}=-\Psi^2-\Psi-7F\Psi-6F^2-\frac{m+1}{m}F.
\end{equation}
Next, according to \cite{dies1_86},  \cite{dies2_86},  we need to examine the singular points
of nonlinear ODE \eqref{eq3_F_Psi} in the plane of \,$(F, \Psi)$ and to analyze the behavior of the integral
curves in these points in their projections onto the Poincar\'e sphere. To return
to the variables \,$(\tau, \,\Phi, \,\Phi^\prime, \,\Phi^{\prime\prime})$, we have to use the rather complex formulas
\eqref{f_Phi}--\eqref{eq2_Phi}. As a result, the analysis of solutions in terms of the
initial physical variables is rather difficult.

An exception is the special cases where, e.g., Eq.\eqref{eq3_F_Psi} has solutions of the form
\begin{equation}\label{F_Psi_l}
 \Psi=A F+B.
\end{equation}
Here, \,$A$ and \,$B$ are constants which are generally dependent on
\,$m$. Then, for determining all the possible values of \,$A$, \,$B$,
and \,$m$, we obtain the equations
\begin{equation*}
B^2+B=0, \quad A+3AB+7B+(m+1)/m=0, \quad 2A ^2+7A+6=0,
\end{equation*}
which give the solutions
\begin{equation}\label{BA_1_2}
B=B_m=0, \quad  A=A_m=-(m+1)/m, \quad m\in\{1; 2\};
\end{equation}
\begin{equation}\label{BA_1/3_1/2}
B=B_m=-1, \quad A=A_m=(1-6m)/(2m), \quad m\in\{1/3; 1/2\}.
\end{equation}
The solutions in terms of the original variables are then derived using
the formulas \eqref{F_Psi}, \eqref{F_Psi_l} which imply
\begin{equation}\label{f_F_Phi_AB}
F =(C/A)|\Phi|^A-B/A, \qquad f(\Phi)=\Phi^2 F(\Phi),
\end{equation}
where \,$C$ is an arbitrary constant. Then we have finally
\begin{equation}\label{eq_Phi_BA}
\frac{d\Phi_m}{d\tau}=\Phi_m ^2\,\left[(C/A_m)|\Phi_m|^{A_m}-B_m/A_m\right],
\quad m\in\{1/3; 1/2; 1; 2\}.
\end{equation}
As a result, we obtain sets of solutions \eqref{phi_05}, \eqref{phi_1_2}, for \,$m \in \{1/2; \,1; \,2\}$, and a
more complex family of implicit solutions \eqref{phi_1/3}, for \,$m=1/3$, setting in \eqref{eq_Phi_BA}
\,$C = -b^3$ (the example to the last case have been discussed in Subsubsection 3.2.3).

(Note that some typewriting errors were made in \cite{dks_07} in the first formula in
\eqref{f_F_Phi_AB} and hence in \eqref{eq_Phi_BA}. They are corrected here.)
\begin{corollary}\label{c_10}. For ODE \eqref{eq3_F_Psi}, the values of \,$A$, \,$B$, and \,$m$ defined
by \eqref{BA_1_2} and \eqref{BA_1/3_1/2} cover all the possible solutions of form \eqref{F_Psi_l}.
Due to ODEs \eqref{f_F_Phi_AB}, \eqref{eq_Phi_BA}, it implies in turn the formulas \eqref{phi_05}, \eqref{phi_1_2}
and \eqref{phi_1/3} respectively, for the solutions to ODE \eqref{eqf_phi}.
\end{corollary}

The notion of the Poincar\'{e} sphere, as well as the principles of the analysis of second order nonlinear dynamic systems
on the Poincar\'{e} sphere, are given, e.g., in \cite{ALGM}. Fig.B1 from \cite{ALGM}, p. 241, gives an illustration to this approach. Figs.B2--B5 from \cite{dies2_86} give some general illustration for "flows on the Poincar\'{e} sphere"\, to nonlinear ODE \eqref{eq3_F_Psi} (also such figures are presented in \cite{dies1_86}, \cite{dies_84}, \cite{dies_85}).
Here the points \,$A$, \,$B$, \,$C$ correspond to the singular points of ODE \eqref{eq3_F_Psi} in $\{F,\Psi\}$ plane
with coordinates \,$(0,0)$, \,$(0,-1)$ and \,$(-(m+1)/(6m), 0)$ respectively; points like \,$Q$, \,$E$, \,$G$ on the circle correspond to the infinite singular points of ODE \eqref{eq3_F_Psi}, whose identification and analysis seem to be highly complicated. Note that, since no problem is specified now for ODE \eqref{eq3_F_Psi} in the nonphysical variables,
one have to examine the behavior of all the trajectories on the Poincar\'{e} sphere and to select those corresponding to sought-for solutions of ODE \eqref{eqf_phi} in the original variables (as noted before, the transition to these variables is nontrivial), i.e., to the solutions of problem \eqref{eqf_phi}--\eqref{con3f_phi} when they exist.
Consequently, illustrations for auxiliary nonphysical entities \,$\{F,\Psi\}$ like those on Figs.B2--B5 seem to be rather difficult to interpret, especially when one passes from them to the description of the flow in the original variables
(flow pictures were not presented in the cited papers). At the same time, the pictures are difficult to distinguish in essence for different \,$m$'s. Probably to us, not being experts in this kind of analysis, the treatment of the problems in the papers cited has its own mathematical interest. We would note once more that regardless of the difficulties we point out here,
many of the facts and formulas in  \cite{dies1_86}, \cite{dies2_86} seem to be new.

Let us remark here: for enumeration of figures in this Appendix B, we use an additional symbol, namely a letter B.
\begin{center}
\includegraphics[width = 10cm,height=6cm]{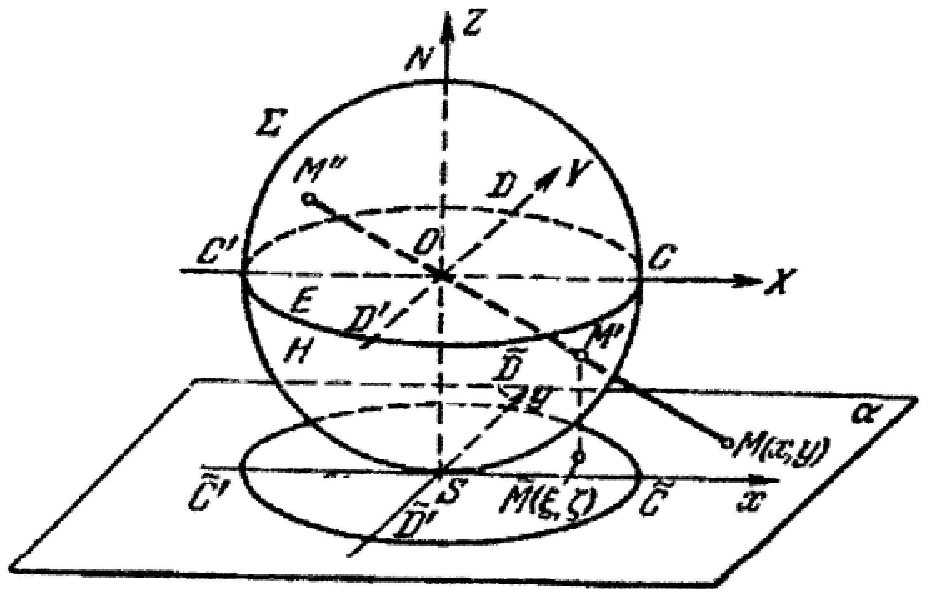}
\end{center}
\begin{center}
Fig.B1 (\cite{ALGM}, p.241)
\end{center}
\begin{center}
\includegraphics[width = 6cm,height=6cm]{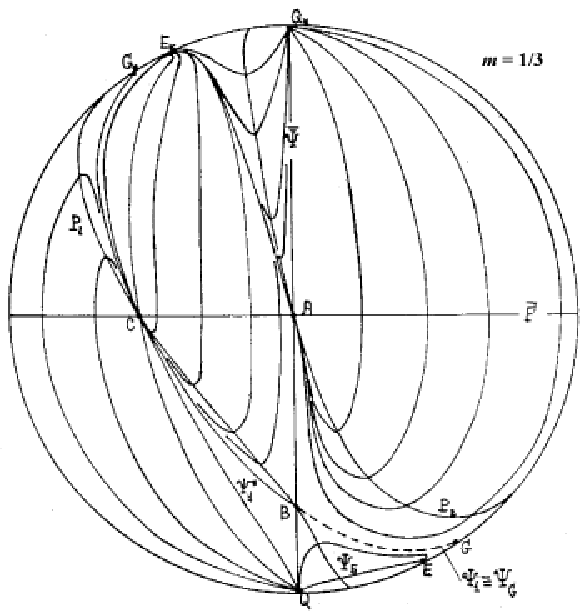}
\end{center}
\begin{center}
Fig.B2 (\cite{dies2_86}; \,$m=1/3$)
\end{center}
\begin{center}
\includegraphics[width = 6cm,height=6cm]{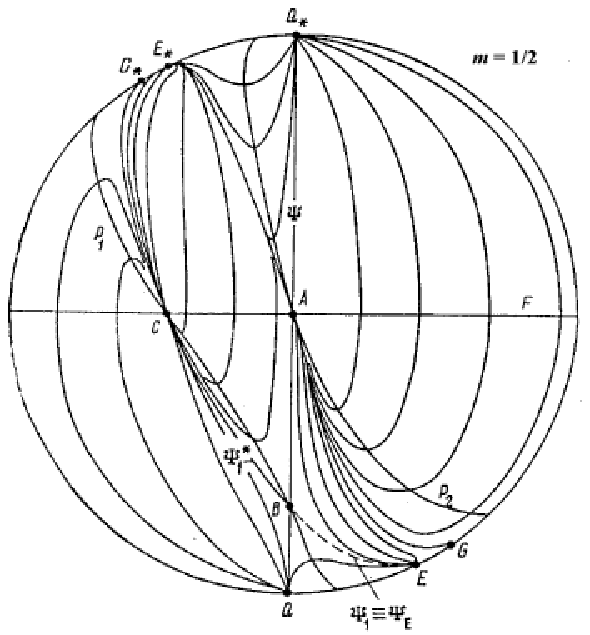}
\end{center}
\begin{center}
Fig.B3 (\cite{dies2_86}; \,$m=1/2$)
\end{center}
\begin{center}
\includegraphics[width = 6cm,height=6cm]{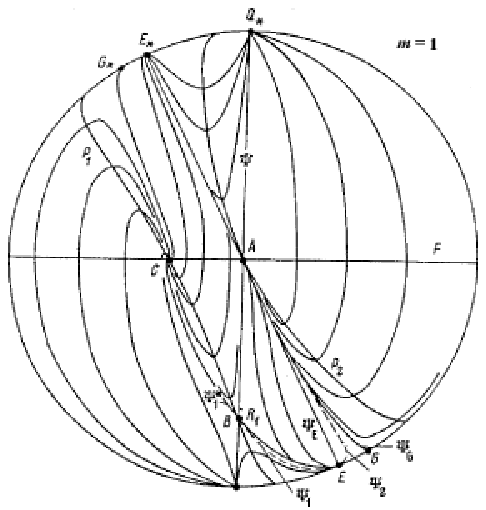}
\end{center}
\begin{center}
Fig.B4 (\cite{dies2_86}; \,$m=1$)
\end{center}
\begin{center}
\includegraphics[width = 6cm,height=6cm]{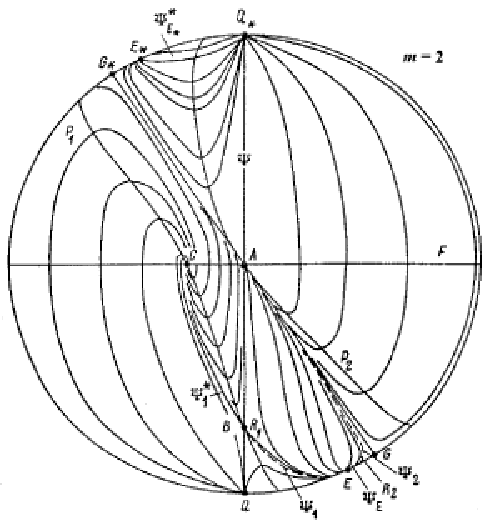}
\end{center}
\begin{center}
Fig.B5 (\cite{dies2_86}; \,$m=2$)
\end{center}
\subsection{An Auxiliary Singular Nonlinear Problem for the Phase Trajectories}
In \cite{dies2_86}, the ODE problem  \eqref{eq_f_Phi} that corresponds to the original problem \eqref{eqf_phi}--\eqref{con3f_phi} is formulated as (in our notation):
\begin{equation}\label{eqf_Phi}
f\ddot{f}+\dot{f}^2+\Phi\dot{f}-[(m-1)/m]f=0, \quad
-a<\Phi<+\infty \quad (a>0, m>0),
\end{equation}
\begin{equation}\label{f_Phia}
f=a(\Phi+a)-\frac{1}{4m} (\Phi+a)^2 + O((\Phi+a)^3),\qquad \Phi\to -a+0,
\end{equation}
\begin{equation}\label{f_Phi_inf}
f=m b^{1/m}\Phi^{(m-1)/m} + \frac{m(m-1)(m-2)}{m+1} b^{2/m} \Phi^{-2/m} + O(\Phi^{-1-3/m}),\,\, \Phi\to +\infty.
\end{equation}

Conditions \eqref{f_Phia}, \eqref{f_Phi_inf} are treated there as boundary conditions, i.e., problem \eqref{eqf_Phi}--\eqref{f_Phi_inf} assumed to be a singular two point BVP. In \cite{dies1_86}, there is one more term
in expansion  \eqref{f_Phi_inf} that is an exponentially decaying function with an arbitrary constant factor
(as we will see later, this term is correct). But in \cite{dies2_86} the term is dropped (possibly because the presence of arbitrary constant does not meet with the treating a problem of \eqref{eqf_Phi}--\eqref{f_Phi_inf} type as a two point BVP).

In effect, as we will show later on, problem \eqref{eqf_Phi}, \eqref{f_Phia} is a singular initial one.
Its setting will be specified later on. The problem has a unique solution for each \,$m\not=0$, and expression
\eqref{f_Phi_inf} gives the main term of asymptotic representation for the solution for \,$m>1/2$ and big \,$\Phi>0$.
\subsubsection{Statement of initial problem for nonlinear ODE degenerating with respect to the phase variable}
Consider ODE \eqref{eq_f_Phi} from the viewpoint of the paper given. Let us remark that \eqref{f_Phi}
and the formula
\begin{equation}\label{df_Phi}
\dot{f}(\Phi(\tau)) = \Phi ^{\prime\prime}(\tau)/\Phi^\prime(\tau)
\end{equation}
hold along the trajectory of ODE \eqref{eqf_phi}.

First, let \,$\Phi(\tau)$ be a solution of singular nonlinear CP \eqref{eqf_phi}, \eqref{ccon1f_phi}. Then taking into
account formulas \eqref{f_Phi}, \eqref{df_Phi}, conditions \eqref{ccon1f_phi} at \,$\tau\to-\infty$, Proposition \ref{p_1} and expansion \eqref{l_ser}, we obtain the limit conditions for solutions to ODE \eqref{eq_f_Phi} as \,${\Phi\to-a+0}$: \,\,$\lim_{\Phi\to-a+0}{f(\Phi)}=0$, \,\,$\lim_{\Phi\to-a+0}{\dot{f}(\Phi)}=a$.

As a result, we obtain a singular nonlinear CP
\begin{equation}\label{CP_eq_f}
(f\dot{f}+\Phi f\dot{)}=[(2m-1)/m] f, \quad \Phi>-a \quad (a>0, \,\,m\not=0),
\end{equation}
\begin{equation}\label{CP_cond_f}
\lim_{\Phi\to a+0}{f(\Phi)}=0, \qquad  \lim_{\Phi\to-a+0}{\dot{f}(\Phi)}=a,
\end{equation}
where ODE \eqref{CP_eq_f} is the same as \eqref{eq_f_Phi} but represented in a more convenient
form.

Now let \,$\Phi(\tau)$ be a solution of singular nonlinear IBVP \eqref{eq2_phi}--\eqref{eq2_0}, \eqref{con3f_phi}, for some fixed \,$b>0$, \,$a=a(b)>0$ and \,$m:\,1/2<m<\infty$ (according to Theorem \ref{t_3}). Then the solution  \,$f(\Phi)$ of singular nonlinear CP \eqref{CP_eq_f}, \eqref{CP_cond_f} must satisfy limit condition
\begin{equation}\label{f_inf1}
\lim_{\Phi\to\infty}{[f(\Phi)/\Phi^{(m-1)/m}]}=m b^{1/m} \qquad (b>0, \,\, 1/2<m<\infty).
\end{equation}

Above all we consider singular nonlinear CP \eqref{CP_eq_f}, \eqref{CP_cond_f}.
For its solutions in particular  the following relation is valid:
\begin{equation}\label{rel_f_Phi}
f(\Phi) \dot{f}(\Phi)+\Phi f(\Phi)=[(2m-1)/m] \int_{-a}^\Phi {f(s)\,ds}.
\end{equation}
\begin{remark}\label{r_11}
For \,$m=\infty$ and \,$m=1/2$, there exist the exact solutions \,$f_m(\Phi,a)$ to singular CP
\eqref{CP_eq_f},  \eqref{CP_cond_f}:
\begin{equation}\label{ff_inf}
f_\infty(\Phi,a)=a\,(\Phi+a),
\end{equation}
\begin{equation}\label{ff_1/2}
f_{1/2}(\Phi,a)=(a^2 - \Phi^2)/2.
\end{equation}

These solutions don't satisfy condition \eqref{f_inf1}. For the initial singular CP \eqref{eqf_phi}, \eqref{ccon1f_phi}, function \eqref{ff_inf} corresponds to the exact solution \,${\Phi_\infty(\tau-\tau_s,a)}$ defined by formula
\eqref{phi_inf}; for \,$a^2>\Phi^2$, function \eqref{ff_1/2} corresponds to the exact solution \,$\Phi_{1/2}(\tau-\tau_s,a)$ defined by \eqref{phi_05}, and, for \,$a^2<\Phi^2$, it corresponds to the exact solution
\,$\Phi^{(2)}_{{\rm sing},1/2}(\tau-\tau_p,a)$ defined by \eqref{phi2_sing}.
\end{remark}
\begin{remark}\label{r_12}
ODE \eqref{CP_eq_f} is invariant with respect to the following change of variables:
\begin{equation*}
f_{\rm new}=f/a^2, \qquad \Phi_{\rm new}=\Phi/a.
\end{equation*}
Then, for problem \eqref{CP_eq_f},  \eqref{CP_cond_f}, it is enough to put \,$a=1$ because
the following relation  holds:
\begin{equation}\label{f_Phi_a}
f(\Phi,a) = a^2f(\Phi/a,1).
\end{equation}

So, while remaining the \,$a$ in the formulas,  we put \,$a=1$  for the numerical examples.
\end{remark}
\subsubsection{Auxiliary singular CP for nonlinear ODE with a regular singular point at
zero and solvability of the initial degenerating problem}
For the singular nonlinear problem \eqref{CP_eq_f},  \eqref{CP_cond_f}, it follows that
ODE \eqref{CP_eq_f} degenerates with respect to the phase variable as
\,$\Phi\to -a+0$. To study this degeneration, we set
\begin{equation}\label{deg_f_a}
f(\Phi) = A(\Phi + a)^\alpha [1 + o(1)],
\quad \dot{f}(\Phi) = \alpha A(\Phi +a )^{\alpha - 1}[1+o(1)],
\end{equation}
\begin{equation}\label{deg1_f_a}
\ddot{f}=\alpha (\alpha - 1) A (\Phi +a)^{\alpha - 2}[1+o(1)],
\qquad \Phi \to -a+0,
\end{equation}
where \,$\alpha>0$.

Substituting \eqref{deg_f_a},  \eqref{deg1_f_a} into \eqref{CP_eq_f} yields
(in principal as \,$\Phi\to -a$):
\begin{equation*}
\alpha (\alpha-1)A^2(\Phi+a)^{2\alpha -2} + \alpha^2 A^2
(\Phi+a)^{2\alpha -2} + \alpha A (\Phi+a)^\alpha  -
\end{equation*}
\begin{equation*}
-a \alpha A (\Phi+a)^{\alpha -1} - [(m-1)/m] A(\Phi+a)^\alpha +
\ldots = 0.
\end{equation*}

To satisfy this equation as \,$\Phi \to -a$, we need to take away the terms with the least powers.
Then we have the relations
\begin{equation*}
2\alpha-2=\alpha-1, \qquad A^2+A(-a)=0,
\end{equation*}
so that
\begin{equation}\label{alpha A}
\alpha=1, \qquad A=a.
\end{equation}

Then, taking into account \eqref{f_Phi_a}--\eqref{alpha A}, we introduce the changes of
variables
\begin{equation}\label{t_chi}
t = \Phi/a + 1 \quad (t>0), \quad f(\Phi) = \widetilde f(t)= a^2\,t\,[1 + \chi (t)],
\end{equation}
and obtain a singular CP for the ODE with respect to the new desired function \,$\chi (t)$,
using \eqref{CP_eq_f}, \eqref{CP_cond_f}, and \eqref{t_chi}; this problem has a regular singular point at
\,$t=0$ (the classification of singular points for systems of linear and nonlinear ODEs can be found, e.g.,
in \cite{kl}, \cite{was}):
\begin{equation}\label{eq_chi}
t^2\ddot{\chi} + 3t\dot{\chi} + \chi =
G(t,\chi,t\dot{\chi}) + \eta(t,m), \qquad t>0,
\end{equation}
\begin{equation}\label{cond_chi}
\lim_{t\to +0}\chi(t)=0, \qquad \lim_{t\to
+0}\left[t\dot{\chi}(t)\right]=0,
\end{equation}
where
\begin{equation}\label{G_chi}
G(t,\chi,t\dot{\chi})=t\dot{\chi}\left[(1 + \chi)^{-1} - 1\right]
-t\dot{\chi}(1 + \chi)^{-1}\left[t + t\dot{\chi}\right], \quad
\end{equation}
\begin{equation}\label{eta_t}
\eta (t,m)= - t/m, \quad m>0,
\end{equation}
so that \,$G(t, \chi, t\dot{\chi})$ is a holomorphic function at the point \,$ (t, \chi, t\dot{\chi}) = (0, 0,0)$,
and
\begin{equation*}
G(t,0,0)\equiv 0, \quad \frac{\partial G}{\partial \chi} (t,0,0)
\equiv 0, \quad \frac{\partial G}{\partial {(t \dot {\chi})}} (0,0,0)=0,
\quad  \eta (0,m) = 0.
\end{equation*}

To find the eigenvalues of the linearized problem, we derive \,${\lambda (\lambda -1) + 3\lambda + 1=0}$, hence
\,$\lambda_{1,2}=-1$. Then the results concerning singular CPs for nonlinear ODEs
(in particular see Theorem 5 in \cite{k_83} and references therein) imply the following assertion.
\begin{proposition}\label{p_8}. For any fixed \,$m\not= 0$, singular
nonlinear CP \eqref{eq_chi}, \eqref{cond_chi} has a unique solution \,$\chi_m(t)$; it is
a holomorphic function at \,$t=0$,
\begin{equation}\label{chi_0_ser}
\chi_m(t)=\sum_{k=1}^{\infty}\chi_{m,k} \,t^k, \quad
|t|\le t_0, \quad t_0>0,
\end{equation}
where the substitution of \eqref{chi_0_ser}  into \eqref{eq_chi} yields:
\begin{equation*}
\chi_{m,1}= -1/(4m), \quad \chi_{m,k}=(k+1)^{-2} \Big\{-\chi_{m,k-1}\, [1+m(k-1)]/m-
\end{equation*}
\begin{equation}\label{chi_rec}
-\sum_{l=1}^{k-1}{[l(k+3)+1]\,\chi_{m,l}\,\chi_{m,k-l}}\Big\},
\quad k=2,3,\ldots
\end{equation}
For \,$m=\infty$ and \,$m=1/2$, there are the exact solutions:
\begin{equation}\label{chi_inf_1/2}
\chi_\infty(t)\equiv 0,\qquad \chi_{1/2}(t)= - t/2.
\end{equation}
\end{proposition}

In particular for \,$m>0$, from Proposition \ref{p_8} we have the relations
\begin{equation*}
\chi_m(0)=0, \qquad \chi^\prime_m(0)=-1/(4m)<0, \qquad
\chi^{\prime\prime}_m(0)=(2m-1)/(36 m^2),
\end{equation*}
wherefrom
\begin{equation*}
\chi^{\prime\prime}_m(0)>0,  \quad m>1/2; \quad \chi^{\prime\prime}_{1/2}(0)=0;\quad
\chi^{\prime\prime}_m(0)<0, \quad  m<1/2.
\end{equation*}
Moreover, for any fixed \,$m\not=0$, as it is easy to check, the solutions of nonlinear ODE
\eqref{eq_chi} have no singularities of a pole-type at finite points \,$t>0$.
\begin{corollary}\label{c_11} For any fixed \,$m>1/2$, the solution \,$\chi_m(t)$ of the
singular nonlinear CP \eqref{eq_chi}, \eqref{cond_chi} exists globally on
\,$\mathbb{R}_+$ and satisfies inequality \,$\chi_m(t)>-1$ \,$\forall t\in \mathbb{R}_+$.
\end{corollary}

As a corollary of Proposition \ref{p_8}, relation \eqref{rel_f_Phi} and the fact that, for any finite
\,$\Phi > -a$, the solutions of nonlinear ODE \eqref{CP_eq_f} have no singularities of a pole-type,
we obtain finally a following assertion.
\begin{theorem}\label{t_4} For any fixed \,$a>0$ and \,$m\not=0$, singular
nonlinear CP \eqref{CP_eq_f}, \eqref{CP_cond_f} has a unique solution \,$f_m(\Phi, a)$; it
is a holomorphic function at \,$\Phi=-a$:
\begin{equation}\label{fm_Phi_a}
f_m(\Phi,a)=a(\Phi+a) \left[1+\sum_{k=1}^\infty
\chi_{m,k}\,(\Phi+a)^k/a^k \right],
\end{equation}
where  \,$\chi_{m,k}$ $(k \ge 1)$ are defined by recurrence relations  \eqref{chi_rec}.

For \,$m=\infty$ and \,$m=1/2$, there are the exact solutions:
\begin{equation}\label{finf_Phi_a}
f_\infty (\Phi,a) = a\,(\Phi +a),
\end{equation}
\begin{equation}\label{f1/2_Phi_a}
f_{1/2}(\Phi,a)=(a^2-\Phi^2)/2=a (\Phi +a)[1-(\Phi + a)/(2a)].
\end{equation}

Moreover, for \,$m: 1/2< m < \infty $, the solution \,$f_m(\Phi, a)$ of singular nonlinear CP \eqref{CP_eq_f},
\eqref{CP_cond_f} exists globally and is a positive function on the interval \,$(-a, \infty)$.
\end{theorem}

From \eqref{fm_Phi_a}, \eqref{chi_rec}, in particular we have as \,$\Phi +a \to 0$:
\begin{equation}\label{fm_Phi_par}
f_m(\Phi,a)=a(\Phi+a)\left [1-\frac {\Phi+a}{4am}+
\frac{(2m -1)(\Phi+a)^2}{72a^2m^2} \right]+O((\Phi+a)^4).
\end{equation}

In \cite{dies1_86}, \cite{dies2_86} two terms of the expansion in a formula of form \eqref{fm_Phi_par} were
given (see here \eqref{f_Phia}), but the character of the representation was not discussed. An
exact assertion follows from the above argument, and the terms of converging series \eqref{fm_Phi_a} can be
found via the formal substitution of this series into \eqref{CP_eq_f}.

We also would note that singular nonlinear CP \eqref{eq_chi}, \eqref{cond_chi} was briefly described
in our paper \cite{dks_07} with indication of two terms of expansion for its solution \,$\chi_m(t)$
(see Remark 5 therein). Unfortunately, a mistake was made in the sign in formula \eqref{G_chi} for
\,$G(t,\chi,t\dot{\chi})$ that resulted in a mistake in the second (resp. third) term in the expansion
for \,$\chi_m(t)$ (resp. \,$f_m(\Phi,a)$); the mistakes are corrected here.
\subsubsection{Asymptotic behavior at infinity of the solutions to the
initial singular problem for different values of \,$m>0$}
Global behavior of the solution of singular nonlinear CP \eqref{CP_eq_f}, \eqref{CP_cond_f} for different
values of a self-similarity parameter \,$m>0$ gives a rather difficult problem.

To clarify the principle behavior of the solutions of ODE \eqref{CP_eq_f} for big \,$\Phi$,
assume
\begin{equation}\label{f_Phi_bet}
f(\Phi)=B\,\Phi^\beta\,[1+o(1)], \qquad \Phi\to \infty.
\end{equation}
Substituting \eqref{f_Phi_bet} into \eqref{CP_eq_f} and keeping the main terms give
\begin{equation}\label{f_Phi_bet1}
B^2 \beta(\beta-1) \Phi^{2\beta-2} + B^2\beta^2\Phi^{2\beta-2} + B \beta\Phi^\beta-[(m-1)/m] B \Phi^\beta +\ldots=0.
\end{equation}
Now, after removing in \eqref{f_Phi_bet1} the most growing terms with respect to \,$\Phi$ big, we get two opportunities:

Case I. The biggest power in \eqref{f_Phi_bet} is \,$\beta$. Then
\begin{equation}\label{bet_1}
\beta=(m-1)/m, \qquad m>0,
\end{equation}
and the parameter \,$B$ (\,$B\not=0$) is free; indeed, in this case, for \,$m>0$, we obtain
\begin{equation*}
2\beta-2=-2/m <(m-1)/m=\beta.
\end{equation*}

Case II. The two biggest powers are equal: $2\beta-2=\beta$ in \eqref{f_Phi_bet}; then
\begin{equation}\label{B_2}
\beta=2, \qquad B=-(m+1)/(6m), \qquad m\not=0.
\end{equation}

We will give a more accurate proof of the existence of the solution families
with a change of the dependent variable in \eqref{CP_eq_f}:
\begin{equation}\label{f_Phi_z}
f(\Phi)=B\,\Phi^\beta\,[1+Z(\Phi)].
\end{equation}
Here for the new function \,$Z(\Phi)$ the conditions
\,$\lim_{\Phi\to\infty}{Z(\Phi)}=\lim_{\Phi\to\infty}{\dot{Z}(\Phi)}=\lim_{\Phi\to\infty}{\ddot{Z}(\Phi)}=0$
must hold.

Now we have a singular CP for \,$Z(\Phi)$ at infinity:
\begin{equation*}
\ddot{Z}+[4\beta/\Phi+\Phi^{1-\beta}/B]\,\dot{Z}
+[(2\beta^2-\beta)/\Phi^2]\,Z +(\Phi^{1-\beta}/B)\, [(1+Z)^{-1}-1]\,\dot{Z}+
\end{equation*}
\begin{equation}\label{Z_CP}
+\dot{Z}^2\,(1+Z)^{-1}+(2\beta^2-\beta)/\Phi^2 + [\beta-(m-1)/m]/(B \Phi^\beta)=0, \quad \Phi\gg 1,
\end{equation}
\begin{equation}\label{Z_cond}
\lim_{\Phi\to\infty}{Z(\Phi)}=\lim_{\Phi\to\infty}{\dot{Z}(\Phi)}=0.
\end{equation}
\begin{remark}\label{r_13}.
For case II, singular nonlinear CP \eqref{Z_CP}, \eqref{Z_cond} at infinity always has the trivial solution  \,${Z\equiv 0}$,
and the same is true for case I with \,$m=1$, when $\beta=0$, and with \,$m=2$, when \,$\beta=1/2$ ($B\not=0$ being arbitrary). Then expression \eqref{f_Phi_z} gives the exact solutions \,$f_{\rm II,m}$ and \,$f_{\rm I,m}$ to ODE \eqref{CP_eq_f} that are not the solutions of singular CP \eqref{CP_eq_f}, \eqref{CP_cond_f}:
\begin{equation}\label{f_beta_2}
f_{\rm II,m}= -\Phi^2 (m+1)/(6m), \qquad m\not=0,
\end{equation}
\begin{equation}\label{f_m_1_2}
f_{\rm I,1}\equiv B, \qquad f_{\rm I,2}=B \sqrt{\Phi}.
\end{equation}
For the initial ODE \eqref{eqf_phi}, function \eqref{f_beta_2} corresponds to the exact blow-up solution
\,$\Phi_{\rm sing, m}^{(1)}(\tau-\tau_p,a)$, defined by formula \eqref{phi1_sing}, and functions \eqref{f_m_1_2}
lead to the exact solutions \,${\Phi_1(\tau-\tau_s,a)}$ and \,${\Phi_2(\tau-\tau_s,a)}$, defined by formulas \eqref{phi_1_2}.
\end{remark}

\textbf{Case I.} For \,$m>0$ and \,$\beta=(m-1)/m$ we have
\begin{equation}\label{caseI}
1-\beta=1/m,  \qquad 2\beta^2-\beta=(m-1)(m-2)/m^2.
\end{equation}

Now we make the change of the independent variable in \eqref{Z_CP}, \eqref{Z_cond}:
\begin{equation}\label{Z_Phi}
x=\Phi^{(m+1)/m}, \qquad \Phi=x^{m/(m+1)}.
\end{equation}

Then, denoting \,$\widetilde{Z}(x)=Z(\Phi(x))$ and taking into account equalities \eqref{caseI}, from
\eqref{Z_CP}, \eqref{Z_cond} we get a singular nonlinear CP for \,$\widetilde{Z}(x)$ at infinity depending on
\,$B\not=0$ and  \,$m>0$ as the parameters (prime stands for differentiating by \,$x$):
\begin{equation*}
\widetilde{Z}^{\prime\prime}+\widetilde{Z}^{\prime}\,\left [\frac{m}{B(m+1)}+\frac{4m-3}{x(m+1)}\right]+\widetilde{Z} \,\frac{(m-1)(m-2)}{(m+1)^2\,x^2}+
\end{equation*}
\begin{equation}\label{Z_CP_x}
+\frac{m}{B(m+1)} \widetilde{Z}^{\prime}\, [(1+\widetilde{Z})^{-1}-1]+\widetilde{Z}^{\prime 2}\,(1+\widetilde{Z})^{-1} + \frac{(m-1)(m-2)}{(m+1)^2\,x^2}=0, \quad x\gg 1,
\end{equation}
\begin{equation}\label{Z_cond_x}
\lim_{x\to\infty}{\widetilde{Z}(x)}=\lim_{x\to\infty}{\widetilde{Z}^\prime(x)}=0.
\end{equation}

According to \cite{was}, nonlinear ODE \eqref{Z_CP_x} has an irregular singular point of a rank \,$1$
as \,$x\to\infty$, and the following assertion is valid.
\begin{proposition}\label{p_9}
For any fixed \,$B\ne 0$ and \,$m>0$, singular nonlinear CP \eqref{Z_CP_x}, \eqref{Z_cond_x} has
a particular solution \,$\theta (x)=\theta (x,m,B)$ that can be represented as a formal series
\begin{equation}\label{teta}
\theta (x,m,B)=\sum_{k=1}^{\infty} \theta_k/x^k,
\end{equation}
where coefficients in the expansion are found by formal substitution of series \eqref{teta}
into ODE \eqref{Z_CP_x}:
\begin{equation}\label{teta_1}
\theta_1= B(m-1)(m-2)/[m(m+1)],
\end{equation}
\begin{equation}\label{teta_2}
\theta_2=B^2(m-1)(m-2)(11-m^2-8m)/[2m^2(m+1)^2],
\end{equation}
\begin{equation*}
\theta_{k-1}=\frac{B(m+1)}{m(k-1)}\Bigl\{\theta_{k-2}\left [(k-1)(k-2)-\frac {4m-3}{m+1}(k-2) +
\frac{(m-1)(m-2)}{(m+1)^2}\right ]+
\end{equation*}
\begin{equation}\label{teta_k-1}
+ F_{k-1}(\theta_{1},\ldots,\theta_{k-2},m, B)\Bigr\}, \qquad k=3,4, \ldots,
\end{equation}
where \,$F_{k-1}$ are derived from \,$\theta_{1},\ldots,\theta_{k-2}$
by addition and multiplication operations only; moreover, there exists a true solution to ODE
\eqref{Z_CP_x} that has series \eqref{teta}--\eqref{teta_k-1} as its asymptotic expansion for large \,$x$.
\end{proposition}

After substituting the difference \,$w(x)=\widetilde{Z}(x)-\theta(x)$ into \eqref{Z_CP_x}, \eqref{Z_cond_x} we get the singular nonlinear CP at infinity for \,$w(x)$:
\begin{equation*}
w^{\prime\prime}+w^{\prime}\,\left[\frac{m}{B(m+1)}+\frac{4m-3}{(m+1)x}\right]+w\frac{(m-1)(m-2)}{(m+1)^2 x^2}+
\end{equation*}
\begin{equation*}
+\frac{m}{B(m+1)}\Bigl\{(w^\prime+\theta^\prime(t))\left[(1+w+\theta(t))^{-1}-1\right]-\theta^\prime(t)\left
[(1+\theta(t))^{-1}-1\right]\Bigr\} +
\end{equation*}
\begin{equation}\label{w_eq_x}
+(w^\prime+\theta^\prime(t))^2(1+w+\theta(t))^{-1}-\theta^\prime(t)^2(1+\theta (t))^{-1}=0,
\end{equation}
\begin{equation}\label{w_x_cond}
\lim_{x\to\infty}{w(x)}=\lim_{x\to\infty}{w^\prime(x)}=0.
\end{equation}
Keeping only main terms in ODE \eqref{w_eq_x} for \,$w(x)$ with coefficients tending to zero
not faster than \,$1/x$ as \,$x\to\infty$, we get a linear ODE:
\begin{equation}\label{tild_w_eq_x}
\widetilde{w}^{\prime\prime}+\widetilde{w}^{\prime}
\,\left [\frac{m}{B(m+1)}+\frac{3m^2+4m-5}{x(m+1)^2}\right]=0.
\end{equation}
For \,$B>0$, we get the representation \,${\widetilde{w}(x)=P\, x^{-a_2}\,\exp{(-a_1 x)}}$
for the one-parameter family of solutions to ODE \eqref{tild_w_eq_x} tending to zero
as \,$x\to\infty$, where \,$P$ is an arbitrary constant,
\begin{equation}\label{a_1_2}
a_1= m/[B(m+1)]>0, \qquad a_2=-(2m^2+4m-4)/[m(m+1)].
\end{equation}
\begin{proposition}\label{p_10}
For any fixed \,$m>0$ and \,$B>0$, singular nonlinear CP \eqref{w_eq_x}, \eqref{w_x_cond} with data at infinity has
a one-parameter set of solutions \,$w(x)=w(x,B,m,P)$; these solutions are represented
by the exponential parametric Lyapunov series
\begin{equation}\label{w_C1}
w\left(x,B,m,P\right)=P\,x^{-a_2}\,\exp{(-a_1 x)}+\sum_{k=2}^\infty{C_k(x) P^k\, x^{-k a_2}\,\exp{(-k a_1 x)}},
\end{equation}
where \,$a_1$ and \,$a_2$ are defined by \eqref{a_1_2}, \,$P$ is a parameter and the functions
\,$C_k(x)$ ($k=2,3,\ldots$) are of at most of power growth for large \,$x>0$.
\end{proposition}

Summarizing and taking into account the changes of variables, we get the following.
\begin{proposition}\label{p_11}
For every fixed \,$m>0$, nonlinear  ODE \eqref{CP_eq_f} has a two-parameter family of solutions \,$f_m(\Phi,B,D)$ that can be represented for large positive \,$\Phi$ in the principal approximation as
\begin{equation*}
f_m(\Phi,B,D)=B\,\Phi^{(m-1)/m)}\,\Bigl\{1+ \theta{\left(\Phi^{(m+1)/m},m,B\right)}+
\end{equation*}
\begin{equation}\label{fm_inf}
+D\,\Phi^{\varkappa_1}\,\exp{\left(-\frac{m}{B(m+1)}\Phi^{(m+1)/m}\right)}\Bigl[1+o(1)\Bigr]\Bigr\},
\quad \Phi\to\infty,
\end{equation}
where  \,$B$ and \,$D$ are parameters, \,$B>0$,
\begin{equation}\label{kapa1}
\varkappa_1=-(2m^2+4m-4)/[m(m+1)],
\end{equation}
and $\theta(x,m,B)$ is defined by Proposition \ref{p_9}.
\end{proposition}

\textbf{Case II.} In this case the exact solution \,$f=-\Phi^2 (m+1)/(6m)$, for \,$m\not=0$,
does not correspond to any global solution.

CP at infinity \eqref{Z_CP}, \eqref{Z_cond} now takes the form:
\begin{equation*}
\Phi^{2}\ddot{Z}+[2(m+4)/(m+1)]\,\Phi\dot{Z}+6\,Z-
\end{equation*}
\begin{equation}\label{Z_CP2}
-[6m/(m+1)] [(1+Z)^{-1}-1]\,\Phi \dot{Z}+ (1+Z)^{-1}\, (\Phi\dot{Z})^2=0, \quad \Phi\gg 1,
\end{equation}
\begin{equation}\label{Z_cond2}
\lim_{\Phi\to\infty}{Z(\Phi)}=\lim_{\Phi\to\infty}{\dot{Z}(\Phi)}=0.
\end{equation}
\subsubsection{The main result for \,${\bf m:\,1/2<m<\infty}$}
In this case, the solution of singular nonlinear CP \eqref{CP_eq_f}, \eqref{CP_cond_f}
exists globally and is positive on \,$(-a, \infty)$ (see Theorem \ref{t_4}).
Hence representation \eqref{fm_inf} holds for the solution for big \,$\Phi$
and some values of \,$B>0$ and \,$D$.

Further, let now \,$\Phi(\tau)$ be a solution of singular nonlinear IBVP \eqref{eq2_phi}--\eqref{eq2_0}, \eqref{con3f_phi},
for some fixed \,$b>0$ and \,$m:\,1/2<m<\infty$ (see Theorem \ref{t_3}). Then  the solution \,$f(\Phi)$ of singular CP \eqref{CP_eq_f},  \eqref{CP_cond_f} must satisfy condition \eqref{f_inf1}, i.e., in \eqref{fm_inf} we get
\begin{equation}\label{B_fin}
B=m b^{1/m}.
\end{equation}
\begin{theorem}\label{t_5}.
For fixed \,$m:\,1/2<m<\infty$, let the function \,$f_m(\Phi)=f_m(\Phi(\tau))$ be a solution of singular nonlinear CP \eqref{CP_eq_f}, \eqref{CP_cond_f}, with the trajectory  \,$\Phi(\tau)=\Phi(\tau,b)$
being the solution of singular IBVP \eqref{eq2_phi}--\eqref{eq2_0}, \eqref{con3f_phi}, for fixed \,$b>0$. Then:

(i) \,$f_m(\Phi)$ satisfies the constraints
\begin{equation}\label{f_rest1}
(a^2-\Phi^2)/2<f_m(\Phi)<a(\Phi+a), \quad -a\le\Phi\le a,
\end{equation}
\begin{equation}\label{f_rest2}
0<f_m(\Phi)<a(\Phi+a), \quad \Phi>a;
\end{equation}

(ii) for finite \,$\Phi$, the representation of the solution \,$f_m(\Phi)$ is given in Theorem \ref{t_4};

(iii) for big positive \,$\Phi$, representation \eqref{fm_inf} is true for the solution \,$f_m(\Phi)$
with \,$B=m b^{1/m}$ and some \,$D=D_m$.
\end{theorem}

Taking into account the relations
\begin{equation*}
x^{-2}\Phi^{(m-1)/m}=\Phi^{-3/m-1}, \quad x^{-3}\Phi^{(m-1)/m}=\Phi^{-4/m-2},
\end{equation*}
\begin{equation}\label{rel_x_Phi}
x^{-4}\Phi^{(m-1)/m}=\Phi^{-5/m-3}, \ldots,
\end{equation}
and formulas  \eqref{teta}--\eqref{teta_k-1}, from Theorem \ref{t_5}, we get, approximately,
\begin{equation*}
f_m(\Phi)=m b^{1/m}\,\Phi^{(m-1)/m}+\frac{m(m-1)(m-2)}{m+1}\,b^{2/m}\,\Phi^{-2/m}+
\end{equation*}
\begin{equation*}
+\frac{m(m-1)(m-2)}{2(m+1)^2}\,b^{3/m}\,\Phi^{-1-3/m}+O(\Phi^{-2-4/m})+
\end{equation*}
\begin{equation}\label{fPhi_inf}
+D_m\,\Phi^{\varkappa_1}\,\exp{\left(-\frac{b^{-1/m}}{m+1}\Phi^{(m+1)/m}\right)}\Bigl[1+o(1)\Bigr],
\quad \Phi\to\infty
\end{equation}
(cf. to \eqref{f_Phi_inf}).

Fig.B6 shows the graphs of the solutions of singular nonlinear CP \eqref{CP_eq_f}, \eqref{CP_cond_f}
in the assumptions of Theorem \ref{t_5}.
\begin{center}
\includegraphics[width=10cm,height=8cm]{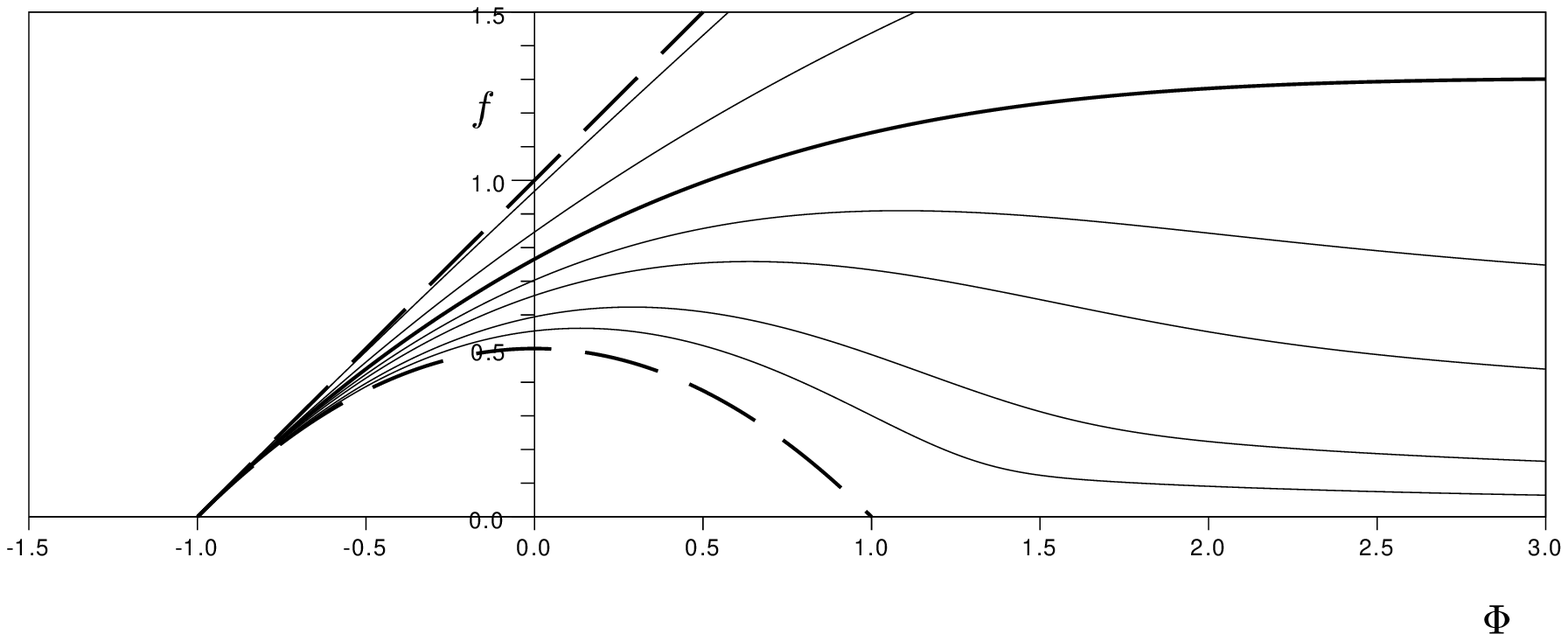}
\end{center}
\begin{center}
Fig.B6 ($a=1$; from bottom to top \,$m=:$\, $0.5$ (bold dotted line);
$0.55, \,0.6, \,0.7, \,0.8$; \,$1$ (bold line); \,$1.5,\,7$;
\,$\infty$ (bold dotted line))
\end{center}
\subsubsection{Notes to the case \,${\bf m:\,0<m<1/2}$}
In this case the solution of singular nonlinear CP \eqref{CP_eq_f}, \eqref{CP_cond_f}
may change the sign at some finite point \,$\Phi=\Phi_{zero}$. Then in a neighborhood of this
point we will have
\begin{equation*}
f\dot{f}=\dot{(f^2/2)}\sim {\rm const}, \quad f\sim \pm \sqrt{C_z|\Phi_{zero}-\Phi|},
\end{equation*}
where \,$c_z$  is a constant. Hence \,$\Phi=\Phi_{zero}$ is a branch point
and the solution becomes a multifunction (not uniquely defined function).
\begin{center}
\includegraphics[width=8cm,height=6cm]{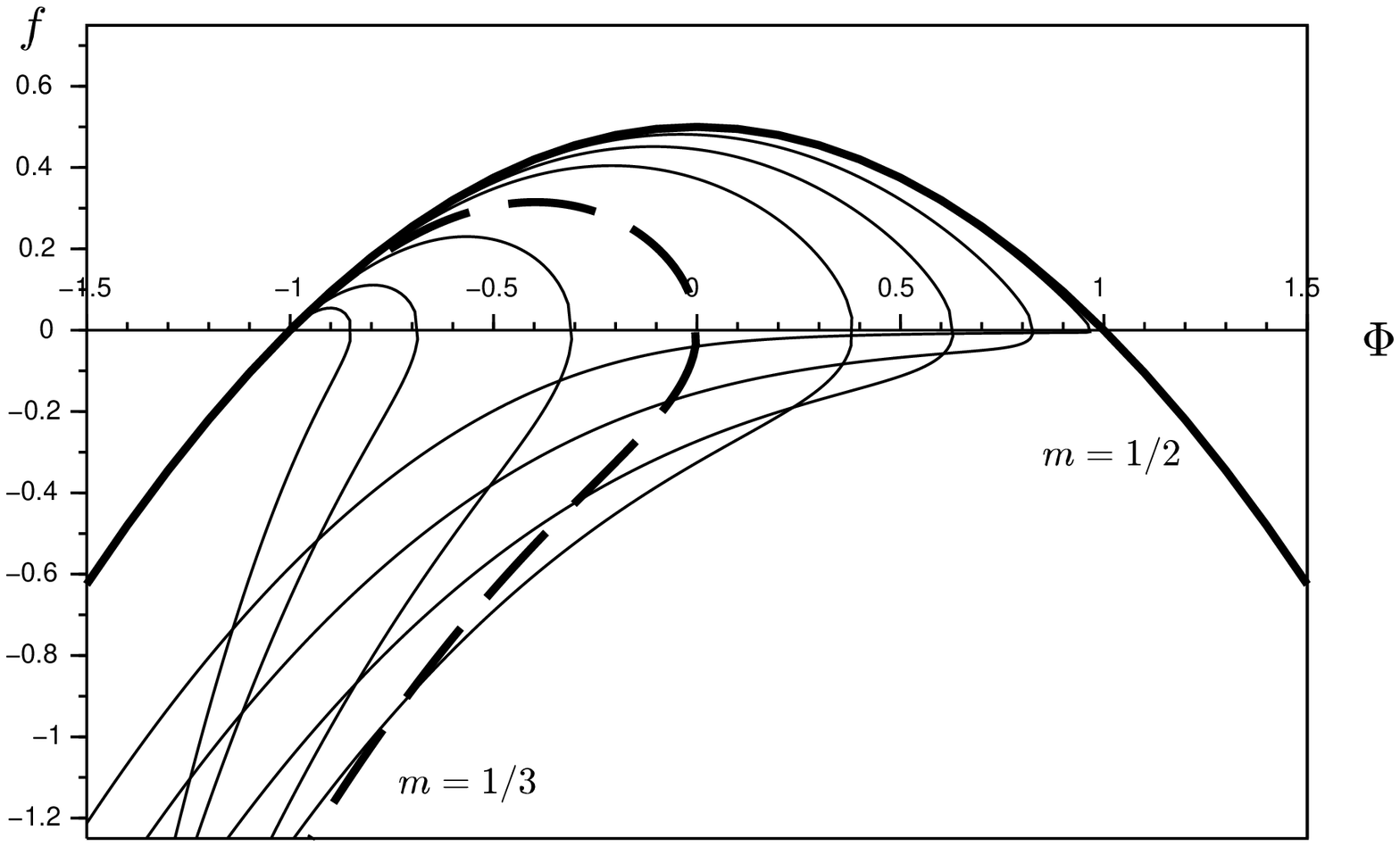}
\end{center}
\begin{center}
Fig.B7 ($a=1$; \,$m=:$ \,$1/2$, $255/512$, $31/64$, $11/24$, $5/12$, $1/3$,
$1/4$, $1/8$, $1/16$)
\end{center}
\section{Conclusions}
To conclude, we note that the analysis of singular problems for nonlinear ODEs arising in natural
science models is associated with great difficulties, which motivates special interest in
problems that are amenable to a comprehensive analysis. In our view, the approach to the hydrodynamic
problem described above applied in \cite{dks_07}, \cite{kss_09} and in this paper and a different one proposed in \cite{dies1_86},  \cite{dies2_86} supplement each other and can be of interest as applied to other problems.

\end{document}